\documentstyle[amssymb,10pt]{amsart}%\input amssym.def
\newtheorem{theorem}{Theorem}[section]
\newtheorem{lemma}[theorem]{Lemma}
\newtheorem{proposition}[theorem]{Proposition}
\newtheorem{corollary}[theorem]{Corollary} 
\theoremstyle{definition}  
\newtheorem{definition}[theorem]{Definition}

\newtheorem{example}[theorem]{Example}

\newtheorem{remark}[theorem]{Remark}

\newcommand{\id}{\operatorname{id}}

\newcommand{\End}{\operatorname{End}}

\newcommand{\Ad}{\operatorname{Ad}}

\newcommand{\nc}{\newcommand}
\nc{\Symm}{{\on{Sym}}}

\newcommand{\on}{\operatorname}   
\newcommand{\eps}{\varepsilon}
 \nc{\cE}{{\cal E}}

\renewcommand{\a}{{\mathfrak a}}
\renewcommand{\b}{{\mathfrak b}}

\renewcommand{\d}{{\mathfrak d}}
\nc{\SL}{{\mathfrak sl}}
\nc{\HH}{{\mathfrak h}}
\newcommand{\g}{{\mathfrak{g}}}
\newcommand{\x}{{\mathfrak{x}}}
\renewcommand{\u}{{\mathfrak{u}}}
\newcommand{\p}{{\mathfrak{p}}}
\newcommand{\vv}{{\mathfrak{v}}}

\newcommand{\z}{{\mathfrak{z}}}
\newcommand{\h}{{\mathfrak{h}}}
\renewcommand{\l}{{\mathfrak{l}}}
\renewcommand{\k}{{\mathfrak{k}}}
\newcommand{\n}{{\mathfrak{n}}}
\newcommand{\m}{{\mathfrak{m}}}
\renewcommand{\b}{{\mathfrak{b}}}

\newcommand{\SG}{{\mathfrak{S}}}

\nc{\wh}{\widehat}\nc{\wt}{\widetilde}

\newcommand{\la}{{\lambda}}

\newcommand{\ben}{\begin{enumerate}}
\newcommand{\een}{\end{enumerate}}
\newcommand{\ad}{{\text{ad}}}

\newcommand{\cO}{{\mathcal O}}

\newcommand{\CC}{{\mathbb{C}}}

\newcommand{\RR}{{\bold R}}
\newcommand{\LL}{{\bold L}}
\renewcommand{\tt}{{\bold t}}

\newcommand{\ZZ}{{\mathbb{Z}}}

\newcommand{\cP}{{\mathcal P}}

\newcommand{\cA}{{\mathcal A}}

\newcommand{\dd}{{\bf d}}

\newcommand{\alg}{{\operatorname{alg}}}
\newcommand{\Irr}{{\operatorname{Irr}}}

\renewcommand{\Ad}{{\operatorname{Ad}}}

\begin{document}

\title[Quantization of some Poisson-Lie dynamical $r$-matrices] 
{Quantization of some Poisson-Lie dynamical $r$-matrices
and Poisson homogeneous spaces}

\begin{abstract} 
Poisson-Lie (PL) dynamical $r$-matrices are generalizations of dynamical 
$r$-ma\-tri\-ces, where the base is a Poisson-Lie group. We prove analogues of 
basic results for these $r$-matrices, namely constructions of (quasi)Poisson
groupoids and of Poisson homogeneous spaces. We introduce a class of PL
dynamical $r$-matrices, associated to nondegenerate Lie bialgebras with a
splitting;
this is a generalization of trigonometric $r$-matrices with an abelian base.
We prove a composition theorem for PL dynamical $r$-matrices, and construct
quantizations of the polarized PL dynamical $r$-matrices. This way, we obtain 
quantizations of Poisson homogeneous structures on $G/L$ ($G$ a semisimple 
Lie group, $L$ a Levi subgroup), thereby generalizing earlier constructions. 
\end{abstract}

\dedicatory{To the memory of Joseph Donin}

\author{Benjamin Enriquez}
\address{IRMA (CNRS), rue Ren\'e Descartes, F-67084 Strasbourg, France}
\email{enriquez@@math.u-strasbg.fr}

\author{Pavel Etingof}
\address{Department of Mathematics, Massachusetts Institute of Technology,
Cambridge, MA 02139, USA}
\email{etingof@@math.mit.edu}

\author{Ian Marshall}
\address{D\'epartement de Math\'ematiques, EPFL, CH-1015 Lausanne, 
Switzerland}
\email{ian.marshall@@epfl.ch}

\maketitle

\section*{Introduction}

In this paper, we continue the study of Poisson-Lie (PL) dynamical
$r$-matrices, which was started in \cite{FM}. We construct new examples of such
$r$-matrices, together with their quantizations. We apply this to the
quantization of Poisson homogeneous spaces, which were introduced in 
\cite{DGS}.

In Section \ref{sect:1}, we define PL dynamical $r$-matrices and give some
examples. As in \cite{Lu}, PL dynamical $r$-matrices give 
rise to (quasi)Poisson groupoids. 

In Section \ref{sect:2}, we introduce the notion of a nondegenerate Lie
bialgebra with a splitting $\g = \l\oplus \u$. We associate to this datum 
a PL dynamical $r$-matrix $\sigma_\l^\g$. This construction is a PL analogue 
of \cite{EE2}
(see also \cite{FGP}, Proposition 1, and \cite{Xu}, Theorem 2.3). The main
example is the inclusion of a Levi subalgebra $\l$ in a simple Lie algebra
$\g$. When $\l$ coincides with the Cartan subalgebra $\h\subset \g$,
$\sigma_\l^\g$ is the standard trigonometric $r$-matrix (see \cite{EV1}). 
In general, $\sigma_\l^\g$ is an ingredient in a composition theorem for PL
dynamical $r$-matrices, generalizing \cite{EE2} (see also \cite{EV1}, Theorem
3.14 and \cite{FGP}, Proposition 1) and \cite{Mu} (who treated the case when 
$\l\subset \g$ is the inclusion of a Cartan in a Levi subalgebra). 

In Section 3, we introduce the notion of a polarized nondegenerate 
Lie bialgebra with a splitting. We construct a quantization $\Psi_\l^\g$ 
of $\sigma_\l^\g$ in this situation. This construction is a generalization 
of \cite{EE2}, and is based on a nonabelian analogue of the inversion of 
the Shapovalov pairing; from a representation theoretic viewpoint, this may 
be formulated in terms of intertwiners, as in \cite{EE2}. This idea is
already present in \cite{DM}; however, in order to carry out the analogue 
of the construction of \cite{EE2}, one needs to construct a left coideal 
$U_\hbar(\u_+) \subset U_\hbar(\g)$ (see Sections \ref{sect:3:1}, 
\ref{sect:3:2}). We also prove a quantum composition formula for the twists
$\Psi_\l^\g$ (Section \ref{sect:3:8}). 
  
The second part of the paper (Sections \ref{sect:4}, \ref{sect:5}, 
\ref{sect:6}) is devoted to applying these constructions to the quantization of
Poisson homogeneous spaces of the form $G/L$, where $\g$ is a simple Lie
algebra and $\l\subset \g$ is a Levi subalgebra.  

In Section \ref{sect:4}, we recall the classification of the Poisson homogeneous
structures on $G/L$ (\cite{DGS}). The set of all these Poisson structures is an
algebraic variety $\cP$. We introduce a Zariski open subvariety $\cP_0 \subset
\cP$, and show that the Poisson structures corresponding to elements of $\cP_0$
are exactly those which may be obtained either using the dynamical $r$-matrix 
$r_\l^\g$, or using its Poisson-Lie analogue $\sigma_\l^\g$. 

In Section \ref{sect:5}, we construct quantizations of all the Poisson
homogeneous structures corresponding to the elements of $\cP_0$. For this, we
prove an algebraicity result for $\Psi_\l^\g$. This result is based on the
computation of the Shapovalov pairing for $U_q(\g)$ (\cite{DCK}), and on the
quantum composition formula for the twists $\Psi_\l^\g$. 

In Section \ref{sect:6}, we compare the quantizations of $G/L$ obtained in 
Section \ref{sect:5} and in \cite{EE2}. In \cite{EE2}, there was constructed the
quantization of a family of Poisson structures, indexed by an analytic open
subset $U$ of $\cP_0$; this work is based on the Knizhnik-Zamolodchikov associator
and is therefore not purely algebraic, contrary to the construction of Section
\ref{sect:5}. We show that when the parameter of the 
Poisson structure belongs to $U$, both quantizations are equivalent. 

The Poisson homogeneous space $G/L$, equipped with a structure from $\cP_0$, 
may be viewed as a dressing orbit of $G^*$. 
Similarly, in the
rational case, the $\g$-invariant homogeneous structures on $G/L$ are the
coadjoint orbits $G/L \hookrightarrow \g^*$. In Section \ref{sect:7}, 
we quantize these $G$-  (or $\g$-) space embeddings (in the rational case,
quantizations were constructed in \cite{AL,EE2}). We also prove that 
the quantum function algebra of $G/L$ is a subalgebra of
$U_\hbar(\g)/\on{Ann}(M_\chi^\hbar)$ (in the rational case, of
$U(\g)/\on{Ann}(M_\chi)[[\hbar]]$), where $M_\chi^\hbar$ (resp., $M_\chi$)
is the generalized Verma module over $U_\hbar(\g)$ (resp., over $U(\g)$).  
Systems of generators for $\on{Ann}(M_\chi)$, whose classical limits
are systems of generators for the defining ideal of the orbit $G\chi \subset 
\g^*$, are constructed in \cite{O} when $\g$ is reductive of classical type.  

%Finally, we construct explicit generators of $\on{Ker}(\pi_\chi)$
%when $\g = \sl_n(\CC)$. 

\bigskip

\noindent
{\bf Acknowledgments.} 
Our motivation in this work came, to an important extent, from the papers
\cite{DGS} and \cite{DM}, as well as from very useful personal
explanations by Joseph Donin. We had
planned to show him the final results of our work and discuss what to do
next. Sadly, this was not to be. We dedicate this paper to his memory.

We would like to thank L. Feh\'er and A. Mudrov for discussions. 

The work of P.E. was partially supported by the NSF grant DMS-9988796.
P.E. and I.M. are indebted to IRMA (Strasbourg) for hospitality. 

\section{Poisson-Lie dynamical $r$-matrices and quasi-Poisson structures:
general facts} \label{sect:1}

In this section, we recall the notion of dressing actions of Poisson-Lie groups.
We define Poisson-Lie (PL) dynamical $r$-matrices and give some examples. 
We show that such $r$-matrices give rise to (quasi)Poisson groupoids, and 
discuss quantization of these constructions. The material of this section is a
PL generalization of basic constructions involving dynamical $r$-matrices. 

Our base field is $\CC$ (although the results of Sections \ref{sect:1}
to Section \ref{sect:3} hold over a field of characteristic $0$). 

\subsection{Poisson-Lie groups and Lie bialgebras}

A Poisson-Lie group $A$ is a Lie group equipped with a Poisson 
structure, such that the product map is a Poisson morphism 
$A\times A \to A$. This definition makes sense in the formal, 
algebraic or analytic categories, and we thus obtain the notions of 
formal, etc., Poisson-Lie groups. 

If $A$ is a (finite dimensional) Poisson-Lie group and $\a$ is its Lie algebra, 
then $\a$ identifies with $(\m/\m^2)^*$, where $\m$ is the 
maximal ideal of the local ring of $A$ at the origin, and the Poisson bracket 
$\wedge^2(\m) \to \m$ induces a map $\wedge^2(\m/\m^2)\to \m/\m^2$, 
whose dual is a map $\delta : \a\to \wedge^2(\a)$. Then $(\a,\delta)$
is a Lie bialgebra, i.e., $(\delta\otimes \id)\circ \delta(a) + $ 
cyclic permutations = 0
for any $a\in\a$, and $\delta([a,b]) = [a \otimes 1 + 1 \otimes a, \delta(b)]
+ [\delta(a),b\otimes 1 + 1 \otimes b]$ for any $a,b\in\a$. 

To a Lie bialgebra $\a$, one associates a formal Poisson-Lie group, 
i.e., a formal series Hopf algebra $\cO(A)$ with a Poisson structure
compatible with the coproduct: $\cO(A) = U(\a)^*$.

\subsection{Dressing actions}

\begin{definition} \label{def1}
{\it Let $(\a,\delta)$ be a Lie bialgebra 
and let $\cO$ be a Poisson algebra. 
A Poisson action of $\a$ on $\cO$ is a 
linear map $\a\to \on{Der}(\cO)$, 
$a\mapsto \theta_a$, such that 
$\theta_a(\{f,g\}) = \{\theta_a(f),g\} + \{a,\theta_a(g)\} 
+ \sum \theta_{a^{(1)}}(f) \theta_{a^{(2)}}(g)$, where 
$f,g\in \cO$, $a\in \a$, and $\delta(a) = \sum a^{(1)} \otimes a^{(2)}$.}
\end{definition}

Let $A$ be a Poisson-Lie group with Lie bialgebra $(\a,\delta)$, 
let $M$ be a Poisson manifold, and let $A\times M\to M$ be a left 
Poisson action of
$A$ on $M$. Then the map $\LL : \a\to \on{Der}(\cO_M)$ defined by 
$(\LL_a f)(x) := {{\on{d}}\over{\on{d}\eps}}_{|\eps = 0}
f(e^{-\eps a}x)$ is a Poisson action of $(\a,\delta)$ on 
$\cO_M$ (here $\cO_M$ is the function algebra of $M$). 
In the same way, if $M\times A \to A$ is a right Poisson action of $A$ on $M$, 
then $\RR : \a\to\on{Der}(\cO_M)$, defined by 
$(\RR_a f)(x) := {{\on{d}}\over{\on{d}\eps}}_{|\eps = 0}
f(x e^{\eps a})$ is a Poisson action of $(\a,\delta)$ on 
$\cO_M$. We also use the notation $(\LL_{a(x)}f)(x)$, $(\RR_{a(x)}f)(x)$
to denote left and right infinitesimal translations by an element of 
$\cO_M \otimes \a$. 

\medskip 

Assume that $\a$ is finite dimensional. 
Let $\d(\a) = \a \oplus \a^*$ be the double Lie bialgebra of $(\a,\delta)$. 
The injections $\a\to \d(\a)$, $\a^{*,\on{cop}} \to \d(\a)$
are Lie bialgebra morphisms (here cop means the Lie bialgebra with opposite
cobracket). Let $A$ be a Poisson-Lie group with Lie bialgebra $\a$. 
Assume that the adjoint action of $\a$ on $\d(\a)$ extends to $A$ 
(this is the case e.g. if $A$ is the formal group of $\a$). The {\it
left and right dressing actions} are Poisson actions 
of $\d(\a)$ on $\cO_A$, defined by 
$$
(\on{dress}^L_a f)(x) = -\RR_{(\Ad(x^{-1})(a))_\a}f(x), 
\quad 
(\on{dress}^R_a f)(x) = -\LL_{(\Ad(x)(a))_\a}f(x), 
$$
where $a\in \d(\a)$, $f\in \cO_A$, and $\alpha\mapsto \alpha_\a$ 
is the projection on
$\a$ parallel to $\a^*$. If $a\in \a$, we have $\on{dress}^L_a = \LL_a$, 
$\on{dress}^R_a = \RR_a$. 

If $f,g\in\cO_A$, we have 
\begin{equation} \label{useful}
\{f,g\} = \sum_i \on{dress}_{\eps^i}^L(f) \LL_{e_i}(g)
= - \sum_i \on{dress}^R_{\eps^i}(f) 
\RR_{e_i}(g),
\end{equation} 
where $(e_i),(\eps^i)$ are dual bases of $\a$ and $\a^*$. 

Assume that $A\to D$, $A^*\to D$ are Poisson-Lie group morphisms 
associated with the Lie bialgebra morphisms $\a\to \d(\a)$, 
$\a^{*,\on{cop}} \to \d(\a)$, and that the maps $A \to D / A^*$, 
$A \to A^* \setminus D$
are Poisson isomorphisms (these assumptions hold e.g. 
if $A,A^*,D$ are the formal groups of $\a,\a^*,\d(\a)$). 
Then $\on{dress}^L$ (resp., $\on{dress}^R$) corresponds to the left 
action of $D$ on $D/ A^*$ (resp., right action of $D$ on 
$A^* \setminus D$), which are Poisson. 

If $\b$ is a Lie algebra and $\a$ is the Lie bialgebra $\b^*$
(with zero bracket), then the simply-connected Lie group 
with Lie algebra $\a$ is $A = \a$. Then  $\on{dress}^L = \on{dress}^R$
and the common restriction of both actions to $\b\subset \d(\a)$
is the adjoint action of $\b$ on $S^\cdot(\b) 
= \cO_{\b^*} = \cO_A$.

\subsection{Poisson-Lie dynamical $r$-matrices}

Let $(\g,\delta)$ be a finite dimensional Lie bialgebra and $\h\subset \g$
be a Lie subbialgebra. Let $\h^*$ be the Lie bialgebra dual to $\h$ and
let $H^*$ be a Poisson-Lie group with Lie bialgebra 
$\h^*$. We denote by 
$\cO_{H^*}$ the function ring of $H^*$ (when $H^*$ is a formal 
group, $\cO_{H^*} = U(\h^*)^*$ will be denoted $\CC[[H^*]]$; 
if $H^*$ is an algebraic group, 
then $\cO_{H^*}$ is a Hopf algebra contained in $U(\h^*)^*$ and will be denoted 
$\CC[H^*]$). 

Define $\on{d}^L : \cO_{H^*} \to \h \otimes \cO_{H^*}$  
by $\on{d}^L(f) = - \sum_i e_i \otimes \LL_{\eps^i}(f)$, where 
$(e_i), (\eps^i)$ are dual bases of $\h$ and $\h^*$. 
If $\rho = \sum_\alpha a_\alpha \otimes b_\alpha \otimes \ell_\alpha$, 
we set $\on{d}^L\rho = \sum_\alpha a_\alpha \otimes b_\alpha
\otimes \on{d}^L(\ell_\alpha)$. 

If $V$ is a vector space, $\wedge^i(V)$ is the space of totally 
antisymmetric tensors in $V^{\otimes i}$. If $x\in V^{\otimes 3}$
is antisymmetric w.r.t. a pair of tensor factors, then 
$\on{Alt}(x)$ is the sum of its cyclic permutations 
$x + x^{2,3,1} + x^{3,1,2}$. 
We set $\on{CYB}(\rho) = [\rho^{1,2},\rho^{1,3}] + [\rho^{1,2},\rho^{2,3}]
+ [\rho^{1,3},\rho^{2,3}]$. 

The dressing actions of $\h$ on $\cO_{H^*}$ extend to actions on 
the localizations $\cO_{H^*}[1/{\bf P}]$, where ${\bf P}\in \cO_{H^*}$ 
is nonzero. 

\begin{definition} \label{def:PL} (see \cite{DM})
{\it Let $Z\in \wedge^3(\g)^\g$ be such that $\on{Alt}(\delta \otimes
\id^{\otimes 2})(Z) =0$. 
A Poisson-Lie dynamical $r$-matrix for $(H^*,\g,Z)$
is an element $\rho\in \wedge^2(\g) \otimes \cO_{H^*}[1/{\bf P}]$, 
which we view as a function $\rho : H^* \to \wedge^2(\g)$, 
such that: 

(a) for any $a\in \h$, $\on{dress}^L_a(\rho) + [a\otimes 1 + 1
\otimes a,\rho] = 0$ 

(b) $\on{CYB}(\rho) + \on{Alt}(\on{d}^L \rho) 
+ \on{Alt}((\delta\otimes \id)(\rho)) = Z$.}
\end{definition}

This definition is the specialization of the notion from \cite{DM}
of a classical dynamical $r$-matrix over a base algebra, when the base is
$H^*$. 

\begin{remark}
If $\rho$ is as in Definition \ref{def:PL}, then $h\mapsto \rho(h^{-1})$
satisfies the analogue of Definition \ref{def:PL}, where $\on{dress}^L$, 
$\on{d}^L$ are replaced by $\on{dress}^R$, $\on{d}^R$ ($\on{d}^R$
is defined by replacing $\LL$ with $\RR$ in the definition of $\on{d}^L$).  
\end{remark} 

\begin{remark}
When $\h = \g$, we may identify $\rho$ with a right-invariant 
2-form $\omega$ on $G^*$. Then the term $\on{Alt}(\on{d}^L \rho) + 
\on{Alt}((\delta\otimes \id)(\rho))$ identifies with its exterior 
derivative $\on{d}\omega$. 
\end{remark}

\begin{example} \label{from:old}
If $\delta = 0$, then we may take $H^*$ to be 
the formal group of $\h^*$. Then $\on{dress}$ is the adjoint action of 
$\h$ on $\cO_{\h^*} = \wh S^\cdot(\h)$. In that case, $\rho \mapsto -\rho$ 
yields a bijection between Poisson-Lie (PL)
dynamical $r$-matrices for $(H^*,\g,Z)$ and (usual) 
dynamical $r$-matrices for $(\h,\g,Z)$. 

More generally, assume that $\delta_{ | \h} = 0$, and that 
$(\g,\delta)$ is coboundary, i.e., we have $\delta(a) = [a\otimes 1 
+ 1 \otimes a,r_0]$, where $r_0\in \wedge^2(\g)$ is such that 
$Z(r_0) := \on{CYB}(r_0) 
\in \wedge^3(\g)^\g$. Then the map $\rho \mapsto r_0 - \rho$
is a bijection between $\{$PL dynamical $r$-matrices 
for $(H^*,\g,Z)\}$ and $\{$usual dynamical $r$-matrices for 
$(\h,\g,Z+Z(r_0))\}$. 

In particular, if $\g$ is a semisimple Lie algebra and  
$\h\subset\g$ is a Cartan subalgebra, then  
$r_\h^\g(\lambda) := - {1\over 2}\sum_{\alpha\in\Delta_+}
(e_\alpha\wedge f_\alpha) \on{coth}{{(\lambda,\alpha)}\over 2}$
is a usual dynamical $r$-matrix for $(\h,\g,Z(r_0))$, 
where $r_0 = {1\over 2} \sum_{\alpha\in\Delta_+} e_\alpha \wedge f_\alpha$.
Therefore $r_0-r(\lambda)$ is a 
PL dynamical $r$-matrix for $(H^*,\g,0)$, where 
$\g$ is equipped with its standard Lie bialgebra structure 
$\delta(a) = [a\otimes 1 + 1 \otimes a,r_0]$ 
(we set $x\wedge y := x\otimes y - y \otimes x$.) 
\end{example}

\begin{example} \label{ex:FM}
(The Balog-Feh\'er-Palla $r$-matrix \cite{BFP,FM}.) 
Let $(\g,r)$ be a factorizable Lie bialgebra. 
Recall that this means that $t := r + r^{2,1}$ is a nondegenerate
element of $S^2(\g)^\g$, $\on{CYB}(r) = 0$. The cobracket of 
$\g$ is given by $\delta(a) = 
[a \otimes 1 + 1 \otimes a,r]$ for any $a\in\g$. 
Set $Z := [t^{1,2},t^{2,3}]$. The relation 
$\on{Alt}(\delta\otimes \id^{\otimes 2})
(Z) =0$ is satisfied because the twist by $r-r^{2,1}$ of the trivial 
quasi-Lie bialgebra structure on $\g$ is a quasi-Lie bialgebra.    

Set $L(\xi) = (\id\otimes\xi)(r)$ and $R(\xi) = - (\xi\otimes
\id)(r)$. Then the maps $R,L : \g^*\to\g$ are Lie algebra morphisms, 
and $R\oplus L : \g^* \to \g\oplus \g$ is a Lie algebra injection. 
We also denote by $R,L$ the corresponding formal group morphisms $G^*\to G$,   
and by $\lambda : G^* \to G$ the map such that $\lambda(g^*) = 
L(g^*) R(g^*)^{-1}$. Finally, we define 
$g^* \mapsto g(g^*)$ to be the map $G^* \to \End(\g)$, such that 
$g(g^*) = \Ad(\lambda(g^*))$. 

Define $\rho_{\on{BFP}} : G^* \to \wedge^2(\g)$ by 
$$
\rho_{\on{BFP}}(g^*) =
\big( \big( \nu {\id + {g(g^*)^{2\nu}}\over{\id - g(g^*)^{2\nu}}} 
- {1\over 2}{{\id+g(g^*)}\over{\id-g(g^*)}}\big) \otimes \id\big)(t).  
$$ 
Then $\rho_{\on{BFP}}$ is a dynamical $r$-matrix for $(G^*,\g,Z_\nu)$, 
where $Z_\nu = (\nu^2 - {1\over 4}) [t^{1,2},t^{2,3}]$, and $\nu\in\CC$
is an arbitrary constant. 
\end{example}

\begin{remark}
A constant PL dynamical 
$r$-matrix for $(G^*,\g,Z)$ is the same as a classical twist 
between the quasi-Lie bialgebras $(\g,-\delta,0)$ and $(\g,-\delta - \delta_\rho,
-Z)$, where $\delta_\rho(x) = [x\otimes 1 + 1 \otimes x,\rho]$
(see \cite{Dr:QH}). \hfill \qed \medskip 
\end{remark}

In Section \ref{sect:sigma}, we will construct new families of examples of 
PL dynamical $r$-matrices.   

\subsection{Quasi-Poisson structures} \label{quasi:poisson}

If $\cO$ is a commutative algebra equipped with an antisymmetric
bracket, satisfying the Leibniz identity w.r.t. each variable, 
a Poisson action of a Lie bialgebra $(\a,\delta)$ on $\cO$
is defined as a Lie algebra morphism $\a\to \on{Der}(\cO)$, 
satisfying the condition of Definition \ref{def1}.  

\begin{proposition} \label{prop:QP}
Let $(\g,\delta)$ be a finite dimensional Lie bialgebra, 
let $Z\in \wedge^3(\g)^\g$, and let $\h\subset\g$ be a Lie subbialgebra. 
Let $\rho = \sum_{\alpha} a_\alpha \otimes b_\alpha \otimes \ell_\alpha
\in\wedge^2(\g)\otimes \cO_{H^*}$. There is a unique bracket on 
$\cO_{H^*} \otimes \cO_G$, antisymmetric and satisfying the Leibniz identity 
w.r.t. each variable, such that 
$$
\forall f_1,f_2\in \cO_{H^*}, \quad
\{f_1\otimes 1,f_2\otimes 1\} = 
\{f_1,f_2\}_{H^*} \otimes 1 ,   
$$
$$
\forall g\in \cO_{G}, \forall f\in \cO_{H^*}, \quad 
\{f\otimes 1,1\otimes g\} = - \sum_i \LL_{\eps^i}(f)\otimes \LL_{e_i}(g), 
$$
$$
\forall g_1,g_2\in\cO_{G}, \quad 
\{1 \otimes g_1, 1 \otimes g_2\} = 1 \otimes \{g_1,g_2\}_{G} 
+ \sum_\alpha \ell_\alpha \otimes \LL_{a_\alpha}(g_1) \LL_{b_\alpha}(g_2) . 
$$
Here $(e^i),(\eps_i)$ are dual bases of $\h,\h^*$, and $\{-,-\}_A$
is the Poisson bracket of a Poisson-Lie group $A$. 

If $\rho$ is a PL dynamical $r$-matrix for $(H^*,\g,Z)$, then 
this bracket satisfies $\{\{\phi_1,\phi_2\},\phi_3\} + \on{c.p.}
= m \big( 
(\id\otimes\RR)^{\otimes 3}(Z)(\phi_1\otimes \phi_2 \otimes \phi_3)
\big)$ for any $\phi_1,\phi_2,\phi_3\in\cO_{H^*} \otimes \cO_G$.  

This structure is equipped with the following 
commuting Poisson actions: 

(a) the action of $(\h,\delta)$ by $\on{dress}^L \otimes \id + 
\id\otimes \LL$, 

(b) the action of $(\g,\delta)$ by $\id\otimes \RR$. 
\end{proposition}

\medskip 

{\em Proof.} Straightforward. \hfill \qed \medskip 

The notion of a quasi-Poisson structure (\cite{AKM}) may be generalized 
as follows: $(\a,\delta_\a,Z_\a)$ is a quasi-Lie bialgebra, $M$ is a manifold 
equipped with a bivector and an action of $\a$ by vector fields, 
such that $\theta_a(\{f,g\}) = \{\theta_a(f),g\} + \{f,\theta_a(g)\} 
+ \sum \theta_{a^{(1)}}(f) \theta_{a^{(2)}}(g)$ and $\{\{f,g\},h\} + \on{c.p.}
= m(\theta^{\otimes 3}(Z_\a)(f \otimes g \otimes h))$. 

In the situation of Definition \ref{def:PL},   
$(\g,\delta,Z)$ is a quasi-Lie bialgebra and 
the structure defined in Proposition 
\ref{prop:QP} is quasi-Poisson under the action of $(\g,\delta,Z)$
by $\id\otimes \RR$.  

\subsection{Poisson groupoids} One checks that $H^*\times G \times H^*$
is equipped with a Poisson groupoid structure, corresponding to the Poisson
bracket on $\cO_{H^*} \otimes \cO_{G} \otimes \cO_{H^*}$
$$
\{f^{(1)},f^{\prime (1)}\} = \{f,f'\}_{H^*}^{(1)}, \; 
\{f^{(3)},f^{\prime (3)}\} = -\{f,f'\}_{H^*}^{(3)}, \; 
\{f^{(1)},f^{\prime (2)}\} = 0, 
$$
$$
\{f^{(1)},g^{(2)}\} = - \sum_i \LL_{\eps^i}(f)^{(1)}
\LL_{e_i}(g)^{(2)}, \; 
\{f^{(3)},g^{(2)}\} = - \sum_i \RR_{\eps^i}(f)^{(3)}
\RR_{e_i}(g)^{(2)}, 
$$
$$
\{g^{(2)},g^{\prime(2)}\} = \{g,g'\}_G^{(2)} 
+ \sum_\alpha
\ell_\alpha^{(1)} \big( \LL_{a_\alpha}(g) \LL_{b_\alpha}(g') \big)^{(2)} 
- \sum_\alpha
\ell_\alpha^{(3)} \big( \RR_{a_\alpha}(g) \RR_{b_\alpha}(g') \big)^{(2)} 
$$
for any $f,f'\in \cO_{H^*}$ and $g,g'\in \cO_G$.  
This construction is carried out in \cite{FM} in the setup of 
Example \ref{ex:FM}.  

\subsection{Quantization of PL dynamical $r$-matrices}

Let $(U_\hbar(\g),\Delta)$ be a quantized universal enveloping (QUE) algebra. 
We say that $\Phi\in U_\hbar(\g)^{\wh\otimes 3}$ is an associator for 
$(U_\hbar(\g),\Delta)$ if 
\begin{equation} \label{conds:Phi:1}
\Phi^{2,3,4} \Phi^{1,23,4} \Phi^{1,2,3} = \Phi^{1,2,34} \Phi^{12,3,4}, 
\end{equation}
\begin{equation} \label{conds:Phi:2}
\forall x\in U_\hbar(\g), \quad [\Delta^{(3)}(x),\Phi] = 0, 
\quad 
\Phi = 1 + O(\hbar), \quad \on{alt}(\Phi) = O(\hbar^2). 
\end{equation}
Here alt$: V^{\otimes 3} \to V^{\otimes 3}$ is defined by 
alt$(x) = \sum_{\sigma\in \SG_3} \on{sgn}(\sigma) \sigma(x)$, and 
$\wh\otimes$ is the $\hbar$-adically completed tensor product. 

Then $(U_\hbar(\g),\Delta,\Phi)$ is a quasi-Hopf algebra, therefore it gives
rise to a quasi-Lie bialgebra $(\g,\delta,Z)$ (see \cite{Dr:QH}), where
$Z = ( \hbar^{-2}\on{alt}(\Phi)$ mod $\hbar)$, and 
\begin{equation} \label{conds}
Z\in \wedge^3(\g)^\g, \; \on{and} \; 
\on{Alt}(\delta\otimes \id\otimes \id)(Z) = 0. 
\end{equation} 

Conversely, if $Z$ satisfies (\ref{conds}), then 
$\Phi\in U_\hbar(\g)^{\wh\otimes 3}$ satisfying (\ref{conds:Phi:1}),
(\ref{conds:Phi:2}) will be called a quantization of $Z$. 

\medskip 

Let $(U_\hbar(\g),\Delta)$ be a QUE algebra, $\Phi\in 
U_\hbar(\g)^{\wh\otimes 3}$ be an associator, let $U_\hbar(\h)\subset
U_\hbar(\g)$ be a QUE subalgebra and let $\Psi \in 
U_\hbar(\g)^{\wh\otimes 2}\wh\otimes U_\hbar(\h)$ be such that 
\begin{equation} \label{conds:Psi:1}
\Psi^{2,3,4}\Psi^{1,23,4} \Phi^{1,2,3} = \Psi^{1,2,34} \Psi^{12,3,4}, 
\end{equation}
\begin{equation} \label{conds:Psi:2}
\hbar^{-1}(\Psi-1) \in U_\hbar(\g)^{\wh\otimes 2}
\wh\otimes \CC[[H^*]]_\hbar, \quad 
\forall y\in U_\hbar(\h), [\Psi,\Delta^{(3)}(y)] = 0. 
\end{equation}
Here $\CC[[H^*]]_\hbar \subset U_\hbar(\h)$ is defined as 
$\{a\in U_\hbar(\h) | 
\forall n\geq 0, (\id - \eta\circ\eps)^{\otimes n} \circ \Delta^{(n)}(f)
= O(\hbar^n)\}$; it is a flat deformation of $\CC[[H^*]] := U(\h^*)^*$ 
(see \cite{Dr:QG,Gav}). 

Set 
\begin{equation} \label{rho:Psi}
\rho := \big( \hbar^{-1}(\Psi - \Psi^{2,1,3})\ \on{mod}\  \hbar \big). 
\end{equation} 
Then $\rho \in \wedge^2(\g) \otimes \CC[[H^*]]$. 

\begin{proposition}
Let $(\g,\delta)$ be the classical limit of $(U_\hbar(\g),\Delta)$
(so for $x\in\g$, $\delta(x) = \big(
\hbar^{-1} ( \Delta(\wt x) - \Delta(\wt x)^{2,1}) 
\on{mod}\ \hbar \big)$ for any $\wt x$ such that 
$(\wt x\ \on{mod}\ \hbar) = x$). 
Then $\rho$ is a PL dynamical $r$-matrix for $(H^*,\g,Z)$. 
\end{proposition}

{\em Proof.} The proof relies on two facts: 

(a) if $f\in \CC[[H^*]]_\hbar$ and if we view $\Delta(f)$  as an element of 
$U_\hbar(\g) \wh\otimes \CC[[H^*]]_\hbar$, then 
$\big( \hbar^{-1}(\Delta(f) - 1 \otimes f)\  \on{mod}\ \hbar\big) 
= \on{d}^L(f_0)$, where $f_0 = (f\ \on{mod}\ \hbar)\in \CC[[H^*]]$;  

(b) if $a\in \CC[[H^*]]_\hbar$, $x,y\in U_\hbar(\g)$, 
$f\in \CC[[H^*]]_\hbar$, then 
\begin{align*}
& \big( \hbar^{-1}[\Delta^{(3)}(a), x\otimes y \otimes f]\ \on{mod} \ \hbar\big)
\\ & 
= \sum_i \big(
[e_i,x_0] \otimes y_0 \otimes f_0 + x_0 \otimes [e_i,y_0] \otimes f_0 
+ x_0 \otimes y_0 \otimes \on{dress}^L_{e_i}(f_0)  \big) 
(1\otimes 1 \otimes \LL_{\eps^i}(a_0)) , 
\end{align*}
where $x_0,\ldots,a_0$ are the reductions modulo $\hbar$
of $x,\ldots,a$.  
 
(a) follows from \cite{EGH}, appendix. (b) follows from (a) and 
(\ref{useful}). 
\hfill \qed \medskip 

An element $\Psi \in U_\hbar(\g)^{\wh\otimes 2} \wh\otimes U_\hbar(\h)$
satisfying (\ref{conds:Psi:1}), (\ref{conds:Psi:2}) and (\ref{rho:Psi})  
will be called a quantization of $\rho$ compatible with $\Phi$. 

\subsection{Relation with twists of QUE algebras} \label{ex:pseudotwists}

Assume that $(\g,\delta)$ is coboundary, and $\h\subset\g$ is such that
$\delta_{|\h} = 0$. Let $r_0\in \wedge^2(\g)$ be such that $\delta(a)
= [a\otimes 1 + 1 \otimes a,r_0]$ and $Z(r_0) \in \wedge^3(\g)^\g$. 
Let $r'$ be a dynamical $r$-matrix for $(\h,\g,0)$; then $r_0 - r'$
is a PL dynamical $r$-matrix for $(H^* = \h^*,\g,-Z(r_0))$
(see Example \ref{from:old}). This correspondence may be quantized as 
follows. 

Assume that $\Phi \in U(\g)^{\otimes 3}[[\hbar]]$ is an 
associator quantizing $Z(r_0)$ and let $J_0 \in U(\g)^{\otimes 2}[[\hbar]]$
be a $\h$-invariant quantization of $r_0$ such that 
$$
J_0^{1,23} J_0^{2,3} = \Phi J_0^{12,3} J_0^{1,2}. 
$$
Let $J'\in U(\g)^{\otimes 2} \otimes U(\l)[[\hbar]]$ be a pseudotwist 
quantizing $r'$, i.e., 
$$
J^{\prime 12,3,4} J^{\prime 1,2,34} = (\Phi^{-1})^{1,2,3} J^{\prime 1,23,4} 
J^{\prime 2,3,4}. 
$$
Define $U_\hbar(\g)$ as the QUE algebra $(U(\g)[[\hbar]],\on{Ad}(J_0^{-1}) 
\circ \Delta_0)$. Then 
\begin{equation} \label{rel:Psi:J'}
\Psi := (J')^{-1} J_0^{1,2}
\end{equation} 
is a quantization of $r_0 - r'$, associated with the quantum group 
$U_\hbar(\g)$. 

In the particular case where $\g$ is semisimple, $\h\subset \g$ is the 
Cartan subalgebra, and $r'$ is as in Example \ref{from:old}, 
a quantization $J'$ of $r'$ was obtained in \cite{EE2}.
On the other hand, a quantization $\Psi$ of the PL dynamical $r$-matrix 
$r_0 - r'$ was obtained in \cite{EV2} (we are going to generalize this 
construction in Section \ref{sect:quant}). In \cite{EE2}, Remark 3.7,
it was also conjectured that the relation between $\Psi$ and $J'$ is given by 
(\ref{rel:Psi:J'}).

\subsection{Quantization of quasi-Poisson structures}

If $M$ is a (formal or algebraic) manifold with a quasi-Poisson structure under the action of a 
quasi-Lie bialgebra $(\a,\delta_\a,Z_\a)$, then a quantization of $M$ is
defined as the data of (a) a quantization $(U_\hbar(\a),\Delta_\a,\Phi_\a)$
of $(\a,\delta_\a,Z_\a)$, (b) a $\CC[[\hbar]]$-module $\cO_{M,\hbar}$, 
isomorphic to $\cO_M[[\hbar]]$,  
equipped with a bilinear map $\mu : \cO_{M,\hbar}^{\wh\otimes 2}
\to \cO_{M,\hbar}$, deforming the product of $\cO_M$, (c) an action 
$\Theta_\a$ of $U_\hbar(\a)$ on $\cO_{M,\hbar}$, deforming the action of 
$\a$ on $\cO_M$, such that 
$\Theta_\a(a) (\mu(f\otimes g)) = \sum \Theta_\a(a^{(1)})(f) 
\Theta_\a(a^{(2)})(g)$, 
and $\mu \circ (\mu\otimes \id) = \mu \circ (\id\otimes \mu) \circ 
\Theta_\a^{\otimes 3}(\Phi_\a)$.  

\medskip 

Assume that $(\g,\delta)$ is a Lie bialgebra, $Z\in \wedge^3(\g)^\g$
is such that $\on{Alt}(\delta\otimes \id\otimes \id)(Z) = 0$, 
$\h\subset\g$ is a Lie subbialgebra, $\rho$ is a 
PL dynamical $r$-matrix for $(H^*,\g,Z)$. 

Assume that $(U_\hbar(\g),\Delta)$ is a quantization of $(\g,\delta)$, 
$\Phi$ is a quantization of $Z$, $U_\hbar(\h)\subset U_\hbar(\g)$
is a quantization of $\h\subset \g$, and $\Psi = \sum_\alpha 
A_\alpha\otimes B_\alpha \otimes L_\alpha$ is a quantization of 
$\rho$. Then $(U_\hbar(\g),\Delta,S^{\otimes 3}(\Phi)^{-1})$
is also a quantization of $(\g,\delta,Z)$ (here $S$ is the antipode of 
$U_\hbar(\g)$). 

Set $\CC[[G]]_\hbar := U_\hbar(\g)^*$; this is an algebra, equipped with 
commuting actions $\LL,\RR$ of $U_\hbar(\g)$, where 
$(\LL(a)f)(x) = f(xa)$, $(\RR(a)(f))(x) = f(S^{-1}(a)x)$. 

\begin{proposition}
Let $A := U_\hbar(\h) \wh\otimes \CC[[G]]_\hbar$, and define 
$\mu : A^{\wh\otimes 2} \to A$ by 
$$
\mu \big( (a\otimes f) \otimes (b\otimes g)\big) 
= \sum_\alpha L_\alpha a^{(2)} b \otimes
\LL_{A_\alpha}(f) \LL_{B_\alpha a^{(1)}}(g) , 
$$
$\RR : U_\hbar(\g) \otimes A \to A$ by $\RR(x)(a\otimes f) = a \otimes 
\RR(x)(f)$. Then $B := \CC[[H^*]]_\hbar \bar\otimes \CC[[G]]_\hbar
\subset A$ is stable under $\mu$ and $\Theta$, and is a quantization of the 
quasi-Poisson structure of Proposition \ref{prop:QP}, compatible with the 
quantization $(U_\hbar(\g),\Delta,S^{\otimes 3}(\Phi)^{-1})$
of $(\g,\delta,Z)$. 
\end{proposition}

Here $\CC[[a_1,\ldots,a_d,\hbar]] \bar\otimes \CC[[b_1,\ldots,b_{d'},\hbar]]$
means $\CC[[a_1,\ldots,b_{d'},\hbar]]$. 

$B$ is also equipped with an action of $U_\hbar(\h)$
by $\Theta'(x)(a\otimes f) = x^{(2)} a S(x^{(3)}) \otimes \LL_{x^{(1)}}(f)$, 
and the actions $\Theta,\Theta'$ commute.

\section{A composition theorem for Poisson-Lie dynamical $r$-matrices}
\label{sect:2}

In this section, we introduce a PL dynamical $r$-matrix $\sigma_\l^\g$,
associated to a nondegenerate Lie bialgebra with a splitting. This construction
is a generalization of \cite{EE2} (see also \cite{FGP}, Proposition 1,
\cite{Xu}, Theorem 2.3). When $\l$ is a Cartan subalgebra of a semisimple Lie
algebra $\g$, then $\sigma_\l^\g$ is a trigonometric $r$-matrix. In general, 
$\sigma_\l^\g$ is the basic ingredient of a composition theorem for 
PL dynamical $r$-matrices (Section \ref{sect:2:4}), 
generalizing \cite{EE2}, Proposition 0.1 (also \cite{FGP}, Proposition 1).  

\subsection{Lie bialgebras with a splitting} \label{LBA:split}
\label{nondeg}

Let $(\g,\delta)$ be a finite dimensional 
Lie bialgebra and let $\l\subset \g$ be a Lie
subbialgebra. Assume that we have a decomposition $\g = \l\oplus \u$, 
such that $[\l,\u] \subset \u$ and $\delta(\u) \subset \big( (\l\otimes\u) 
\oplus (\u\otimes \l) \big) \oplus (\u\otimes \u)$. Equivalently, 
$\g^* = \l^* \oplus \u^*$, where $\l^*$ is a Lie subalgebra and 
$\u^*$ is an ideal of 
$\g^*$. We will call such a quadruple $(\g,\delta,\l,\u)$ a {\it Lie 
bialgebra with a splitting.}

Assume further that $\g = \l\oplus \u$ is nondegenerate as a Lie algebra
with a splitting (see \cite{EE2}); in other words, we assume that for
$\lambda\in\l^*$ generic, the pairing $\omega_0(\lambda) : 
\u\times \u\to\CC$, $(u,v)\mapsto \lambda([u,v]_\l)$ is nondegenerate 
(here $z\mapsto z_\l$ is the projection on $\l$ parallel to $\u$). 

As we explained in \cite{EE2}, the nondegeneracy assumption may be expressed as
follows. We associate an element $p_\l^\g\in S^d(\l)$ ($d:= \on{dim}(\u)$)
to $\g = \l\oplus \u$, well defined up to multiplication by a nonzero scalar:
fix a linear isomorphism $i : \u \to \u^*$, define $a(\lambda)
\in \on{End}(\u)$ by $a(\lambda)(u) = 
i^{-1}\big( \omega_0(\lambda)(u,-)\big)$ for any $\lambda\in \g^*$, 
then $p_\l^\g(\lambda)
= \on{det}(a(\lambda))$. The nondegeneracy assumption is that 
$p_\l^\g\neq 0$. 

\begin{example} \label{example:21}
$(\g,\delta)$ is a semisimple Lie algebra with its standard 
Lie bialgebra structure, $\l\subset\g$ is a Levi subalgebra, $\u = \u_+ \oplus
\u_-$, 
$\u_\pm = \oplus_{\alpha\in \Delta_+(\g) \setminus
\Delta_+(\l)}\g_{\pm \alpha}$. 
More generally, $\g = \l\oplus \u_+ \oplus \u_-$ is a polarized nondegenerate
Lie algebra (see \cite{EE2}), and $\delta(x) = [x\otimes 1 + 1 \otimes x,r]$, 
where $r\in (\l\otimes\l) \oplus(\u_+\otimes \u_-)$ is such that 
$r + r^{2,1}$ is invariant and $\on{CYB}(r) = 0$.  
\end{example}

\subsection{Localizations}

Assume that $L^* \subset G^*$ is a subgroup corresponding to 
$\l^* \subset \g^*$. E.g., $L^*,G^*$ may be the formal groups associated 
with these Lie algebras; we denote by $\CC[[L^*]] = U(\l^*)^*$ the function 
ring of the formal group $L^*$. If $x\in L^*$, define 
$\omega_x : \u \times \u \to\CC$ by 
\begin{equation} \label{def:omega}
\omega_x(u,v) = \langle \big( \on{Ad}(x^{-1})(u) \big)_\g, 
\big( \on{Ad}(x^{-1})(v) \big)_{\g^*}
\rangle, 
\end{equation}
where $\langle -,-\rangle$ is the canonical bilinear form of 
$\d(\g) = \g \oplus \g^*$, and $x\mapsto x_\g,x_{\g^*}$ are the projections on 
$\g,\g^*$. Since the bilinear form on $\d(\g)$ is Ad-invariant and 
$\u$ is an isotropic subspace of $\d(\g)$, $\omega_x$ is antisymmetric. 

If $x\in L^*$,\footnote{if 
$L^*$ is a formal group, this means that we fix an Artinian 
$\CC$-ring $\Lambda$, and $x\in L^*(\Lambda) = \{$ring homomorphisms 
$U(\l^*)^* \to \Lambda\}$} 
define $a(x) \in \on{End}(\u)$ by 
$a(x)(u)  = i^{-1}(\omega_x(u,-))$ and set 
\begin{equation} \label{def:P}
P_\l^\g(x):= \on{det}(a(x)).
\end{equation} 
Then $P_\l^\g(x)$ is an element of $\CC[[L^*]]$, well-defined up to multiplication 
by a nonzero scalar. Let $\m_{L^*}$ be the maximal ideal of $\CC[[L^*]]$
(it coincides with the kernel of the counit map). Then $P_\l^\g\in 
\m_{L^*}^d$. Recall that the associated graded ring of $\CC[[L^*]]$ 
w.r.t. its filtration by the powers of $\m_{L^*}$ is $S^\cdot(\l)$. 
Then the class of $P_\l^\g$ in $\m_{L^*}^d / \m_{L^*}^{d+1}$ coincides with 
$p_\l^\g$. Therefore $P_\l^\g$ is nonzero. 

The localized ring $\CC[[L^*]][1/P_\l^\g]$ is also filtered by 
$\on{val}(f/(P_\l^\g)^k) = \on{val}(f) - kd$, and the associated 
graded ring is $S^\cdot(\l)[1/p_\l^\g]$. 

Any linear lift of the projection map $\m_{L^*} \to \m_{L^*} / \m_{L^*}^2 = 
\l$ induces a (noncanonical) identification of $\CC[[L^*]]$ with 
the formal series completion $\wh S^\cdot(\l)$ of the symmetric algebra 
of $\l$, and an injection of $\CC[[L^*]][1/P_\l^\g]$ in the degree completion 
of $S^\cdot(\l)[1/p_\l^\g]$, which extends to an identification of the valuation
completion of $\CC[[L^*]][1/P_\l^\g]$ with this degree completion.

\subsection{The $r$-matrix $\sigma_\l^\g$} \label{sect:sigma}

If $x\in L^*$ is such that $P_\l^\g(x)$ is invertible\footnote{
in the formal case, we fix a pro-Artinian ring without zero divisors
$\Lambda = \lim_{\leftarrow}\Lambda_n$, $x\in L^*(\Lambda)
= \lim_{\leftarrow} L^*(\Lambda_n)$, and $P_\l^\g(x)$
is assumed nonzero in $\Lambda$; then $\sigma_\l^\g(x) \in \wedge^2(\u)
\otimes \on{Frac}(\Lambda)$}, 
then $\omega_x \in \wedge^2(\u)^*$
is invertible, and we set $\sigma_\l^\g(x) := (\omega_x)^{-1}$. Then 
$$
\sigma_\l^\g \in \wedge^2(\u) \otimes \CC[[L^*]][1/P_\l^\g]. 
$$
$\sigma_\l^\g$ is uniquely determined by the following equivalent 
conditions (we set 
$\sigma_\l^\g = \sum_\alpha a_\alpha \otimes b_\alpha \otimes \ell_\alpha$):  

1) if $u,v\in \u$, then $\omega_x(u,v) = \sum_\alpha \omega_x(u,a_\alpha)
\omega_x(v,b_\alpha) \ell_\alpha(x)$ (equality in $\CC[[L^*]]$)

2) for any $v\in\g$, $\sum_\alpha a_\alpha \otimes
\omega_x(v,b_\alpha)\ell_\alpha(x) = v \otimes 1$
(equality in $\u\otimes \CC[[L^*]]$).

\begin{theorem}
$\sigma_\l^\g$ is a PL dynamical $r$-matrix for $(L^*,\g,0)$. 
\end{theorem}

{\em Proof.}
Let us first prove the $\l$-invariance of $\sigma_\l^\g$. 
We first prove the identity: 
\begin{equation} \label{pre:invce}
\forall u,v\in \u, \forall a\in \l, 
\quad 
\on{dress}_a^L \omega_x(u,v) + 
\omega_x([u,a],v) + \omega_x(u,[v,a]) = 0. 
\end{equation} 
Set $a' = \on{Ad}(x^{-1})(a)$ (so $a'\in \l\oplus \l^*$), 
$u'= \Ad(x^{-1})(u)$, $v' = \Ad(x^{-1})(v)$, then 
$$
\on{dress}_a^L \omega_x(u,v) = 
\langle ([a'_{\l^*},u'])_\g, v'_{\g^*} \rangle 
+ \langle u'_\g, ([a'_{\l^*},v'])_{\g^*}\rangle, 
$$ 
and
$$
\omega_x([u,a],v) + \omega_x(u,[v,a]) =
\langle ([u',a'])_\g , v'_{\g^*}\rangle 
+ \langle u'_\g, ([v',a'])_{\g^*}\rangle. 
$$ 
The sum of these terms is 
\begin{equation} \label{tobe0}
\langle [u',a'_\l]_\g,v'_{\g^*} \rangle 
+ \langle u'_\g, [v',a'_\l]_{\g^*} \rangle.
\end{equation}  
Here $a'_\l$ is the projection of $a'$ on $\l$ parallel to $\l^*$. 
Let us prove that (\ref{tobe0}) is zero. We have 
$[\l,\u\oplus\u^*] \subset \u\oplus \u^*$, therefore 
$u',v'\in\u\oplus\u^*$.  
If $v'\in\u$, (\ref{tobe0}) is clearly $0$; 
if $u'\in\u^*$, (\ref{tobe0}) is 
$\langle [u',a'_\l]_\g,v'_{\g^*} \rangle 
= \langle [u',a'_\l],v'_{\g^*} \rangle 
= \langle a'_\l,[v'_{\g^*},u'] \rangle$, 
which is zero since $\u^*$ is an ideal of $\g^*$. 
If $u'\in\u$ and $v'\in\u^*$, then the vanishing
of (\ref{tobe0}) follows from the invariance of 
$\langle-,-\rangle$. This proves (\ref{pre:invce}). 

Let us now prove that 
\begin{equation} \label{wanted:invce}
\forall a\in \l, \quad 
\on{dress}_a^L \sigma_\l^\g + [a^1 + a^2,\sigma_\l^\g] = 0. 
\end{equation}
Write $\sigma_\l^\g = \sum_\alpha a_\alpha \otimes b_\alpha 
\otimes \ell_\alpha$. 
Then (\ref{pre:invce}) implies that for any $u,v\in\u$, we have 
$$
\sum_\alpha \omega_x([u,a],a_\alpha) \omega_x(v,b_\alpha) \ell_\alpha(x) 
+  \omega_x(u,a_\alpha) \omega_x([v,a],b_\alpha) \ell_\alpha(x) 
= - \on{dress}^L_a \omega_x(u,v) ,  
$$
and using again (\ref{pre:invce}), we get 
\begin{align} \label{interm}
& \sum_\alpha \big( - \on{dress}^L_a \big( \omega_x(u,a_\alpha) \big) 
\omega_x(v,b_\alpha) 
\ell_\alpha(x) 
- \omega_x(u,a_\alpha) \on{dress}^L_a \big( \omega_x(v,b_\alpha) \big)  
\ell_\alpha(x) 
\\ & \nonumber + \omega_x(u,[a,a_\alpha]) \omega_x(v,b_\alpha) \ell_\alpha(x)
+ \omega_x(u,a_\alpha) \omega_x(v,[a,b_\alpha]) \ell_\alpha(x) \big) = 
- \on{dress}^L_a \omega_x(u,v).
\end{align} 
Applying $\on{dress}^L_a$ to the identity 
$\sum_\alpha \omega_x(u,a_\alpha) \omega_x(v,b_\alpha) \ell_\alpha(x) 
= \omega_x(u,v)$, and adding the resulting identity to 
(\ref{interm}), we get 
\begin{equation} \label{wanted}
\sum_\alpha \omega_x(u,a_\alpha) \omega_x(v,b_\alpha) \on{dress}^L_a
\ell_\alpha(x) + \omega_x(u,[a,a_\alpha]) \omega_x(v,b_\alpha) \ell_\alpha(x)
+ \omega_x(u,a_\alpha) \omega_x(v,[a,b_\alpha]) \ell_\alpha(x) = 0.
\end{equation} 
This is the result of pairing (\ref{wanted:invce}) with $u\otimes v$
using $\omega_x \otimes \omega_x$. This proves (\ref{wanted:invce}).

\medskip 

Let us now show that $\sigma_\l^\g$ satisfies the Poisson-Lie CDYBE.
The nontrivial components of this identity lie in the direct sum of 
$\wedge^3(\u) \otimes \CC[[L^*]]$ and of $ \big( (\l\otimes \wedge^2(\u)) 
\oplus (\u\otimes \l\otimes \u) \oplus (\wedge^2(\u) \otimes \l)\big) 
\otimes \CC[[L^*]]$. Since the l.h.s. of the PL CDYBE is obviously invariant
under the cyclic permutations of the first three tensor factors, 
it suffices to show the identities 
\begin{equation} \label{u:u:u}
(p_\u^{\otimes 3} \otimes \id)(\on{l.h.s.\ of\ the\ PL\ CDYBE}) = 0 
\end{equation}
and 
\begin{equation} \label{l:u:u}
(p_\l\otimes p_\u^{\otimes 2} \otimes \id)
(\on{l.h.s.\ of\ the\ PL\ CDYBE}) = 0 . 
\end{equation}
Here $p_\l,p_\u$ are the projections $\g\to\l$, $\g\to\u$ corresponding to 
$\g = \l\oplus\u$.

\medskip 

Let us prove (\ref{u:u:u}). 
Let $u,v,w$ be arbitrary elements of $\u$. For $z\in \{u,v,w\}$, 
we set $\overline z := \Ad(x^{-1})(z)$, which we decompose as 
$\overline z = z' + z''$, $z'\in\g$, $z''\in\g^*$.  
Then we have $\langle [\overline u,\overline v], \overline w\rangle
= \langle [u,v], w\rangle = 0$ (since $\langle-,-\rangle$ is 
Ad-invariant and $\g$ is isotropic). On the other hand,  
$\langle [\overline u,\overline v], \overline w\rangle$ can be 
expanded (again using isotropies and invariances) as 
$\langle [u',v'],w'' \rangle + \langle [u'',v''], w'\rangle 
+ \on{cycl.\ perm}$. Therefore
\begin{equation} \label{tobeusedlater}
\langle [u',v'],w'' \rangle + \langle [u'',v''], w'\rangle 
+ \on{cycl.\ perm.} = 0. 
\end{equation}
Let us now extend $\omega_x$ to a bilinear map $\wt\omega_x : 
\d(\g)\times \d(\g)\to \CC$ using formula (\ref{def:omega}). 
Since $\on{Ad}(x^{-1})(\l) \subset 
\l\oplus \l^*$ and $\on{Ad}(x^{-1})(\u) \subset \u\oplus \u^*$, 
we get 
\begin{equation} \label{pty:wt:omega}
\wt\omega_x(\ell,a) = 0 \on{\ whenever\ } 
\ell\in\l\on{\ and\ }a\in\u.
\end{equation}

Now 
\begin{align*}
& \sum_\alpha \sum \omega_x(p_\u(a_\alpha^{(1)}),u) 
\omega_x(p_\u(a_\alpha^{(2)}),v) \omega_x(b_\alpha,w) \ell_\alpha(x)
\\ & 
=  \sum \omega_x(p_\u(w^{(1)}),u) \omega_x(p_\u(w^{(2)}),v)
=  \sum \wt\omega_x(w^{(1)},u) \wt\omega_x(w^{(2)},v) \on{\ (using\ 
(\ref{pty:wt:omega}))}
\\ & = \sum  \langle \big( \Ad(x^{-1})(w^{(1)}) \big)_\g, 
\big( \Ad(x^{-1})(u) \big)_{\g^*}\rangle 
\langle \big( \Ad(x^{-1})(w^{(2)}) \big)_\g, 
\big( \Ad(x^{-1})(v) \big)_{\g^*}\rangle \\ 
& 
=  \sum \langle \Ad(x^{-1})(w^{(1)}),  
\big( \Ad(x^{-1})(u) \big)_{\g^*}\rangle 
\langle \Ad(x^{-1})(w^{(2)}),  
\big( \Ad(x^{-1})(v) \big)_{\g^*}\rangle \on{(since\ }\g^*
\on{\ is\ isotropic)}\\ 
& 
=  \sum \langle w^{(1)},  
\Ad(x) \circ p_{\g^*} \big( \Ad(x^{-1})(u) \big)\rangle 
\langle w^{(2)},  
\Ad(x) \circ p_{\g^*}\big( \Ad(x^{-1})(v) \big)\rangle \on{(since\
}\langle -,-\rangle
\on{\ is\ Ad-invariant)}\\ 
& 
=  
\langle w, \big[   
\Ad(x) \circ p_{\g^*} \big( \Ad(x^{-1})(u) \big), 
\Ad(x) \circ p_{\g^*}\big( \Ad(x^{-1})(v) \big) \big] \rangle \on{(since\
}\Ad(x)(\g^*) \subset \g^*)
\\ 
& 
=  
\langle \Ad(x^{-1})(w), \big[   
p_{\g^*} \big( \Ad(x^{-1})(u) \big),
p_{\g^*}\big( \Ad(x^{-1})(v) \big) \big] \rangle \on{(since\
}\Ad(x) \on{\ is\ an\ automorphism)}
\\ & = \langle \overline w, [u'',v'']\rangle = \langle w',[u'',v'']\rangle 
\on{\ (since\ }\g^*\on{\ is\ isotropic)};   
\end{align*}
here $p_{\g^*}$ is the projection $\g\oplus \g^*\to\g^*$, 
$p_\u$ is the projection $\g = \l\oplus \u \to \u$, and we write $\delta(w)
= w^{(1)} \otimes w^{(2)}$. 

On the other hand, 
\begin{align*}
& \sum_{\alpha,\beta}
\omega_x(p_\u([a_\alpha,a_\beta]),u) \omega_x(b_\alpha,v) 
\omega_x(b_\beta,w) \ell_\alpha\ell_\beta(x)
\\ & 
= \omega_x(p_\u([v,w]),u) = \wt\omega_x([v,w],u) \on{\ (using\ 
(\ref{pty:wt:omega}))}
\\ & 
=  \langle [\Ad(x^{-1})(v),\Ad(x^{-1})(w)], 
p_{\g^*}\circ \Ad(x^{-1})(u) \rangle 
\on{\ (since\ }\g^*\on{\ is\ isotropic)}
\\ & = \langle [\overline v,\overline w], u''\rangle
= 
\langle [v',w'], u''\rangle
+ \langle v',[w'',u''] \rangle + \langle w',[u'',v''] \rangle
\end{align*}
(using the isotropies of $\g$ and $\g^*$, and the invariance of
$\langle-,-\rangle$). Taking the difference of
these results, we get  
\begin{align*}
& \sum_{\alpha,\beta}
\omega_x(p_\u([a_\alpha,a_\beta]),u) \omega_x(b_\alpha,v) 
\omega_x(b_\beta,w) \ell_\alpha\ell_\beta(x)
- \sum_\alpha \omega_x(p_\u(a_\alpha^{(1)}),u) 
\omega_x(p_\u(a_\alpha^{(2)}),v) \omega_x(b_\alpha,w) \ell_\alpha(x)
\\ & 
% \omega_x(p_\u([v,w]),u) - \omega_x(p_\u(w^{(1)}),u) 
% \omega_x(p_\u(w^{(2)}),v)
= \langle [v',w'], u''\rangle + \langle v',[w'',u''] \rangle. 
\end{align*} 
Therefore
\begin{align} \label{preprefinal} 
& \nonumber 
\sum_{\alpha,\beta}
\omega_x(p_\u([a_\alpha,a_\beta]),u) \omega_x(b_\alpha,v) 
\omega_x(b_\beta,w) \ell_\alpha\ell_\beta(x)
- \sum_\alpha \omega_x(p_\u(a_\alpha^{(1)}),u) 
\omega_x(p_\u(a_\alpha^{(2)}),v) \omega_x(b_\alpha,w) \ell_\alpha(x)
\\ & 
+ \on{\ cycl.\  perm.\ } 
% \omega_x(p_\u([v,w]),u) - \omega_x(p_\u(w^{(1)}),u) 
% \omega_x(p_\u(w^{(2)}),
= \langle [v',w'], u''\rangle + \langle v',[w'',u''] \rangle
+ \on{\ cycl.\  perm.\ } = 0,  
\end{align} 
where the last equality follows from (\ref{tobeusedlater}). 
Now the l.h.s. of (\ref{preprefinal}) is the pairing 
of 
$(p_\u^{\otimes 3}\otimes \id)(\on{CYB}(\rho) 
+ \on{Alt}(\delta\otimes \id)(\rho))$
with $u\otimes v \otimes w$ (using $\omega_x^{\otimes 3}$). 
Therefore this pairing is zero, for any $u,v,w\in\u$. It follows that 
$$
(p_\u^{\otimes 3}\otimes \id)(\on{CYB}(\rho) 
+ \on{Alt}(\delta\otimes \id)(\rho)) = 0. 
$$
This proves (\ref{u:u:u}). 

\medskip 

Let us prove (\ref{l:u:u}). For this, we will prove that for any $u,v\in\u$, the
pairing (using $\wt\omega_x^{\otimes 2}$) of the r.h.s. of (\ref{l:u:u})
with $\id\otimes u \otimes v$ is zero. This pairing is 
\begin{align} \label{pairing}
& \nonumber 
\sum_{\alpha,\beta} p_\l([a_\alpha,a_\beta]) \otimes 
\omega_x(b_\alpha,u) \omega_x(b_\beta,v) \ell_\alpha\ell_\beta(x)
- \sum_{i,\alpha} e_i \otimes \omega_x(a_\alpha,u) \omega_x(b_\alpha,v) 
\LL_{\eps_i} \ell_\alpha(x)
\\ &  
+ \sum_\alpha p_\l(a_\alpha^{(1)}) \otimes \omega_x(p_\u(a_\alpha^{(2)}),u)
\omega_x(b_\alpha,v) \ell_\alpha(x)
- \sum_\alpha p_\l(a_\alpha^{(1)}) \otimes 
\omega_x(p_\u(a_\alpha^{(2)}),v) \omega_x(b_\alpha,u) \ell_\alpha(x). 
\end{align}
The first sum of (\ref{pairing}) is equal to $p_\l([u,v]) \otimes 1$. 

The second sum of (\ref{pairing}) is equal to 
\begin{align} \label{two:next}
& \nonumber 
- \sum_{i,\alpha} e_i \otimes \LL_{\eps^i} \big( \omega_x(a_\alpha,u)
\omega_x(b_\alpha,v) \ell_\alpha(x)\big)
+ \sum_{i,\alpha} e_i \otimes \LL_{\eps^i}(\omega_x(a_\alpha,u))
\omega_x(b_\alpha,v) \ell_\alpha(x) 
\\ & \nonumber  
+ \sum_{i,\alpha} e_i \otimes \omega_x(a_\alpha,u)
\LL_{\eps^i}(\omega_x(b_\alpha,v)) \ell_\alpha(x) 
= - \sum_{i,\alpha} e_i \otimes \LL_{\eps^i}(\omega_x(u,v))
\\ & \nonumber
+ \sum_{i,\alpha} e_i \otimes \wt\omega_x([\eps^i,a_\alpha],u)
\omega_x(b_\alpha,v)\ell_\alpha(x) 
+ \sum_{i,\alpha} e_i \otimes \wt\omega_x(a_\alpha,[\eps^i,u])
\omega_x(b_\alpha,v)\ell_\alpha(x) 
\\ & \nonumber 
+ \sum_{i,\alpha} e_i \otimes \omega_x(a_\alpha,u)
\wt\omega_x([\eps^i,b_\alpha],v)\ell_\alpha(x) 
+ \sum_{i,\alpha} e_i \otimes \omega_x(a_\alpha,u)
\wt\omega_x(b_\alpha,[\eps^i,v])\ell_\alpha(x) 
\\ & \nonumber
= - \sum_i e_i \otimes \wt\omega_x([\eps^i,u],v) - \sum_i e_i \otimes 
\wt\omega_x(u,[\eps^i,v])
- \sum_i e_i \otimes \wt\omega_x([\eps^i,v],u) - \sum_i e_i \otimes 
\wt\omega_x(v,[\eps^i,u])
\\ & 
+ \sum_i e_i \otimes \wt\omega_x([\eps^i,u],v) + \sum_i e_i \otimes 
\wt\omega_x(u,[\eps^i,v])
= - \sum_i e_i \otimes \wt\omega_x([\eps^i,v],u) - \sum_i e_i \otimes 
\wt\omega_x(v,[\eps^i,u])
\end{align}
where both equalities follows from the identity 
$$
\forall \eps\in\h^*,\forall\alpha,\beta\in \d(\g), \; 
\LL_\eps(\wt\omega_x(\alpha,\beta)) 
= \wt\omega_x([\eps,\alpha],\beta) + \wt\omega_x(\alpha,[\eps,\beta]). 
$$
The two last sums of (\ref{pairing}) give the contribution 
$$
- \sum p_\l(v^{(1)}) \otimes \omega_x(p_\u(v^{(2)}),u) 
+ \sum p_\l(u^{(1)}) \otimes \omega_x(p_\u(u^{(2)}),v). 
$$ 
Now if $x\in \u$, then 
$$
\sum x^{(1)} \otimes x^{(2)} = \delta(x) = 
\sum_i [x,e_i] \otimes \eps^i + \sum_i e_i \otimes [x,\eps^i]
+ \sum_j [x,u_j] \otimes \lambda^j + \sum_j u_j \otimes [x,\lambda^j] ,
$$
where $(u_j),(\lambda^j)$ are dual bases of $\u$ and $\u^*$. Therefore 
$$
\sum p_\l(x^{(1)}) \otimes p_\u(x^{(2)}) = \sum_i e_i \otimes p_\u([x,\eps^i]), 
$$
so the contribution of the two last sums of (\ref{pairing}) is 
\begin{equation} \label{two:last}
- \sum_i e_i \otimes \omega_x(p_\u([v,\eps^i]),u) 
+ \sum_i e_i \otimes \omega_x(p_\u([u,\eps^i]),v). 
\end{equation}
The sum of (\ref{two:next}) and (\ref{two:last}) is equal to 
\begin{equation} \label{4:last}
\sum_i e_i \otimes \wt\omega_x(p_{\u^*}([v,\eps^i]),u)
+ \sum_i e_i \otimes \wt\omega_x(p_\u([u,\eps^i]),v)
- \sum_i e_i \otimes \wt\omega_x(v,[\eps^i,u]). 
\end{equation}
Now $\wt\omega_x(\alpha,\beta) = 0$ whenever $\alpha\in\g^*$, so the 
first term of (\ref{4:last}) vanishes; we have 
\begin{align*}
& \wt\omega_x(p_\u([u,\eps^i]),v) = \langle 
p_\g \circ \Ad(x^{-1}) \circ p_\u([u,\eps^i]), p_{\g^*} \circ \Ad(x^{-1})(v)
 \rangle 
\\ & 
= \langle p_\g \circ \Ad(x^{-1}) ([u,\eps^i]), p_{\g^*} \circ \Ad(x^{-1})(v)
\rangle 
\end{align*}
(since $\Ad(x^{-1})(\g^*) \subset \g^*$ and $p_\g(\g^*) = 0$), and 
$$
- \wt\omega_x(v,[\eps^i,u]) = \langle p_\g\circ \Ad(x^{-1})(v), 
p_{\g^*} \circ \Ad(x^{-1})([u,\eps^i])\rangle; 
$$
adding up these equalities, we get 
$$
(\ref{4:last}) = \sum_i e_i \otimes 
\langle \Ad(x^{-1})(v),\Ad(x^{-1})([u,\eps^i])\rangle = 
\sum_i e_i \otimes \langle [v,u],\eps^i \rangle = p_\l([v,u])\otimes 1. 
$$
Therefore $(\ref{pairing}) = 0$, which proves (\ref{l:u:u}). 
\hfill \qed \medskip

\begin{remark} The fact that $\sigma_\l^\g$ is equivariant implies that 
$\on{dress}^L_a(P_\l^\g) = \chi(a) P_\l^\g$, where $\chi : \l \to \CC$ 
is the character
defined by $\chi(a) = \on{tr}(\on{ad}(a)_{|\u})$. 
\end{remark}

\begin{remark}
Let $\CC[[L^*]][1/P_\l^\g]_{\geq i}$ be the valuation $\geq i$ part of 
$\CC[[L^*]][1/P_\l^\g]$. Then $\sigma_\l^\g$ belongs to $\wedge^2(\u)
\otimes \CC[[L^*]][1/P_\l^\g]_{\geq -1}$. Its image in 
$\wedge^2(\u) \otimes \on{gr}(\CC[[L^*]][1/P_\l^\g])_{-1}
= \wedge^2(\u) \otimes S^\cdot(\l)[1/p_\l^\g]_{-1}$ is $-r_\l^\g$, where 
$r_\l^\g$ is the $r$-matrix introduced in \cite{EE2}. Taking the images in 
$\on{gr}(\CC[[L^*]][1/P_\l^\g])_{-1}$ (resp., 
$\on{gr}(\CC[[L^*]][1/P_\l^\g])_{-2}$)
of the invariance (resp., CDYBE) identities for $\sigma_\l^\g$, 
we find the invariance (resp., CDYBE) identities for $r_\l^\g$, i.e., 
Proposition 1.1 in \cite{EE2}. 
\end{remark}

\begin{example}
If $(\g,\delta)$ is as in Example \ref{example:21} and $\l = \h$ is the Cartan
subalgebra of $\g$, then $\sigma_\l^\g$ is given by 
$$
\forall \lambda\in\h^*, \quad 
\sigma_\h^\g(e^\lambda) = 
\sum_{\alpha\in\Delta_+(\g)} {{e_\alpha \wedge f_\alpha}\over{
1 - e^{-(\lambda,\alpha)} }}, 
$$
which is the solution of Example \ref{from:old}. 
\end{example}

\begin{example}
Assume that $(\g,\delta)$ is semisimple and $\l\subset\g$ is a Levi 
subalgebra as in Example \ref{example:21}. The Cartan decomposition of $\l$
is $\l = \h\oplus \m_+ \oplus \m_-$, where 
$\m_+ = \oplus_{\alpha\in \Delta_+(\l)} \CC e_\alpha$, 
$\m_- = \oplus_{\alpha\in \Delta_+(\l)} \CC f_\alpha$. 
Let $G$ be a (formal or algebraic) 
group with Lie algebra $\g$, $H \subset L_\pm \subset L \subset G$
be subgroups corresponding to $\h\subset (\h \oplus \m_\pm) \subset \l \subset 
\g$. Then $L^* = \{(x_-,x_+)\in L_- \times L_+ | h_\pm(x_+)h_\pm(x_-) = 1\}$, 
where 
$h_\pm : L_\pm \to H$ are the group morphisms corresponding to the projections 
$(\h\oplus \m_\pm) \to \h$. If $x\in L$, then $\Ad(x)$ restricts to an
automorphism of $\u$. For a generic $x\in L$, $\Ad(x)-\id$ is invertible, 
when viewed as an endomorphism of $\u$. Then $\sigma_\l^\g : L^* \to 
\wedge^2(\u)$ is defined by 
$$
\sigma_\l^\g(x_+,x_-) = \sum_{\alpha\in \Delta_+(\g) \setminus \Delta_+(\l)}
\big( \id - \Ad(x_+x_-^{-1}) \big)^{-1}(e_\alpha) \wedge f_\alpha. 
$$
\end{example}

\subsection{The composition theorem} \label{sect:2:4}

\begin{theorem} \label{thm:comp}
Let $(\g,\l,\u,\delta)$ be a nondegenerate Lie bialgebra 
with a splitting (see Section \ref{LBA:split}). Let $(\a,\delta)$ 
be a Lie bialgebra equipped with a Lie bialgebra inclusion 
$\g\subset \a$. Let $Z\in \wedge^3(\a)^\a$. Then the map 
$\rho \mapsto \rho_{|L^*} + \sigma_\l^\g$ takes $\{$PL dynamical $r$-matrices
for $(G^*,\a,Z)\}$ to $\{$PL dynamical $r$-matrices for $(L^*,\a,Z)\}$. 
\end{theorem}

{\em Proof.} This is a consequence of the following statement: if $\tau : G^*
\to \wedge^2(\a)$ is such that $\on{dress}_a^L(\tau) + [a^1 + a^2,\tau] = 0$
for any $a\in\g$, then 
\begin{equation} \label{rest:form}
(\on{d}^L\tau)_{|L^*} = \on{d}^L(\tau_{|L^*})
+ [(\tau_{|L^*})^{1,2},(\sigma_\l^\g)^{1,3} + (\sigma_\l^\g)^{2,3}]. 
\end{equation}
Here the map $\g^{\otimes 3} \otimes \CC[[G^*]] \to \g^{\otimes 3} \otimes 
\CC[[L^*]]$, $\omega \mapsto \omega_{|L^*}$ is the tensor product of 
$\id^{\otimes 3}$ with the restriction map $\CC[[G^*]] \to \CC[[L^*]]$. 

Let us prove (\ref{rest:form}). Let $(u_j),(\lambda^j)$ be dual bases of $\u$
and $\u^*$. Set $\tau = \sum_\alpha x_\alpha \otimes y_\alpha \otimes
m_\alpha$. Then 
$$
(\on{d}^L\tau)_{|L^*} - \on{d}^L(\tau_{|L^*})
= - \sum_{\alpha,j} x_\alpha \otimes y_\alpha \otimes u_j 
\otimes \LL_{\lambda^j}(m_\alpha)_{|L^*}
= \sum_\alpha x_\alpha \otimes y_\alpha \otimes u_j \otimes
\RR_{\on{Ad}(x^{-1})(\lambda^j)}(m_\alpha)_{|L^*}. 
$$
On the other hand, 
\begin{align*}
& [(\tau_{|L^*})^{1,2},(\sigma_\l^\g)^{1,3} + (\sigma_\l^\g)^{2,3}] 
= \sum_{\alpha,\beta} [x_\alpha,a_\beta] \otimes y_\alpha \otimes b_\beta 
\otimes (m_\alpha \ell_\beta)_{|L^*} + x_\alpha \otimes [y_\alpha,a_\beta]
\otimes b_\beta \otimes (m_\alpha \ell_\beta)_{|L^*}
\\ & 
= \sum_{\alpha,\beta} x_\alpha \otimes y_\alpha
\otimes b_\beta \otimes \big( \ell_\beta 
\on{dress}_{a_\beta}^L(m_\alpha)\big)_{|L^*} \on{\ (using\ the\ 
equivariance\ of\ }\tau\on{)}
\\ & = - \sum_{\alpha,\beta} x_\alpha \otimes y_\alpha \otimes b_\beta \otimes 
\big( \ell_\beta \RR_{(\Ad(x^{-1})(a_\beta))_{\g^*}}
(m_\alpha)\big)_{|L^*}. 
\end{align*}
Therefore (\ref{rest:form}) follows from 
\begin{equation} \label{sbk}
\sum_\beta \Big( \big( \Ad(x^{-1})(a_\beta) \big)_{\g^*} 
\otimes b_\beta\Big) \ell_\beta(x) 
= - \sum_j \Ad(x^{-1})(\lambda^j) \otimes u_j. 
\end{equation}  
Let us prove (\ref{sbk}). Both sides of this equation belong to 
$\g^* \otimes \u \otimes \CC[[L^*]]$. We will show that the pairings of both sides
of (\ref{sbk}) with $b\otimes \id\otimes \id$ are the same, for any $b\in\g$. As
$x$ is formal near the origin, the map $p_\g \circ \Ad(x^{-1}) : \g \to \g$ is
bijective, so we have $b = p_\g \circ \Ad(x^{-1})(a)$, where $a\in\g$. 
Pairing the l.h.s. of (\ref{sbk}) with $b\otimes \id \otimes \id$, 
we get 
$$
\sum_\beta b_\beta \otimes \langle p_\g^* \circ \Ad(x^{-1})(a_\beta), 
p_\g \circ \Ad(x^{-1})(a) \rangle \ell_\beta(x) = \sum_\beta 
b_\beta \otimes \wt\omega_x(a,a_\beta) \ell_\beta(x) = - p_\u(a) \otimes 1, 
$$
and pairing the r.h.s. of (\ref{sbk}) with the same element, we get 
\begin{align*}
& - \sum_j u_j \otimes \langle \Ad(x^{-1})(\lambda^j), 
p_\g \circ \Ad(x^{-1})(a) \rangle  
= 
- \sum_j u_j \otimes \langle \Ad(x^{-1})(\lambda^j), \Ad(x^{-1})(a) \rangle 
\on{\ (as\ } \g^* \on{\ is\ isotropic)} 
\\ & 
= - \sum_j \langle \lambda^j, a\rangle u_j \otimes 1 = - p_\u(a) \otimes 1.
\end{align*}
This proves (\ref{sbk}). \hfill \qed \medskip 

\begin{remark} \label{rem:comp}
Assume that $\a = \g \oplus \vv$ is a nondegenerate Lie bialgebra with a 
splitting. Then $\sigma_\l^\a = (\sigma_\g^\a)_{|L^*} + \sigma_\l^\g$, 
in other words, the map of Theorem \ref{thm:comp} takes $\sigma_\g^\a$
to $\sigma_\l^\a$ (here $Z=0$). 
\end{remark}

\begin{example}
Applying Theorem \ref{thm:comp} when $(\g,\delta)$ is semisimple,  
$\l\subset \g$ is a Levi subalgebra as in Example \ref{example:21}, and $\rho$
is the Balog-Feh\'er-Palla $r$-matrix, we get the following result. 
Define $\rho_{\on{BFP}}^{\l,\g} : L^* \to \wedge^2(\g)$ by 
\begin{equation} \label{rho:l:g}
\rho_{\on{BFP}}^{\l,\g}(x_+,x_-) = 
\nu ({{\id + g^{2\nu} }\over{\id - g^{2\nu}}} \otimes \id)(t)
- {1\over 2} ({{\id + g}\over{\id - g}}\otimes \id)(t_\l)
+ {1\over 2}s - {1\over 2}s^{2,1} , 
\end{equation} 
where $g = \Ad(x_+x_-^{-1}) \in \End(\g)$, $t\in S^2(\g)^\g$, $t_\l
\in S^2(\l)^\l$ are the Casimir elements of $\g,\l$, and $s =
\sum_{\alpha\in \Delta_+(\g) \setminus \Delta_+(\l)} 
e_\alpha \otimes f_\alpha$,  
so that $t = t_\l + s + s^{2,1}$. Then $\rho_{\on{BFP}}^{\l,\g}$
is a PL dynamical $r$-matrix for $(L^*,\g,Z_\nu)$.  
\end{example}

\begin{example}
Rewriting (\ref{rho:l:g}) using $t = t_\l + s + s^{2,1}$, we get in the case 
$\l = \h$ (the Cartan subalgebra of $\g$) that 
$$
\rho_{\on{BFP}}^{\h,\g}(e^\lambda) = 
{1\over 2} \sum_{\alpha\in\Delta_+(\g)}
\big( 1 + 2\nu \on{coth} \nu(\lambda,\alpha)\big) 
(e_\alpha \wedge f_\alpha)
$$
is a PL dynamical $r$-matrix for $(H^*,\g,(4\nu^2 - 1)Z(r_0))$. 
Under the correspondence of 
Example \ref{from:old}, this solution corresponds to the usual dynamical 
$r$-matrix 
$\lambda \mapsto 2\nu r_\h^\g(2\nu\lambda)$ for $(\h,\g,4\nu^2 Z(r_0))$. 
\end{example}

\subsection{A PL dynamical $r$-matrix for $(G^*,\g,(4\nu^2-1)Z(r_0))$}

The proof of Theorem \ref{thm:comp} implies the following result: 
let $\rho : G^* \to \wedge^2(\g)$ be a $\g$-equivariant function, 
such that $\rho_{|L^*} + \sigma_\l^\g$ is a PL dynamical $r$-matrix for 
$(L^*,\g,Z)$, then $\rho$ is a PL dynamical $r$-matrix for 
$(G^*,\g,Z)$. This result leads to the following $r$-matrix (an 
analogue of the $r$-matrix in \cite{EE2}, Remark 1.8). 

Let $\g$ be a semisimple Lie algebra and $t\in S^2(\g)^\g$ be nondegenerate. 
Let $\h$ be a Cartan subalgebra of $\g$ and $\b_\pm$ be opposite
Borel subalgebras containing $\h$, $G$ a (formal or algebraic)
Lie group with Lie algebra 
$\g$, $H\subset B_\pm\subset G$ Lie subgroups corresponding to $\h\subset 
\b_\pm \subset \g$. Then $G^* = \{(x_-,x_+)\in B_- \times B_+ | 
h_-(x_-)h_+(x_+)
=1\}$. The map $G^* \to G$, $g^* = (x_-,x_+) \mapsto x_+x_-^{-1} = 
\lambda(g^*)$ is equivariant
for the left dressing action of $\g$ on $G^*$ and its adjoint action on $G$. 

Let $G_{\on{reg}} \subset G(\RR)$ be the subset of regular real 
semisimple elements, and 
define $G^*_{\on{reg}} := \lambda^{-1}(G_{\on{reg}})$. Let $\theta : 
]0,1[\cup ]1,\infty[ \to \{1,-1\}$ be defined by 
$\theta_{|]0,1[} = -1$, $\theta_{|]1,\infty[} = +1$. 
If $g'\in G_{\on{reg}}$, let $\h(g') = \{x\in \g | \on{Ad}(g')(x) = 0\}$
be the Cartan subalgebra of $\g$ associated to $g'$, and if $\h'\subset\g$
is a Cartan subalgebra, let $t_{\h'}$ be the part of $t$ corresponding to 
$\h'$. For $\nu\in\CC$, define $\rho_\nu : G^*_{\on{reg}} \to \wedge^2(\g)$ by 
$$
\rho_\nu(g^*) = -\big( \big(  
{1\over 2} {{\Ad(g) + \id}\over{\Ad(g)-\id}}
+ \nu \theta(\Ad(g))\big) \otimes \id \big) (t - t_{\h(g)}), 
$$
where $g = \lambda(g^*)$. Here $\theta(\Ad(g))$ is the operator acting
by $\theta(\lambda)$ on the $\lambda$-eigenspace of $\on{Ad}(g)$.  

\begin{proposition}
$\rho_\nu$ is a PL dynamical $r$-matrix for $(G^*,\g,(4\nu^2-1)Z(r_0))$. 
\end{proposition}

{\em Proof.} $\rho_\nu$ is clearly equivariant,  
$(\rho_{\nu})_{|H^*} + \sigma_\h^\g = (1+2\nu)r_0$, 
and $\on{CYB}((1+2\nu)r_0) + \on{Alt}(\delta\otimes \id)(
(1+2\nu)r_0) = (4\nu^2 -1)Z(r_0)$. \hfill \qed \medskip

\section{Quantization of $\sigma_\l^\g$} \label{sect:quant} \label{sect:3}

We introduce the notion of a polarized Lie bialgebra with a splitting. 
We will quantize $\sigma_\l^\g$ when $(\g,\delta,\l,\u)$ is polarized. 
This yields a quantization in the case of Example \ref{example:21}. 

\subsection{A construction for Hopf algebras} \label{sect:quot} 
\label{sect:3:1}

Let $A,B$ be Hopf algebras, $i  : B\hookrightarrow A$ and $\pi : A 
\twoheadrightarrow B$ be Hopf algebra morphisms, such that 
$\pi \circ i = \id_B$. Define maps 
$$
\alpha : A \to A, \; \beta : A \to A \otimes B 
$$
by 
$\alpha(a) = \sum a^{(1)} (i\circ \pi \circ S)(a^{(2)})$, and 
$\beta(a) = (\id\otimes \pi) \circ \Delta(a) - a \otimes 1$. 

\begin{proposition} (see \cite{Rad}, Theorem 3)

1) $\on{Im}(\alpha) = \on{Ker}(\beta)$. Set $C := \on{Im}(\alpha)
= \on{Ker}(\beta)$. 

2) $C$ is a subalgebra and a left coideal of $A$. 

3) $C$ is stable under the adjoint action of $A$ on itself. 
\end{proposition}

\begin{example} \label{ex:toldot}
Let $\p,\l$ be Lie algebras, and assume that 
we have Lie algebra morphisms $\l \stackrel{i}{\hookrightarrow} \p$
and $\p \stackrel{\pi}{\twoheadrightarrow} \l$, such that 
$\pi \circ i = \id_\l$. This means that $\p = \l\oplus \u$, $\l$
is a Lie subalgebra of $\p$, $\u$ is an ideal of $\p$. If we set 
$A = U(\p)$, $B = U(\l)$, we get $C = U(\u)$. This is proved using the 
isomorphism $U(\u) \otimes U(\l) \to U(\p)$, $u\otimes l\mapsto ul$. 
\end{example}

\begin{example} \label{ex:vayera}
In the situation of Example \ref{ex:toldot}, let $A = \CC[[P]]$, 
$B = \CC[[L]]$ (here $\CC[[X]] = U(\x)^*$ if $\x$ is a Lie algebra, and 
$X$ is the corresponding formal group). 
Then $C = \CC[[P/L]]$ is the ring of right $\l$-invariant formal 
functions on $P$. The natural map $U \to P/L$ is an isomorphism of formal
manifolds, so we will denote $\CC[[P/L]]$ as $\CC[[U]]$. 
\end{example}

\subsection{Lie bialgebras with a splitting} \label{sect:3:2}
\label{sect:q:lba}

Let $(\l,\delta)$ and $(\p,\delta)$ be Lie bialgebras, equipped with 
Lie bialgebra morphisms $i : \l\hookrightarrow \p$ and 
$\pi : \p\twoheadrightarrow \l$, such that $\pi\circ i = \id_\l$. 
We set $\u := \on{Ker}(\pi)$. Then $\u$ is a Lie algebra ideal of 
$\p$, such that 
$\delta(\u) \subset \big( (\l\otimes \u) \oplus (\u\otimes \l) \big)
\oplus \wedge^2(\u)$. 

\begin{example} \label{example:small}
$\p = \p_+$, a parabolic subalgebra of a semisimple Lie 
bialgebra $\g$, and $\l \subset \p$ is the corresponding 
Levi subalgebra. 
\end{example}

Assume that we have
constructed a pair of QUE algebras $U_\hbar(\l)$, $U_\hbar(\p)$ quantizing
$\l,\p$, and a sequence of morphisms $U_\hbar(\l) \stackrel{i}{\hookrightarrow}
U_\hbar(\p) \stackrel{\pi}{\twoheadrightarrow} U_\hbar(\l)$, whose composition 
is the identity, quantizing $\l\hookrightarrow \p \twoheadrightarrow \l$. This 
can be achieved: (a) in the general case, by applying a quantization functor 
to $\l\hookrightarrow \p \twoheadrightarrow \l$ (see \cite{EK}); 
(b) in the situation of Example \ref{example:21}, let $I \subset [1,r]$ 
be the set of
indices of $\p$. Then $U_\hbar(\p)$ is the subalgebra of $U_\hbar(\g)$
generated by the $h_i,e_i,i\in [1,r]$, and $f_i,i\in I$; 
$U_\hbar(\l) \subset U_\hbar(\p)$ is generated  by the $e_i,f_i,i\in I$, 
and $h_i,i\in [1,r]$. The morphism $\pi$ is defined by 
$\pi(x) = x$ for $x = e_i,f_i$, $i\in I$, and $x=h_i,i\in [1,r]$, 
and $\pi(e_i) = 0$ if $i\notin [1,r]$. 

\medskip 

Then we are in the situation of Section \ref{sect:quot}, with 
$A = U_\hbar(\p)$, $B = U_\hbar(\l)$. 
Define $U_\hbar(\u)$ to be the $C$ defined in this section. 

\begin{proposition} \label{prop:iashar} \label{prop:nemo}
$U_\hbar(\u) \subset U_\hbar(\p)$ is a flat deformation of 
$U(\u) \subset U(\p)$. 
\end{proposition}

{\em Proof.} $U_\hbar(\u)$ is defined as $\on{Ker}(\beta)$, therefore 
it is a complete, saturated submodule of $U_\hbar(\p)$ (here saturated means
that if $x\in U_\hbar(\p)$ is such that $\hbar x\in U_\hbar(\u)$, then 
$x\in U_\hbar(\u)$). So it remains to show that the morphism of 
reduction  modulo $\hbar$, $(\on{mod\ }\hbar) : U_\hbar(\p) \to U(\p)$, 
satisfies $(\on{mod}\ \hbar)(U_\hbar(\u)) = U(\u)$. 

Let $\alpha_0,\beta_0$ be the reductions modulo $\hbar$ of $\alpha,\beta$.  
According to Example \ref{ex:toldot}, $\on{Im}(\alpha_0) = 
\on{Ker}(\beta_0) = U(\u)$.  

Now $(\on{mod\ }\hbar)(U_\hbar(\u)) = (\on{mod\ }\hbar)(\on{Im}(\alpha)) \supset 
\on{Im}(\alpha_0)$, and $(\on{mod\ }\hbar)(U_\hbar(\u)) 
= (\on{mod\ }\hbar)(\on{Ker}(\beta)) \subset \on{Ker}(\alpha_0)$. 
Therefore $(\on{mod}\ \hbar)(U_\hbar(\u)) = U(\u)$. 
\hfill \qed \medskip 

We now describe a quantized formal series Hopf (QFSH) algebra
version of the above constructions. Recall that if $\a$ is a finite 
dimensional Lie bialgebra, and if $U_\hbar(\a)$
is a QUE algebra quantizing $\a$, then the corresponding quantized formal 
series Hopf (QFSH) algebra is 
$$
\CC[[A^*]]_\hbar := \{x\in U_\hbar(\a) | \forall n\geq 0, 
\delta_n(x) \in \hbar^n U_\hbar(\a)^{\wh\otimes n}\}. 
$$
Here $\delta_n = (\id - \eta \circ \eps)^{\otimes n} \circ \Delta^{(n)}$. 
Then $\CC[[A^*]]_\hbar$ is a flat deformation of $\CC[[A^*]] = U(\a^*)^*$, 
the function ring of the formal group with Lie algebra $\a^*$. 

Let $\l,\p$ be as in the beginning of this section. Then we have 
Lie algebra morphisms $\l^* \hookrightarrow \p^*$, $\p^*\twoheadrightarrow
\l^*$ such that $\l^* \hookrightarrow \p^* \twoheadrightarrow \l^*$ 
is the identity, so $\u^*$ is an ideal of $\l^*$.  
We are again 
in the situation of Section \ref{sect:quot}, with $A = \CC[[P^*]]_\hbar$, $B =
\CC[[L^*]]_\hbar$. Define $\CC[[U^*]]_\hbar$ to be the algebra 
$C$ defined in this
section.

\begin{proposition}
$\CC[[U^*]]_\hbar \subset \CC[[P^*]]_\hbar$ is a flat deformation of 
$\CC[[U^*]] \subset \CC[[P^*]]$. 
\end{proposition}

{\em Proof.} Similar to the proof of Proposition \ref{prop:iashar}, 
using Example \ref{ex:vayera}. 
\hfill \qed \medskip

\subsection{Polarized Lie bialgebras with a splitting}
\label{sect:3:3} \label{sect:PBW}

\begin{definition}
{\it $(\g,\delta,\l,\u)$ is a polarized Lie bialgebra with a splitting if 
$(\g,\delta)$ is a Lie bialgebra, such that $\l\subset \g$
is a Lie subbialgebra, $\u = \u_+ \oplus \u_-$, where $\u_\pm$
are Lie subalgebras of $\g$, and $\delta(\u_\pm) \subset \big( 
(\l\otimes \u_\pm)  
\oplus (\u_\pm \otimes \l) \big)  \oplus \wedge^2(\u_\pm)$. 
}\end{definition}

\begin{example}
See Example \ref{example:21}. 
\end{example}

Set $\p_\pm = \l\oplus \u_\pm$, then $\p_\pm$ are Lie subbialgebras of 
$(\g,\delta)$. 
Then we have Lie bialgebra morphisms $\l\subset \p_\pm \subset \g$
and $\p_\pm \twoheadrightarrow \l$, such that the composed map 
$\l \hookrightarrow \p_\pm \twoheadrightarrow \l$ is the identity. 
Applying a quantization functor to this situation, we get QUE algebras 
$U_\hbar(\l)$, $U_\hbar(\p_\pm)$ and $U_\hbar(\g)$ quantizing $\l$, 
$\p_\pm$ and $\g$, together with QUE algebra inclusions 
$U_\hbar(\l) \subset U_\hbar(\p_\pm) \subset U_\hbar(\g)$, and QUE algebra
morphisms $U_\hbar(\p_\pm) \stackrel{\pi_\pm}{\twoheadrightarrow} U_\hbar(\l)$, 
such that the composed map $U_\hbar(\l) \hookrightarrow U_\hbar(\p_\pm)
\stackrel{\pi_\pm}{\twoheadrightarrow} U_\hbar(\l)$  is the identity. 
In the situation of Example \ref{example:21}, $U_\hbar(\l) \subset 
U_\hbar(\p_+) \subset U_\hbar(\g)$ is as above, $U_\hbar(\p_-)$
is the subalgebra of $U_\hbar(\g)$ generated by the $h_i,f_i, i\in [1,r]$, 
and $e_i,i\in I$ and $\pi_\pm : U_\hbar(\p_\pm) \to U_\hbar(\l)$ are as above. 

Let us now set 
$$
U_\hbar(\u_\pm) := \{x\in U_\hbar(\p_\pm) | (\id\otimes \pi_\pm) 
\circ \Delta(x) = x \otimes 1\}, 
$$
then we have seen that $U_\hbar(\u_\pm)$
is a subalgebra and a left coideal of $U_\hbar(\p_\pm)$, and 
$U_\hbar(\u_\pm) \subset U_\hbar(\p_\pm)$ is a flat deformation of 
$U(\u_\pm) \subset U(\p_\pm)$. 

Then the tensor product of inclusions followed by the product
$U_\hbar(\u_+) \wh\otimes U_\hbar(\p_-) \to U_\hbar(\g)$ is an isomorphism
of $\CC[[\hbar]]$-modules. Therefore 
$$
\on{PBW} : U_\hbar(\u_+) \wh\otimes U_\hbar(\l) \wh\otimes U_\hbar(\u_-)
\to U_\hbar(\g)
$$
is an isomorphism of $\CC[[\hbar]]$-modules. We define $H : U_\hbar(\g)
\to U_\hbar(\l)$ as the composed map $(\eps \otimes \id \otimes \eps) \circ 
\on{PBW}^{-1}$. Then we have 
$$
\forall x_\pm\in U_\hbar(\u_\pm), \forall x_0 \in U_\hbar(\l), \quad 
H(x_+ x_0 x_-) = \eps(x_+) \eps(x_-) x_0
$$
and since $U_\hbar(\u_\pm)$ is stable under the adjoint action of 
$U_\hbar(\p_\pm)$ on itself, 
$$
\forall x_\pm \in U_\hbar(\p_\pm), \forall x \in U_\hbar(\g), \quad 
H(x_+ x) = \eps(x_+) H(x), \ H(xx_-) = \eps(x_-) H(x). 
$$

In the same way, we set 
$$
\CC[[U^*_\pm]]_\hbar := \{x\in \CC[[P^*_\pm]]_\hbar | (\id\otimes \pi_\pm) 
\circ \Delta(x) = x \otimes 1\} = U_\hbar(\u_\pm) \cap \CC[[P^*_\pm]], 
$$
then $\CC[[U^*_\pm]]_\hbar$
is a subalgebra and a left coideal of $\CC[[P^*_\pm]]_\hbar$, and 
$\CC[[U^*_\pm]]_\hbar \subset \CC[[P^*_\pm]]_\hbar$ is a flat deformation of 
$\CC[[U^*_\pm]] \subset \CC[[P^*_\pm]]$. Therefore $\on{PBW}$ restricts to 
an isomorphism of $\CC[[\hbar]]$-modules 
$$
\on{PBW} : \CC[[U_+^*]]_\hbar \bar\otimes \CC[[L^*]]_\hbar
\bar\otimes \CC[[U_-^*]])_\hbar \to \CC[[G^*]]_\hbar. 
$$
$H$ restricts to a map $\CC[[G^*]]_\hbar \to \CC[[L^*]]_\hbar$
with the same properties as above (replacing $U_\hbar(\x)$ by
$\CC[[X^*]]_\hbar$ everywhere).

\subsection{Localization of QFSH algebras}

Let $\a$ be a finite dimensional Lie bialgebra, let $U_\hbar(\a)$
be a QUE algebra quantizing $\a$ and let $\CC[[A^*]]_\hbar$ be the 
corresponding QFSH algebra.  Recall that $\CC[[A^*]]_\hbar$ is a 
flat deformation of $\CC[[A^*]] = U(\a^*)^*$. 

Let $D\in S^d(\a)$ be a nonzero element. We denote the 
localization of $\CC[[A^*]]$ w.r.t. $D$ by $\CC[[A^*]]_D$.  
This is a Poisson
algebra, equipped with left and right coproduct morphisms 
\begin{equation} \label{ext:copdt}
\CC[[A^*]]_D \to \CC[[A^*]] \bar\otimes \CC[[A^*]]_D, \; 
\CC[[A^*]]_D \to \CC[[A^*]]_D \bar\otimes \CC[[A^*]], 
\end{equation}
satisfying natural compatibility rules; in particular, $\CC[[A^*]]_D$
is a Poisson base algebra over $\CC[[A^*]]$ (see \cite{DM}, Section 4.1). 
Here $\bar\otimes$ has 
the following meaning: if $V$ is a vector space, then 
$V\bar\otimes \CC[[a_1,\ldots,a_d]] = \CC[[a_1,\ldots,a_d]] \bar\otimes V 
=  V[[a_1,\ldots,a_d]]$.  

\begin{proposition} \label{prop:3:9}
1) There is a unique quantization $\CC[[A^*]]_{D,\hbar}$ of $\CC[[A^*]]_D$, 
containing $\CC[[A^*]]_\hbar$. 

2) The coproduct morphism 
$\CC[[A^*]]_\hbar \to \CC[[A^*]]_\hbar^{\bar\otimes 2}$ 
extends uniquely to morphisms 
$\Delta_l : \CC[[A^*]]_{D,\hbar} \to \CC[[A^*]]_{D,\hbar} 
\bar\otimes \CC[[A^*]]_\hbar$ and 
$\Delta_r : \CC[[A^*]]_{D,\hbar} \to \CC[[A^*]]_\hbar
\bar\otimes \CC[[A^*]]_{D,\hbar}$, 
quantizing (\ref{ext:copdt}), satisfying natural axioms (see \cite{DM}). 
\end{proposition}

{\em Proof.}
In the same way as in \cite{ER}, Proposition 3.1, one proves the following
statement. Let $R$ be a Poisson ring and $S\subset R$ be a multiplicative part
of $S$. Then the localization $R_S$ has a unique Poisson structure, extending 
the Poisson structure of $R$. Let $R_\hbar$ be a quantization of $R$ 
(i.e., $R_\hbar$ is isomorphic 
to $R[[\hbar]]$ as a $\CC[[\hbar]]$-module; it is a topological algebra, with
classical limit $R$). Then there is a unique quantization $R_{S,\hbar}$ of 
$R_S$, containing $R_\hbar$. If $s\in S$ and $\bar s\in R_\hbar$ is a 
lift of $s$ in $R_\hbar$, then $\bar s$ is invertible in $R_{S,\hbar}$. 
The algebra $R_{S,\hbar}$ may be characterized by the following universal 
property: if 
$A_\hbar$ is a quasi-commutative algebra (i.e., $A_\hbar$ is a topological 
$\CC[[\hbar]]$-algebra, whose reduction modulo $\hbar$ is commutative), 
any morphism $R_\hbar \to A_\hbar$ of topological algebras, whose
reduction modulo $\hbar$ extends to a morphism $R_S \to A_\hbar/\hbar A_\hbar$,
extends in a unique way to a morphism $R_{S,\hbar} \to A_\hbar$ of topological
algebras. We call $R_{S,\hbar}$ the localization of $R_\hbar$ w.r.t. $S$. 
Then $\CC[[A^*]]_{D,\hbar}$ corresponds to $R = \CC[[A^*]]_\hbar$, 
$S = \{$powers of $D\}$. 

Let us now construct $\Delta_l,\Delta_r$. We have $\Delta(D) = 
D \otimes 1 + \sum_i a_i \otimes b_i$, where $b_i\in \m(A^*)_\hbar = 
\on{Ker}(\eps : \CC[[A^*]]_\hbar \to \CC[[\hbar]])$. Therefore $\Delta(D)$
is invertible in $\CC[[A^*]]_{D,\hbar} \bar\otimes \CC[[A^*]]_\hbar$, with 
image 
$$
D^{-1} \otimes 1  - \sum_i D^{-1} a_i D^{-1} \otimes b_i
+ \sum_{i,j} D^{-1} a_i D^{-1} a_j D^{-1} \otimes b_i b_j  - \cdots;  
$$
this sum is convergent since $b_{i_1}\cdots b_{i_n} \in \m(A^*)_\hbar^n$. 
Applying the universal property to $\Delta : \CC[[A^*]]_\hbar \to 
\CC[[A^*]]_\hbar^{\bar\otimes 2} \to \CC[[A^*]]_{D,\hbar} \bar\otimes
\CC[[A^*]]_\hbar$, we obtain the existence of $\Delta_l$. 
The morphism $\Delta_r$ is constructed in the same way. 
\hfill \qed \medskip 

In the case when $U_\hbar(\a) = U(\a)[[\hbar]]$ (so the Lie cobracket of $\a$
is zero) and $D\in \CC[[A^*]]$ is homogeneous, we now relate the algebra 
$\CC[[A^*]]_{D,\hbar}$ to the algebra
$\wh{U(\a)}$ introduced in  \cite{EE2}. $\wh{U(\a)}$ is the microlocalization 
of $U(\a)$ corresponding to $D$. It is a $\ZZ$-filtered algebra, whose 
associated graded is isomorphic to $S^\cdot(\a)_D$. 

\begin{lemma} \label{3:10}
Let $\cA$ be the subspace of $\wh{U(\a)}((\hbar))$ consisting in 
all series $\sum_{i\in \ZZ} \hbar^i x_i$, where $x_i$ vanishes for 
$i \ll 0$, and $\on{deg}(x_i) \leq i$ for any $i\in\ZZ$. Then $\cA$
is a subalgebra of $\wh{U(\a)}((\hbar))$, isomorphic to $\CC[[A^*]]_{D,\hbar}$.  
\end{lemma}

Here $X((\hbar)) = X[[\hbar]][\hbar^{-1}]$. 

\medskip

{\em Proof.} $\cA$ is a topologically free $\CC[[\hbar]]$-module and a
subalgebra of $\wh{U(\a)}((\hbar))$; one checks that $\cA / \hbar \cA$
is isomorphic to the degree completion of $S^\cdot(\a)_D$. 
$\CC[[A^*]]_\hbar \subset U(\a)[[\hbar]]$ identifies with the 
space of all series $\sum_{i\geq 0} \hbar^i x_i$, with $\on{deg}(x_i) \leq i$. 
So we have an injection $\CC[[A^*]]_\hbar \subset \cA$. One checks that a 
lift of $D$ in $\CC[[A^*]]_\hbar$ is taken to an invertible element of 
$\cA$. This shows that $\CC[[A^*]]_\hbar \subset \cA$ extends to an injection 
$\cO(A^*)_{D,\hbar} \subset \cA$, whose reduction modulo $\hbar$ is the
identity. This implies that $\CC[[A^*]]_{D,\hbar} \subset \cA$ is an 
isomorphism. 
\hfill \qed \medskip

\begin{remark}
We do not know a quantization of $\wh{U(\a)}$ without assuming 
$U_\hbar(\a) = U(\a)[[\hbar]]$ and $D$ homogeneous. 
\end{remark}

\subsection{Construction of $\Upsilon$}

Let $\g$ be a nondegenerate polarized Lie bialgebra with a splitting. 
Let $p\in S^d(\l)$ be as in Section \ref{nondeg}. Then there exists a unique 
$\wt\sigma_\l^\g \in \u_+ \otimes \u_- \otimes \CC[[L^*]]_{P_\l^\g}$, such that 
$$
\sigma_\l^\g = \wt\sigma_\l^\g - (\wt\sigma_\l^\g)^{2,1}. 
$$
Let $U_\hbar(\l)$, $U_\hbar(\u_\pm)$, $U_\hbar(\p_\pm)$, $U_\hbar(\g)$; 
$\CC[[L^*]]_\hbar$, $\CC[[U_\pm^*]]_\hbar$, $\CC[[P_\pm^*]]_\hbar$, 
$\CC[[G^*]]_\hbar$, and the morphisms between these algebras, 
be as in Section \ref{sect:q:lba}. 

\begin{theorem} \label{thm:Ups} \label{thm:3:12}
There exists a unique element $\Upsilon \in 
\big( U_\hbar(\u_+) \wh\otimes U_\hbar(\u_-)\big) \bar\otimes 
\CC[[L^*]]_{P_\l^\g,\hbar}$, which we denote $\Upsilon = \sum_i A_i \otimes B_i 
\otimes L_i$, such that for any $x_\pm \in U_\hbar(\u_\pm)$, we have 
$$
\sum_i H(x_- A_i) L_i H(B_i x_+) = H(x_- x_+). 
$$
Here we use the injection $U_\hbar(\l) \subset 
\CC[[L^*]]_{P_\l^\g,\hbar}[\hbar^{-1}]$. 
Then this identity also holds when $x_\pm \in U_\hbar(\p_\pm)$. 

$\Upsilon$ is a series $\sum_{n\geq 0} {{\hbar^n}\over {n!}}\Upsilon_n$, where 
$\Upsilon_n \in \big( U_\hbar(\u_+) \wh\otimes U_\hbar(\u_-)\big) \bar\otimes 
\CC[[L^*]]_{P_\l^\g,\hbar}$ is such that its reduction modulo $\hbar$ lies in 
$U(\u_+)_{\leq n} \otimes U(\u_-)_{\leq n} \otimes \CC[[L^*]]_{P_\l^\g}$, 
and the image of
this reduction in $S^n(\u_+) \otimes S^n(\u_-) \otimes\CC[[L^*]]_{P_\l^\g}$ 
is equal to $(\wt\sigma_\l^\g)^n$. 

Moreover, $\Upsilon$ (resp., $\Upsilon^{2,1,3}$) belongs to 
$U_\hbar(\u_\pm) \bar\otimes \CC[[U_\mp^*]]_\hbar \bar\otimes 
\CC[[L^*]]_{P_\l^\g,\hbar}$. 
\end{theorem}

If $V = E[[\hbar]]$ and $W = \CC[[a_1,\ldots,a_d]][1/p][[\hbar]]$, 
then $V\bar\otimes W = E[[a_1,\ldots,a_d]][1/p][[\hbar]]$.  

\medskip 
{\em Proof.} We first prove some properties of 
$H : \CC[[G^*]]_\hbar \to \CC[[L^*]]_\hbar$. 
Let $\m(U_\pm^*)_\hbar \subset \CC[[U_\pm^*]]_\hbar$
be the kernel of the composed map $\CC[[U_\pm^*]] \to \CC[[\hbar]]
\to \CC$. Let $\m(U_\pm)$ be the kernel of $\CC[[U_\pm^*]] \to \CC$. 
Then $\m(U_\pm^*) \hookrightarrow \m(U_\pm^*)_\hbar / \hbar \m(U_\pm^*)_\hbar$. 

\begin{lemma} \label{lemma:H}
Define a bilinear map $(-,-) : \CC[[U_-^*]]_\hbar 
\times \CC[[U_+^*]]_\hbar \to \CC[[L^*]]_\hbar$ by $(f_-,f_+) := H(f_- f_+)$.
Then if $f_\pm \in \m(U_\pm^*)_\hbar^{n_\pm}$, we have 
$(f_-,f_+) = O(\hbar^{\on{max}(n_+,n_-)})$. 

Thus, $\hbar^{-n}(-,-)$ restricts to a bilinear map $\m(U_-^*)_\hbar^n \times
\m(U_+^*)_\hbar^n \to \CC[[L^*]]_\hbar$. This map induces a 
bilinear map $\m(U_-^*)^n \times \m(U_+^*)^n \to \CC[[L^*]]$, 
which factors through a bilinear map $(-,-)_n : S^n(\u_-) \times S^n(\u_+)
\to \CC[[L^*]]$ (as $S^n(\u_\pm) = \m(U_\pm^*)^n / \m(U_\pm^*)^{n+1}$). 
$(-,-)_n$ is the $n$th symmetric power of the bilinear map
$(-,-)_1$, taking $(u_-,u_+)\in \u_- \times \u_+$ to the formal 
function $\ell^* \mapsto \omega_{\ell^*}(u_+,u_-)$ (an element of $\CC[[L^*]]$). 
\end{lemma}

{\em Proof of Lemma.} We have $\g^* = \l^* \oplus \u_+^* \oplus \u_-^*$, 
where $\l^*$ is a Lie subalgebra, $[\l^*,\u_\pm^*] \subset \u_\pm^*$, 
$[\u_+^*,\u_-^*] = 0$. We have $\p_\pm^* = \l^* \oplus \u_\pm^*$, 
and the natural projections $\g^* \to \p_\pm^*$ are Lie algebra morphisms. 

Let compute the classical limit of the maps $\on{PBW}$ and $H$ 
(relative to QFS algebras). The classical limit of $\on{PBW}$
is an algebra morphism $\on{PBW}_0 : \CC[[U_+^*]] \bar\otimes \CC[[L^*]] 
\bar\otimes 
\CC[[U^*_-]] \to \CC[[G^*]]$, induced by the inclusions 
$\CC[[U_\pm^*]] = \CC[[P_\pm^*/L^*]] \subset \CC[[P_\pm^*]] \subset 
\CC[[G^*]]$
followed by multiplication. It follows that $\on{PBW}$ is the 
dual of the coalgebra isomorphism $U(\g^*) \to U(\u_+^*) \otimes 
U(\l^*) \otimes U(\u_-^*)$, whose inverse takes $x_+ \otimes x_0 \otimes x_-$
to $x_+ x_- x_0$. Therefore the classical limit of $H$ is the algebra 
morphism $H_0 : \CC[[G^*]] \to \CC[[L^*]]$
induced by the restriction ($\l^*$ is a Lie subalgebra of $\g^*$, thus 
$L^*$ is a formal subgroup of $G^*$). 

Let us now prove the Lemma. 
Let $f_\pm\in \m(U^*_\pm)_\hbar$ and $\bar f_\pm\in \m(U^*_\pm)$ 
be the image of $f_\pm$ in $\CC[[U^*_\pm]]$. Then $H(f_- f_+) = H(f_+ f_-) 
+ H([f_-,f_+])$. We have $H(f_-f_+) = O(\hbar^2)$. On the other 
hand, $(\on{mod}\ \hbar)(\hbar^{-1} [f_-,f_+]) = \{\bar f_-,\bar f_+\}$. 
Therefore, $H(f_-f_+) = \hbar H_0(\{\bar f_-,\bar f_+\}) + O(\hbar^2)$. 
This proves the first statement when $n_\pm = 1$.  

In the general case, the first statement of the Lemma follows from a 
computation based on the fact that $\CC[[G^*]]_\hbar$ is quasicommutative. 
The same computation shows that if $f_{\pm,i} \in \m(U^*_\pm)_\hbar$, 
$i \in [1,r]$, then 
$$
(f_{-,1}\cdots f_{-,n}, f_{+,1}\cdots f_{+,n}) = \sum_{\sigma\in\SG_n} 
\prod_{i=1}^n (f_{-,i},f_{+,\sigma(i)}) + O(\hbar^{1-n}). 
$$

It remains to compute the map $\m(U^*_-) \otimes \m(U^*_+) \to \CC[[L^*]]$, 
$\bar f_- \otimes \bar f_+ \mapsto H_0(\{\bar f_-,\bar f_+\})$. 
This map clearly factors through $(\m(U_-^*)/\m(U_-^*)^2) \otimes 
(\m(U_+^*)/\m(U_+^*)^2) = \u_- \otimes \u_+$. Let $u_\pm \in \u_\pm$
be the images of $\bar f_\pm$. Let $x\in \l^*$, then 
$$
\langle \{\bar f_-,\bar f_+\} , x^n \rangle = 
\langle \bar f_- \otimes \bar f_+, \delta(x^n) \rangle . 
$$
Now $\langle \bar f_\pm, U(\g^*)\l^* \rangle =0$, so the only 
contribution is that of 
$\langle \bar f_- \otimes \bar f_+, \Delta_0(x)^{n-1}\delta(x) \rangle$, 
i.e., $\langle \bar f_- \otimes \bar f_+, 
\ad(x \otimes 1 + 1 \otimes x)^{n-1}(\delta(x))\rangle$.
Let $x$ be formal near $0$ in $\l^*$, then 
\begin{align*}
& \{\bar f_-,\bar f_+\}(e^x) = \sum_{n\geq 0}{1\over {n!}} 
\langle u_- \otimes u_+, \ad(x \otimes 1 + 1 \otimes x)^{n-1}(\delta(x))
\rangle
\\ & = \sum_{k,l\geq 0} {1\over{(k+l+1) k! l!}} \langle \big[ 
\big( \ad(-x)^k(u_+)\big)_{\g}, \big( \ad(-x)^l(u_-)\big)_{\g}
\big] , x\rangle
\\ & = 
\sum_{k,l\geq 0} {{(-1)^{l+1}}\over{(k+l+1) k! l!}}
\langle u_+, \on{ad}(x)^k \circ p_{\g^*} \circ \on{ad}(x) \circ p_\g \circ 
\on{ad}(x)^l (u_-) \rangle,  
\end{align*}
where $p_\g,p_{\g^*}$ are the projections of $\d(\g) = \g \oplus \g^*$
on $\g$ and $\g^*$. 

On the other hand, 
\begin{align*}
\omega_{e^x}(u_+,u_-) & = \sum_{k,l\geq 0} {{(-1)^{k+l}}\over{k!l!}}
\langle p_\g \circ \on{ad}(x)^k(u_+), p_{\g^*} \circ \on{ad}(x)^l(u_-)\rangle 
\\ & 
= \sum_{k,l\geq 0} {{(-1)^{l}}\over{k!l!}}
\langle u_+, \on{ad}(x)^k \circ p_{\g^*} \circ \on{ad}(x)^l(u_-)\rangle . 
\end{align*}

We have $\on{ad}(x)(\g) \subset \g\oplus \g^*$, $\on{ad}(x)(\g^*) \subset \g^*$,
therefore 
$$
p_{\g^*} \circ \on{ad}(x)^l(u_-) = \sum_{\alpha = 0}^{l-1}
\on{ad}(x)^\alpha \circ p_{\g^*} \circ \on{ad}(x) \circ p_{\g}
\circ \on{ad}(x)^{l-1-\alpha}(u_-). 
$$
So 
\begin{align*}
& \omega_{e^x}(u_+,u_-) = \sum_{k,l,\alpha\geq 0}
{{(-1)^l}\over{k!l!}} \langle u_+, \on{ad}(x)^{k+\alpha} \circ p_{\g^*}
\circ \on{ad}(x) \circ p_{\g} \circ \on{ad}(x)^{l-1-\alpha}(u_-)\rangle
\\ &
= \sum_{k',l'\geq 0} 
\big( \sum_{\alpha\geq 0} {{(-1)^{l'+\alpha+1}}\over
{(k'-\alpha)!(l'+\alpha+1)!}}\big) 
\langle u_+, \on{ad}(x)^{k'} \circ \on{ad}(x) \circ 
p_\g \circ \on{ad}(x)^{l'} (u_-)\rangle   . 
\end{align*}
Now 
\begin{align*}
\sum_{\alpha\geq 0} 
{{ (-1)^{l'+\alpha +1} }\over{(k'-\alpha)! (l'+\alpha+1)!}}
& = 
{{(-1)^{l'+1}}\over{(k'+l'+1)!}} 
\big( \pmatrix k'+l'+1 \\ l'+1\endpmatrix 
- \pmatrix k'+l'+1 \\ l'+2\endpmatrix + \cdots 
+ (-1)^{k'} \pmatrix k'+l'+1 \\ k'+l'+1\endpmatrix 
 \big) 
\\ & =  
{{(-1)^{k'+l'+1}}\over{(k'+l'+1)!}} 
\big( \pmatrix k'+l'+1 \\ 0\endpmatrix 
- \pmatrix k'+l'+1 \\ 1\endpmatrix + \cdots 
+ (-1)^{k'} \pmatrix k'+l'+1 \\ k'\endpmatrix 
 \big) 
\\ & =  
{{(-1)^{k'+l'+1}}\over{(k'+l'+1)!}} (-1)^{k'}
\pmatrix k'+l' \\ k'\endpmatrix 
= {{(-1)^{l'+1}}\over{(k'+l'+1) k'! l'!}}. 
\end{align*}
It follows that $\{\bar f_-,\bar f_+\}(e^x) = \omega_{e^x}(u_+,u_-)$. 

This ends the proof of the Lemma. 
\hfill \qed \medskip 

\begin{lemma}
The bilinear map $(-,-)$ of Lemma \ref{lemma:H} extends to a bilinear map 
$(-,-) : \CC[[U_-^*]]_\hbar \times U_\hbar(\u_+) \to \CC[[L^*]]_\hbar$, 
whose reduction modulo $\hbar$ is a bilinear map
$(-,-)_0 : \CC[[U_-^*]] \times U(\u_+) \to \CC[[L^*]]$, uniquely defined 
by the conditions that $(fg,x)_0 = \sum (f,x^{(1)})_0(g,x^{(2)})_0$ for any 
$x\in U(\u_+)$ and $f,g\in\CC[[U_-^*]]$, and $(f,u_+)_0$ is the formal function 
$\ell^* \mapsto \omega_{\ell^*}(u_+,\bar f)$ if $f\in \m(U_-^*)$ and
$\bar f$ is its image in $\u_-$.   

In particular, we have $(f,x)_0 = 0$ if $\on{val}(x) > \on{deg}(f)$. 
When $\on{val}(x) = \on{deg}(f) = n$, the pairing $(-,-)_0$
factors through a bilinear
map $S^n(\u_-) \times S^n(\u_+) \to \CC[[L^*]]$, which coincides with the 
map described in Lemma \ref{lemma:H}.
\end{lemma}

{\em Proof.} The first statement follows from the fact that 
$\CC[[U_-^*]]_\hbar$ is the $\hbar$-adic completion of 
$\sum_{n\geq 0} \hbar^{-n} \m(U_-^*)_\hbar$. The rest follows from 
the formula 
$$
(f,u_1\cdots u_n)_0 = \{\{f,f_{u_1}\},\ldots, f_{u_n}\}_{|L^*}, 
$$
where $f_u\in \m(U_+^*)$ is a lift of $u\in\u_+ = \m(U_+^*) / \m(U_+^*)^2$ 
and the results of Lemma \ref{lemma:H}. \hfill \qed \medskip 

\medskip 
{\em End of proof of Theorem \ref{thm:Ups}.} 
The properties of $(-,-)_0$, together with the nondegeneracy 
of $\omega_{\ell^*}(u_+,u_-)$, imply that there exists a unique element 
$\Upsilon_0 \in U(\u_+) \bar\otimes \CC[[L^*]]_{P_\l^\g} 
\bar\otimes \CC[[U_-^*]]$, 
which we write $\Upsilon_0 = \sum_i A_{i,0} \otimes L_{i,0} \otimes B_{i,0}$,
such that for any $f_0\in \CC[[U_-^*]]$, $x_0\in U(\u_+)$, we have 
$$
\sum_i (f_0,A_{i,0})_0 L_{i,0} (B_{i,0},x_0)_0 = (f_0,x_0)_0. 
$$
$\Upsilon_0$ satisfies the relation $\Upsilon_0^{12,3,4} = 
\Upsilon_0^{1,3,4} \Upsilon_0^{2,3,4}$, which implies that it has 
the form $\Upsilon_0 = \on{exp}(\upsilon_0)$, with $\upsilon_0
\in \u_+ \otimes \CC[[L^*]]_{P_\l^\g} \bar\otimes \m(U_-^*)$; moreover, 
the properties of $(-,-)_0$ imply that the projection of 
$\upsilon_0$ in $\u_+ \otimes \CC[[L^*]]_{P_\l^\g} \otimes \u_-$ is 
$(\wt\sigma_\l^\g)^{1,3,2}$. 
 
Let $\Upsilon_1$ be a lift of $\Upsilon_0$ to 
$U_\hbar(\u_+) \bar\otimes \CC[[L^*]]_{P_\l^\g,\hbar} 
\bar\otimes \CC[[U_-^*]]_\hbar$, 
which we write $\Upsilon_1 = \sum_i A_{i,1} \otimes L_{i,1} \otimes B_{i,1}$.
Define $\varphi : \CC[[U_-^*]]_\hbar \to \CC[[L^*]]_{P_\l^\g,\hbar} 
\bar\otimes 
\CC[[U_-^*]]_\hbar$ by 
$$
\varphi(f) := \sum_i (f,A_{i,1}) L_{i,1} \otimes B_{i,1}. 
$$ 
Then the reduction modulo $\hbar$ of $\varphi$ is 
$\CC[[U_-^*]] \to \CC[[L^*]]_{P_\l^\g} \bar\otimes \CC[[U_-^*]]$, 
$f_0\mapsto f_0 \otimes 1$. 

$\varphi$ gives rise to an endomorphism 
$\wt\varphi$ of $\CC[[L^*]]_{P_\l^\g,\hbar} \bar\otimes \CC[[U_-^*]]_\hbar$,
defined uniquely by  $\wt\varphi(\ell\otimes f) 
= (\ell \otimes 1) \varphi(f)$
for $\ell \in \CC[[L^*]]_{P_\l^\g,\hbar}$, $f\in \CC[[U_-^*]]_\hbar$. 
The reduction modulo $\hbar$ of $\wt\varphi$ is the identity, 
so $\varphi$ is invertible. Let $\psi$ be its inverse. Define $\Upsilon$
by  
$$
\Upsilon^{1,3,2} := \sum_i (A_{i,1} \otimes L_{i,1} \otimes 1)
(1\otimes \psi(1\otimes B_{i,1})). 
$$
We write $\Upsilon = \sum_i A_i \otimes L_i \otimes B_i$. Then 
if $f \in \CC[[U_-^*]]_\hbar$, we have 
\begin{align*}
& \sum_i (f,A_i)L_i \otimes B_i = \sum_i 
\big( (f,A_{i,1}) L_{i,1} \otimes 1 \big) \psi(1\otimes B_{i,1})
=  \sum (\varphi(f)^{(1)} \otimes 1)\psi(1\otimes \varphi(f)^{(2)}) 
\\ & = \psi\circ\varphi(f)
= \psi \circ \wt\varphi(1\otimes f) = 1\otimes f.   
\end{align*}
Here we write $\varphi(f) = \sum \varphi(f)^{(1)} \otimes \varphi(f)^{(2)}$; 
the third equality follows from the fact that $\psi$ satisfies 
$\psi(\ell \otimes f) = (\ell \otimes 1)\psi(1\otimes f)$. 
The other properties of $\Upsilon$ follow from the fact that 
$(\on{mod\ }\hbar)(\psi)$ is the identity.  
  
Finally, the fact that $\Upsilon^{2,1,3}$ is in $U_\hbar(\u_-)
\bar\otimes \CC[[U^*_+]] \bar\otimes \CC[[L^*]]_{P_\l^\g,\hbar}$ is proven 
by exchanging the indices $+,-$.  
\hfill \qed \medskip

\subsection{The dynamical twist $\Psi_\l^\g$}

Set now 
$$
\Psi_\l^\g = \big( 
\sum_i A_i \otimes B_i S(L_i^{(2)}) \otimes S(L_i^{(1)}) \big)^{-1}.  
$$

Then $\Psi_\l^\g$ lies in $U_\hbar(\u_+) \wh\otimes 
\big( \CC[[P_-^*]]_\hbar \bar\otimes \CC[[L^*]]_{P_\l^\g,\hbar}\big) $, 
therefore in $\big( U_\hbar(\u_+) \wh\otimes U_\hbar(\p_-)\big) 
\bar\otimes \CC[[L^*]]_{P_\l^\g,\hbar}$.  

\begin{theorem} \label{thm:Psi} \label{thm:3:15}
$\Psi_\l^\g$ satisfies the dynamical twist equation 
\begin{equation} \label{dyn:twist}
(\Psi_\l^\g)^{1,2,34} (\Psi_\l^\g)^{12,3,4} = 
(\Psi_\l^\g)^{2,3,4} (\Psi_\l^\g)^{1,23,4}. 
\end{equation}
Moreover $\Psi_\l^\g = 1 + O(\hbar)$, and $(\on{mod}\
\hbar)(\hbar^{-1}(\Psi_\l^\g - 1)) = \wt\sigma_\l^\g$. So $\Psi_\l^\g$
is a quantization of $\sigma_\l^\g$. 
\end{theorem}

Here $(\Psi_\l^\g)^{1,2,34} := (\id\otimes \id\otimes \Delta_r)(\Psi_\l^\g)$, 
see Proposition \ref{prop:3:9}. 

\medskip 

{\em Proof.} Similar to the proof of Proposition 2.5 in \cite{EE2}.   
\hfill \qed \medskip 

\begin{remark}
When the cobracket of $\l$ is zero, $\l^* = L^*$ and $P_\l^\g$ defined by 
(\ref{def:P}) is a homogeneous element of $\CC[[L^*]] = \wh S^\cdot(\l)$. 
So Lemma \ref{3:10} implies that the construction of Theorems \ref{thm:3:12}, 
\ref{thm:3:15} are generalizations of \cite{EE2}, Propositions 0.1 and 1.1. 
Theorem \ref{thm:Psi} may also be viewed as a generalization of the 
construction of \cite{EV2} to the case of a nonabelian base. 
\end{remark}

\subsection{Localized Harish-Chandra maps}

Let $\cO(X)_\hbar$ be QFS algebra, i.e., a flat deformation of a formal series
algebra $\cO(X)$. To each $D\in \cO(X)$, we associate the completed 
localization $\cO(X)_{D,\hbar}$. 

\begin{proposition} 
Let $\bar D\in\cO(X)_\hbar$ be a lift of $D$. For each $\alpha >0$, 
the subspace $\cO_\alpha(\bar D) := \{\sum_{i\geq 0} \hbar^i x_i | 
x_i \in \bar D^{-i-\alpha} \cO(X)_\hbar\} \subset \cO(X)_{D,\hbar}$
is independent of a choice of $\bar D$. We denote it $\cO_\alpha$. 
Set $\cO_0 := \cO(X)_\hbar$. Then we have $\cO_0 \subset \cO_1 \subset \cO_2
\subset \ldots$, and $\cO_\alpha \cO_\beta \subset \cO_{\alpha + \beta}$, 
for any $\alpha,\beta \geq 0$.  
\end{proposition} 

{\em Proof.} For $a,b\in\cO(X)_\hbar$, set $\on{ad}_\hbar(a)(b)
= \hbar^{-1} [a,b]$. Then $\on{ad}_\hbar(a)(b) \in \cO(X)_\hbar$. 
Then we have: 
\begin{equation} \label{useful:id}
\forall a\in \cO(X)_\hbar, 
\; 
a \bar D^{-1} = \bar D^{-1} a - \hbar \bar D^{-2} \on{ad}_\hbar(\bar D)(a)
+ \hbar^2 \bar D^{-3} \on{ad}_\hbar(\bar D)^2(a) - \cdots. 
\end{equation} 
Using this identity, one shows that $\cO_\alpha(\bar D) \cO_\beta(\bar D)
\subset \cO_{\alpha + \beta}(\bar D)$ for any $\alpha,\beta$. 

Let us now prove the independence statement. Let $\wt D$ be another lift of $D$.
Then $\wt D = \bar D + \hbar a$, where $a\in \cO(X)_\hbar$. Then: 
$$
\wt D^{-1} = \bar D^{-1} - \hbar \bar D^{-1} a \bar D^{-1} 
+ \hbar^2 \bar D^{-1} a \bar D^{-1} a \bar D^{-1} - \cdots.   
$$
Using (\ref{useful:id}), we see that $\wt D^{-1} \in \cO_1(\bar D)$. In the same
way, one shows that if $x\in \cO(X)_\hbar$, then $\hbar^i \wt D^{-1-i}x
\in \cO_1(\bar D)$, therefore $\cO_1(\wt D) \subset \cO_1(\bar D)$. 
So $\cO_1(\wt D) = \cO_1(\bar D)$. In the same way, one shows that
$\cO_\alpha(\wt D) = \cO_\alpha(\bar D)$. 
\hfill \qed \medskip

Let $\p = \l\oplus \u$ be a Lie bialgebra with a splitting. Then we have a Hopf
algebra inclusion $\CC[[L^*]] \subset \CC[[P^*]]$ and an algebra inclusion
$\CC[[U^*]] \subset \CC[[P^*]]$. Let $\Pi\in \CC[[L^*]]$ be a nonzero element. 

Let $U_\hbar(\l) \hookrightarrow U_\hbar(\p) \twoheadrightarrow U_\hbar(\l)$ be
a quantization of $\l\hookrightarrow \p \twoheadrightarrow \l$. It gives rise to
a QFSH algebra inclusion $\CC[[L^*]]_\hbar \subset \CC[[P^*]]_\hbar$ and a 
left coideal inclusion $\CC[[U^*]]_\hbar \subset \CC[[P^*]]_\hbar$. The tensor
product of inclusions followed by the product yields isomorphisms of topological
vector spaces
$$
\CC[[L^*]]_\hbar \bar\otimes \CC[[U^*]]_\hbar \to \CC[[P^*]]_\hbar, \; 
\CC[[U^*]]_\hbar \bar\otimes \CC[[L^*]]_\hbar \to \CC[[P^*]]_\hbar. 
$$ 

\begin{proposition}
The tensor product of inclusions followed by the product also 
gives rise to isomorphisms 
\begin{equation} \label{new:maps}
\alpha : \CC[[L^*]]_{\Pi,\hbar} \bar\otimes \CC[[U^*]]_\hbar \to 
\CC[[P^*]]_{\Pi,\hbar},\; 
\beta : \CC[[U^*]]_\hbar \bar\otimes \CC[[L^*]]_{\Pi,\hbar} \to 
\CC[[P^*]]_{\Pi,\hbar}. 
\end{equation}
\end{proposition}

{\em Proof.} The inclusion $\CC[[L^*]]_\hbar \to \CC[[P^*]]_{\Pi,\hbar}$
obviously extends to an inclusion $\CC[[L^*]]_{\Pi,\hbar} \to 
\CC[[P^*]]_{\Pi,\hbar}$. 
Then the reductions modulo $\hbar$ of the maps (\ref{new:maps}) 
coincide with the product map $\CC[[L^*]]_\Pi \bar\otimes \CC[[U^*]] \to 
\CC[[P^*]]_\Pi$, 
and with its composition with the permutation, which are isomorphisms. 
Therefore the maps (\ref{new:maps}) are isomorphisms. 
\hfill \qed \medskip 

Let now $\g = \l\oplus \u_+ \oplus \u_-$ be a polarized Lie bialgebra with a
splitting. Let $\Pi\in \CC[[L^*]]$ be a nonzero element and let 
$U_\hbar(\l) \subset U_\hbar(\p_\pm) \subset U_\hbar(\g)$, 
$U_\hbar(\p_\pm) \to U_\hbar(\l)$ be quantizations of the natural Lie bialgebra
morphisms. Define a product on $\CC[[L^*]]_{\Pi,\hbar} \bar\otimes 
\big( \CC[[U_+^*]]_\hbar \bar\otimes \CC[[U_-^*]]_\hbar \big)$ by 
\begin{align*}
& \mu = (132) \circ (m_{\CC[[U_+^*]]_\hbar} \otimes 
m^{(3)}_{\CC[[L^*]]_{\Pi,\hbar}}
\otimes m_{\CC[[U_-^*]]_\hbar}) \circ 
(\id\otimes e_+ \otimes \id \otimes e_-^{-1}
\otimes \id) 
\\ & \circ (\id\otimes \id \otimes \pi \otimes \id\otimes \id) 
\circ ((132)\otimes (132)), 
\end{align*} 
where $m_A$ is the product map of an algebra $A$, $m_A^{(3)} = 
(m_A \otimes \id) \circ m_A$, $e_\pm : \CC[[L^*]]_{\Pi,\hbar} \bar\otimes 
\CC[[U_\pm^*]]_\hbar \to \CC[[U_\pm^*]]_\hbar \bar\otimes 
\CC[[L^*]]_{\Pi,\hbar}$ 
is the composition $\beta_\pm^{-1} \circ \alpha_\pm$ (see (\ref{new:maps})), 
$$
\pi : \CC[[U^*_-]]_\hbar \bar\otimes \CC[[U^*_+]]_\hbar \to 
\CC[[U_+^*]]_\hbar \bar\otimes \CC[[L^*]]_\hbar \bar\otimes \CC[[U_-^*]]_\hbar
$$
is the composition of the product $\CC[[U_-^*]]_\hbar \bar\otimes 
\CC[[U_+^*]]_\hbar \to \CC[[G^*]]_\hbar$ with the inverse of 
$\on{PBW}$ (see Section \ref{sect:PBW}). 

\begin{proposition} 
$\mu$ is an associative product, extending the transport of the product on
$\CC[[G^*]]_\hbar$ by $\on{PBW}$. Let $\Pi'\in \CC[[G^*]]$ 
be a nonzero element, such
that $\Pi'_{|L^*} = \Pi$. Then there is a unique algebra morphism 
 $\wh{\on{PBW}^{-1}} : \CC[[G^*]]_{\Pi',\hbar} \to \big( \CC[[L^*]]_{\Pi,\hbar}
 \bar\otimes (\CC[[U_+^*]] \bar\otimes \CC[[U_-^*]]_\hbar),\mu\big)$, extending
 $\on{PBW}^{-1}$. 
\end{proposition}

{\em Proof.} The maps $\alpha_\pm,\beta_\pm$ satisfy identities implying the
associativity of $\mu$. The reduction of $\mu$ modulo $\hbar$ is the standard
product on $\CC[[L^*]]_P \bar\otimes \CC[[U_+^*]] \bar\otimes \CC[[U_-^*]]$. 
So it remains to prove that the image of $\Pi'\in\CC[[G^*]]$ in this algebra
is invertible. This image is $\Pi \otimes 1 \otimes 1 + \sum_i 
a_i \otimes b_i \otimes c_i$, where $b_i \otimes c_i$ belongs to the maximal
ideal $\m(U_+^* \times U_-^*)$ of $\CC[[U_+^*]] \wh\otimes \CC[[U_-^*]]$. 
Then the inverse of this image is $\sum_{\alpha\geq 0} (-1)^\alpha
(\Pi^{-1}\otimes 1 \otimes 1)^{\alpha + 1} (\sum_i a_i \otimes b_i \otimes
c_i)^\alpha$; the term corresponding to the index $\alpha$ belongs to 
$\m(U_+^* \otimes U_-^*)^\alpha$, so the series converges. 
\hfill \qed \medskip

\begin{remark}
Set $H := (\eps\otimes \id \otimes \eps) \circ \on{PBW}^{-1}$, then $H$ is a
flat deformation of the restriction map $\CC[[G^*]] \to \CC[[L^*]]$. 
$H$ is also an
analogue of the Harish-Chandra map. Let us define $\wh H$ as the composed map 
$\CC[[G^*]]_{\Pi',\hbar} \to \CC[[L^*]]_{\Pi,\hbar} \bar\otimes  \big( 
\CC[[U_+^*]]_\hbar \bar\otimes \CC[[U_-^*]]_\hbar \big) \stackrel{\id\otimes
\eps\otimes \eps}{\to} \CC[[L^*]]_{\Pi,\hbar}$, then $\wh H$ is an extension of 
$H$. We will denote $H,\wh H$ by $H_\l^\g$, $\wh H_\l^\g$ to avoid confusions.   
\end{remark}

\subsection{The quantum composition formula} \label{sect:3:8}

We now prove a Poisson-Lie analogue of Proposition 2.15 in \cite{EE2}. 
Namely, we prove a composition formula for the twists $\Psi_\l^\g$.  
It may be viewed as a quantization of Remark \ref{rem:comp}. We will use it to 
determine $\Psi_\l^\g$ using 
$\Psi_\h^\g$ and $\Psi_\h^\l$ in the situation of 
Example \ref{example:21}. 

Assume that $\g = \l\oplus \u_+ \oplus \u_-$ and $\l = \k \oplus \m_+ \oplus
\m_-$ are nondegenerate polarized Lie bialgebras. Set 
$\vv_\pm = \u_\pm \oplus \m_\pm$. 
Then $\g = \k\oplus \vv_+ \oplus \vv_-$ is a nondegenerate polarized Lie
bialgebra (the nondegeneracy follows from \cite{EE2}, Lemma 2.13).
We fix compatible quantizations of these polarized Lie bialgebras. 

We denote by $P_\l^\g\in\CC[[L^*]]$ and $P_\k^\l,P_\k^\g\in \CC[[K^*]]$
the elements associated to the decompositions $\g = \l \oplus \u_+ \oplus \u_-$,
$\l = \k \oplus \m_+ \oplus \m_-$, $\g = \k \oplus \vv_+ \oplus \vv_-$
according to (\ref{def:P}). Then $P_\k^\g = (P_\l^\g)_{|K^*} P_\k^\l$. It
follows that we have natural injections 
$$ 
\CC[[K^*]]_{(P_\l^\g)_{|K^*},\hbar} \hookrightarrow \CC[[K^*]]_{P_\k^\g,\hbar}
\;\on{and}\;
\CC[[K^*]]_{P_\k^\l,\hbar} \hookrightarrow \CC[[K^*]]_{P_\k^\g,\hbar}
$$
Define the linear map 
$$
\eta : \big( U_\hbar(\u_+) \wh\otimes U_\hbar(\p_-)\big) 
\bar\otimes \CC[[L^*]]_{P_\l^\g,\hbar}
\to 
\big( U_\hbar(\vv_+) \wh\otimes U_\hbar(\p_-)\big) 
\bar\otimes \CC[[K^*]]_{P_\k^\g,\hbar}
$$
by 
$\eta(\alpha \otimes \beta \otimes \lambda) = 
\sum_i \sum \alpha S(\lambda_{+,i}^{(2)}) \otimes \beta S(\lambda_{+,i}^{(1)})
\otimes \eps(\lambda_{-,i}) \lambda_{0,i}$, where $\wh{\on{PBW}^{-1}}(\lambda) 
= \sum_i \lambda_{+,i} \otimes \lambda_{0,i} \otimes \lambda_{-,i}$, 
with $\lambda_{\pm,i}
\in \CC[[M_\pm^*]]_\hbar$ and $\lambda_{0,i} \in 
\CC[[K^*]]_{(P_\l^\g)_{|K^*},\hbar} \hookrightarrow \CC[[K^*]]_{P_\k^\g,\hbar}$, 
so that $\lambda = \sum_i \lambda_{+,i} \lambda_{0,i} \lambda_{-,i}$. 

\begin{proposition} \label{quantum:comp}
We have 
$$
(\Psi_\k^\g)^{-1} = \eta((\Psi_\l^\g)^{-1})(\Psi_\k^\l)^{-1}
$$
(equality in $\big( U_\hbar(\vv_+) \wh\otimes U_\hbar(\p_-)\big) 
\bar\otimes \CC[[K^*]]_{P_\k^\l,\hbar}$).  
\end{proposition}

{\em Proof.} Similar to the proof of \cite{EE2}, Proposition 2.15. \hfill \qed
\medskip 

\begin{remark}
As in \cite{EE2}, Remark 2.16, this formula enables one to recover $\Psi_\l^\g$
from $\Psi_\k^\l,\Psi_\k^\g$. Indeed, let 
$$
\eta' : \big( U_\hbar(\vv_+) \wh\otimes U_\hbar(\p_-)\big) 
\bar\otimes \CC[[K^*]]_{P_\k^\g,\hbar} \to 
\big( U_\hbar(\u_+) \wh\otimes U_\hbar(\p_-)\big) 
\bar\otimes \CC[[K^*]]_{P_\k^\g,\hbar}[\hbar^{-1}]
$$
be the linear map taking 
$u_+ \lambda_+ \otimes p_- \otimes k$ to 
$u_+ \otimes p_- S(\lambda_+) \otimes k$, where $u_+\in U_\hbar(\u_+)$, 
$\lambda_+ \in U_\hbar(\m_+)$, 
$p_-\in U_\hbar(\p_-)$, $k\in \CC[[K^*]]_{P_\k^\g,\hbar}$. Then 
$\eta' \circ \eta = \id\otimes \id \otimes H_\k^\l$. 
Now $\Psi_\l^\g$ may be recovered uniquely from its 
image $\id\otimes \id\otimes H_\k^\l$ using its $\l$-invariance, as
in \cite{EE2}. 
\end{remark}

\section{Poisson homogeneous structures on $G/L$} \label{quant:homog}
\label{sect:4}

We first complement Section \ref{sect:1} by studying reductions 
of the (quasi)Poisson
structures associated with dynamical $r$-matrices. 

We then restrict ourselves to the case when $G$ is simple
and $L\subset G$ is a Levi subgroup. The Poisson
homogeneous structures on $G/L$ were classified in \cite{DGS}. The set of all
these Poisson structures is an algebraic variety $\cP$, 
which we study in Section
\ref{sect:DGS}. We introduce a Zariski open subset $\cP_0 \subset \cP$
and show that the Poisson structures corresponding to elements of $\cP_0$ 
can be expressed using PL dynamical $r$-matrices, and can therefore be
viewed as examples of the construction of Section \ref{sect:basic}
(Section \ref{sect:4:2}).

\subsection{Reductions of Poisson structures}
\label{sect:basic}

Let $\l\subset\g$ be an inclusion of finite dimensional Lie bialgebras. 
Let $L^*_{\on{alg}}$ be a Poisson-Lie 
algebraic group with Lie bialgebra $\l^*$, let $\CC[L^*]$ be its Hopf
algebra of regular functions. Let
$\sigma = \sum_\alpha a_\alpha \otimes b_\alpha \otimes \ell_\alpha 
\in \wedge^2(\g) \otimes \CC[L^*][1/{\bf P}]$ be a PL dynamical $r$-matrix for
$(L^*,\g,0)$. Then the quasi-Poisson structure of Section \ref{quasi:poisson}
is actually Poisson. 

Define $G/L$ as the quotient of the formal group $G$ by $L$; its
function algebra is 
$\CC[[G/L]] = \{f\in \CC[[G]] | \forall a\in \l, {\bf R}_a(f) 
= 0\}$. 

The Poisson structure of Section \ref{quasi:poisson} may be reduced to 
define a Poisson homogeneous structure 
on $G/L$, as follows. Assume that $\chi_{L^*}\in L^*_{\on{alg}}$ is a point 
where the Poisson bivector vanishes and where $\sigma$ is defined
($\chi_{L^*}$ may be viewed as a character 
of $\CC[L^*]$, such that $\chi_{L^*}(\{f,g\}) = 0$ 
for any $f,g\in\CC[L^*]$). 
Set $\sigma(\chi_{L^*}) :=  \sum_\alpha \ell_\alpha(\chi_{L^*}) 
a_\alpha \otimes b_\alpha$. 
Define a bilinear map on $\CC[[G]]$ by 
$$
\{g_1,g_2\}_{\chi_{L^*}} := \{g_1,g_2\}_G 
- m({\bf R}^{\otimes 2}(\sigma(\chi_{L^*}))(g_1\otimes g_2)). 
$$
This map restricts to a bilinear map $\CC[[G/L]]^2 \to \CC[[G/L]]$, 
because for any $a\in \l$, the vector field $\on{dress}^R_a$ vanishes 
at $\chi$. One also checks that this is a Poisson structure on $G/L$, 
which is Poisson homogeneous under the right action of the Poisson-Lie 
group $G$.

The reduction outlined here is a generalization of \cite{Lu}.

\subsection{The Poisson homogeneous structures on $G/L$} 
\label{sect:DGS}

\begin{proposition} \label{prop:basic}
Let $(\g,r_0)$ be a coboundary Lie bialgebra, i.e., $r_0\in \wedge^2(\g)$
is such that $Z(r_0) = \on{CYB}(r_0) \in \wedge^3(\g)^\g$. Let $\l\subset \g$
be a Lie subbialgebra and assume that $\g = \l\oplus \u$, where $[\l,\u]\subset 
\u$. Let $\cP$ be the set of all elements $\rho \in \wedge^2(\u)^\l$, such that the equality 
$\on{CYB}(\rho) = Z(r_0)$ holds in $\wedge^3(\g/\l)$. The Poisson homogeneous 
structures on $G/L$ correspond bijectively to the elements of $\cP$. 
\end{proposition}

{\em Proof.} The proof is straightforward. The Poisson structure 
corresponding to $\rho$ is $\{f,g\}_{\rho} = 
- m \big( (\LL^{\otimes 2}(r_0) + \RR^{\otimes 2}(\rho))(f\otimes g)\big)$. 
\hfill \qed \medskip 

\begin{remark} \label{rem:get:rho}

1) If $r(\lambda)$ is a dynamical $r$-matrix for $(\l,\g,Z(r_0))$, 
and $\chi\in\l^*$ is a character such that 
$r(\chi)$ is defined and belongs to $\wedge^2(\u)$, then 
$r(\chi)$ belongs to $\cP$. 

2) The cobracket of $\g$ is $\delta(a) = [a^1 + a^2,r_0]$.  
Assume that $r_0 = r_{0,\l} + r_{0,\u}$, where $r_{0,\l}\in \wedge^2(\l)$, 
$r_{0,\u}\in \wedge^2(\u)$. 
If $L^*_{\on{alg}}$ is an algebraic PL group with Lie bialgebra $\l^*$; if 
$\sigma : L^*_{\on{alg}} \to \wedge^2(\g)$ is a PL
dynamical $r$-matrix for $(L^*_{\on{alg}},\g,0)$; if $\chi_{L^*}
\in L^*_{\on{alg}}$ is a point where 
the Poisson structure vanishes, such that $\sigma(\chi_{L^*})$ 
is defined and belongs to $\wedge^2(\u)$, 
then $\sigma(\chi_{L^*}) - r_{0,\u}$ also belongs to $\cP$. 
\hfill \qed \medskip \end{remark}

When $\g$ is a simple Lie algebra, $\l\subset \g$ is a Levi subalgebra and 
$r_0$ is the standard $r$-matrix of $\g$, the authors of \cite{DGS}
described the set of all Poisson homogeneous structures on $G/L$
explicitly as follows. 

Let $\h\subset \l\subset \g$ be a Cartan subalgebra. Let $I$ be the
index set of $\l$ and let $\z$ be the center of $\l$. Then $\z\subset \h$. 

Denote by $\Delta(\l)\subset \Delta(\g) \subset \h^*$ the sets of roots of $\l$
and $\g$. Define $\Delta(\g,\l) \subset \z^*$ as the set of 
all elements of $\z^*$, which are restrictions to $\z$ of elements of
$\Delta(\g)$, and are nonzero. The elements of $\Delta(\g,\l)$ are called
quasi-roots.  

\begin{lemma}
For any $\alpha\in \Delta(\g)$, the restriction of $\alpha$ to 
$\z$ is nonzero iff $\alpha\notin \Delta(\l)$. 
\end{lemma}

{\em Proof.} Let us fix a system $\alpha_s,s\in [1,r]$ of simple roots.
Recall that $I \subset [1,r]$ is the set of indices of $\l$. 
Then $\h^* = \oplus_{i=1}^r \CC \alpha_i$ decomposes as 
$(\oplus_{i\in I} \CC \alpha_i) \oplus (\oplus_{i\notin I} \CC \alpha_i)$. 
The annihilator of the first space is $\z$, so the second space identifies with 
$\z^*$. Let $\alpha\in \Delta(\g)$ be such that $\alpha_{|\z} = 0$. Let us show
that $\alpha\in \Delta(\l)$. Set $\alpha = \sum_{s\in [1,r]}
n_s \alpha_s$. If $s\in S$, then $(\alpha_s)_{|\z} = 0$. Therefore 
$\alpha_{|\z} = \sum_{s\notin I} n_s (\alpha_s)_{|\z}$. Since 
$((\alpha_s)_{|\z})_{s\notin I}$ is a basis of $\z^*$, we get $n_s = 0$
for any $s\notin I$, so $\alpha \in \Delta(\l)$. 
\hfill \qed \medskip 

It follows that $\Delta(\g;\l)$ is the image of $\Delta(\g) \setminus 
\Delta(\l)$ by the restriction map, which is therefore a
surjection $\Delta(\g) \setminus \Delta(\l) \to \Delta(\g;\l)$. 

We fix a system $\alpha_s,s\in [1,r]$ of simple roots and denote by 
$\Delta_+(\g)$, $\Delta_+(\l)$ the corresponding systems of positive 
roots of $\g,\l$. We define $\Delta_+(\g;\l)$ as the image of 
$\Delta_+(\g) \setminus \Delta_+(\l)$ by the restriction map.

\begin{lemma} 
$\Delta(\g;\l)$ is the disjoint union 
$\Delta_+(\g;\l) \coprod (-\Delta_+(\g;\l))$.  
\end{lemma}

{\em Proof.} It suffices to show that the union is disjoint. 
Assume that $\bar\alpha,\bar\beta\in \Delta_+(\g;\l)$ are such that $\bar\alpha 
+ \bar\beta = 0$. Let $\alpha,\beta\in \Delta_+(\g) \setminus \Delta_+(\l)$ 
be preimages of
$\bar\alpha,\bar\beta$. Set $\alpha = \sum_{s=1}^r a_s \alpha_s$, 
$\beta = \sum_{s=1}^r b_s \alpha_s$, then $(\alpha + \beta)_{|\z} = 0$, 
which means that $a_s + b_s = 0$ for any $s\notin I$. Since the $a_s,b_s,s\in I$
are all $\geq 0$, we get $a_s = b_s = 0$ for any $s\notin I$. 
So $\alpha,\beta\in \Delta_+(\l)$, a contradiction. 
\hfill \qed \medskip   

If $\alpha\in \Delta(\g) \setminus \Delta(\l)$, we denote by $\bar\alpha
\in \Delta(\g;\l)$ its image by the restriction map.  

For $\alpha\in \Delta_+(\g)$, we fix elements $e_\alpha\in 
\g_\alpha$, $f_\alpha\in\g_{-\alpha}$ such that
$\langle e_\alpha,f_\alpha \rangle = 1$ ($\langle-,-\rangle$ is 
the Killing form). 

\begin{proposition} (see \cite{DGS})
$\cP$ consists of all the elements of the form 
$$
\rho = 
\sum_{\theta\in \Delta_+(\g;\l)} c(\theta) \sum_{\alpha\in \Delta_+(\g)
\setminus \Delta_+(\l) | \bar\alpha = \theta}
e_\alpha \wedge f_\alpha ,
$$ 
where $c : \Delta_+(\g;\l) \to \CC$ is such that for any 
$\theta,\theta'\in\Delta_+(\g;\l)$ such that $\theta + \theta' \in 
\Delta_+(\g;\l)$, $c(\theta + \theta')
\big( c(\theta) + c(\theta')\big) = c(\theta)c(\theta') + {1\over 4}$. 
\end{proposition}

$\cP$ therefore identifies with an algebraic subvariety of 
$\CC^{\Delta_+(\g;\l)}$. 

Define $\cP_0\subset \cP$ as the Zariski open subset defined by the 
condition that for any $\theta\in \Delta_+(\g;\l)$, we have 
$c(\theta) \neq \pm 1/2$. 

\begin{lemma}
Any element of $\Delta_+(\g;\l) \subset \z^*$ is a linear 
combination of the $\bar\alpha_s,s\notin I$; the 
corresponding coefficients are nonnegative integers, which do not 
vanish simultaneously. 
\end{lemma}

{\em Proof.} If $\theta\in \Delta_+(\g;\l)$
and $\alpha\in \Delta_+(\g) \setminus \Delta_+(\l)$ is such that 
$\bar\alpha = \theta$, then $\bar\alpha = \sum_{s\in [1,r]} n_s \alpha_s$, 
where $n_s\geq 0$ and the $n_s,s\notin I$ are not all zero. 
Therefore $\theta = \sum_{s\notin I} n_s \bar\alpha_s$, where the 
$n_s$ are $\geq 0$, integers, and not all zero. 
\hfill \qed \medskip

Let $\bar I := [1,r] \setminus I$. If $\theta\in \Delta_+(\g,\l)$, we denote by
$(n_s(\theta))_{s\in \bar I}$ the integers such that $\theta = \sum_{s\in \bar
I} n_s(\theta) \bar \alpha_s$. We denote by $\chi_\theta : (\CC^\times)^{\bar
I} \to \CC^\times$ the group morphism taking $(t_s)_{s\in \bar I}$
to $\prod_{s\in \bar I} t_s^{n_s(\theta)}$. 
If $t\in \CC^\times \setminus \{1\}$, we set $\kappa(t) := {1\over 2} {{1 + t}
\over {1-t}}$. 

\begin{proposition} \label{isom:vars}
There is a unique isomorphism of algebraic varieties 
$$
\iota : (\CC^\times)^{\bar I} \setminus \bigcup_{\theta\in \Delta_+(\g,\l)}
\on{Ker}(\chi_\theta) \to \cP_0, 
$$ 
taking $(t_s)_{s\in\bar I}$ to $(c(\theta))_{\theta\in \Delta_+(\g,\l)}$, 
where $c(\theta) = \kappa \circ \chi_\theta((t_s)_{s\in \bar I})$. 
\end{proposition}

{\em Proof.} Let us first show that $\iota$ is well defined. 
If $\tt = (t_s)_{s\in \bar I}$ belongs to 
$(\CC^\times)^{\bar I} \setminus \cup_{\theta\in \Delta_+(\g,\l)}
\on{Ker}(\chi_\theta)$, then $\chi_\theta(\tt)\neq 1$, so 
$\kappa \circ \chi_\theta(\tt)\in \CC$ is well defined. 
Moreover, the family $(\kappa \circ \chi_\theta(\tt))_{\theta\in 
\Delta_+(\g;\l)}$ satisfies the defining relations of $\cP$, because 
of the identity 
$\kappa(x+x') \big( \kappa(x) + \kappa(x')\big) = {1\over 4} 
+ \kappa(xx')$ for any $x,x'\in \CC^\times$ such that $x,x',xx'\neq 1$. 
Finally, $\kappa\circ \chi_\theta(\tt)$ belongs to the image of 
$\kappa$, hence is not equal to $1/2$ or $-1/2$. It follows that 
$(\kappa \circ \chi_\theta(\tt))_{\theta\in \Delta_+(\g;\l)}$
belongs to $\cP_0$. 

$\iota$ is clearly a morphism of algebraic varieties. To show 
that is it an isomorphism, we will construct the inverse morphism. 
We claim that there is a unique morphism  
$$
\iota' : \cP_0 \to (\CC^\times)^{\bar I} \setminus 
\bigcup_{\theta\in \Delta_+(\g;\l)}
\on{Ker}(\chi_\theta), 
$$ 
taking $(c(\theta))_{\theta\in \Delta_+(\g;\l)}$ to
$\tt = (t_s)_{s\in \bar I}$, defined by $t_s = \kappa^{-1}(c(\bar \alpha_s))$
(we have $\kappa^{-1}(y) = {{2y-1}\over{2y+1}}$). Let us show that 
$\iota'$ is well defined. Since each $c(\bar\alpha_s)$ is $\neq
\pm 1/2$, each $t_s$ is defined and belongs to $\CC^\times$. 
Moreover, using the fact that  
$(c(\theta))_{\theta\in \Delta_+(\g;\l)}$ obeys the defining relations of 
$\cP$, one proves by induction on $\on{ht}(\theta) = \sum_{s\in \bar I}
n_s(\theta)$ that $1 + \chi_\theta(\tt) = 2(1 - \chi_\theta(\tt))c(\theta)$. 
So $\chi_\theta(\tt) \neq 1$ for any $\theta\in \Delta_+(\g;\l)$. 
So $\iota'$ is well defined. 

It is straightforward to check that $\iota$ and $\iota'$ are inverse to each
other. 
\hfill \qed \medskip 

If $\tt \in (\CC^\times)^{\bar I} \setminus \cup_{\theta\in 
\Delta_+(\g;\l)} \on{Ker}(\chi_\theta)$, we denote by $\rho_{\tt}$
the element of $\cP_0$ corresponding to $\tt$ under $\iota$. Then 
$$
\rho_{\tt} = {1\over 2} \sum_{\theta\in\Delta_+(\g;\l)}
{{1+\chi_\theta(\tt)}\over{1-\chi_\theta(\tt)}}
\sum_{\alpha\in \Delta_+(\g) \setminus \Delta_+(\l) | \bar\alpha = \theta}
e_\alpha \wedge f_\alpha. 
$$

\subsection{The Poisson structures from $\cP_0$ are reductions}
\label{sect:4:2}

In this section, we show that the elements of $\cP_0$ can be obtained 

(a) using a dynamical $r$-matrix $r_\l^\g$ and a character $\chi\in \l^*$, 
as in Remark \ref{rem:get:rho}, 1) (this result was also obtained in 
\cite{KMST}, Theorem 15); 

(b) using the PL dynamical $r$-matrix $\sigma_\l^\g$ and an element
$\chi_{L^*} \in L^*_{\on{alg}}$, as in Remark \ref{rem:get:rho}, 2). 
 
This shows that any Poisson structure from $\cP_0$ is 

(a) a reduction of a Poisson structure on $\l^* \times G$, as in \cite{Lu}
 
(b) a reduction of a Poisson structure on $L^*_{\on{alg}}\times G$, 
of the type described in Section \ref{sect:basic}.

\subsubsection{$\cP_0$ and $r_\l^\g$}

It follows from \cite{FGP,EE2} that there is a unique dynamical $r$-matrix
$r_\l^\g : \l^* \to \wedge^2(\g)$ for $(\l,\g,Z(r_0))$ (see Example 
\ref{from:old}), such that 
$$
\forall \lambda\in \h^*, \; 
r_\l^\g(\lambda) = - {1\over 2}
\sum_{\alpha\in \Delta_+(\g) \setminus \Delta_+(\l)}
\on{coth}{{(\lambda,\alpha)}\over 2} e_\alpha \wedge f_\alpha 
+ \sum_{\alpha\in \Delta_+(\l)}
\varphi(\lambda,\alpha) e_\alpha \wedge f_\alpha,  
$$
where $\varphi(x) = - {1\over 2} \on{coth}(x/2) + 1/x$.

Recall that $\{$characters of $\l\} = \z^*$. If $\lambda$ is a character of
$\l$, then $(\beta,\lambda) = 0$ for any $\beta\in \h_\l^*$, therefore
$(\alpha,\lambda) = 0$ if $\alpha\in \Delta_+(\l)$, and if $\alpha\in
\Delta_+(\g) \setminus \Delta_+(\l)$, then $(\alpha,\lambda)$ depends only 
on $\bar\alpha\in \Delta_+(\g;\l)$. We denote it by $(\bar\alpha,\lambda)$. 

\begin{definition}
{\it A character $\lambda$ of $\l$ is called regular iff for any $\theta\in
\Delta_+(\g;\l)$, $(\theta,\lambda) \notin 2\pi i \ZZ$. } 
\end{definition}

There is a unique map $\{$regular characters of $\l\} \to (\CC^\times)^{\bar I}
\setminus \cup_{\theta\in \Delta_+(\g;\l)} \on{Ker}(\chi_\theta)$, taking
$\lambda$ to $\tt(\lambda) = (t_s)_{s\in \bar I}$, where $t_s = 
e^{(\bar\alpha_s,\lambda)}$. The map $\lambda \mapsto \tt(\lambda)$
is a covering with group $\ZZ^{\bar I}$. 

\begin{lemma}
If $\chi$ is a regular character of $\l$, we have $\rho_{\tt(\chi)}
= r_\l^\g(\chi)$. 
\end{lemma}

{\em Proof.} We have $(\alpha,\chi) = 0$ if $\alpha\in \Delta_+(\l)$, 
which simplifies the expression of $r_\l^\g(\chi)$. \hfill \qed \medskip 

It follows from \cite{Lu} that the bilinear map on $\CC[[G]])$ defined by 
$$
\{g_1,g_2\}_{\on{Lu}} = m({\bf L}^{\otimes 2}(r_0)(g_1 \otimes g_2))
+ m({\bf R}^{\otimes 2}(r_\l^\g(\chi))(g_1 \otimes g_2))
$$
restricts to a Poisson structure on $\CC[[G/L]]$. 

\begin{corollary}
The Poisson structures $\{-,-\}_{\on{Lu}}$ and $\{-,-\}_{\rho_{\tt}}$
coincide. 
\end{corollary}

\subsubsection{$\cP_0$ and $\sigma_\l^\g$}

Let us restrict $\sigma_\l^\g : L^*_{\on{alg}} \to \wedge^2(\g)$
to $H^*_{\on{alg}} \subset L^*_{\on{alg}}$. $H^*_{\on{alg}}$ is a torus 
with Lie algebra $\h^*$. We have 
$$
\sigma_\l^\g(e^\chi)  = \sum_{\alpha\in \Delta_+(\l) \setminus \Delta_+(\g)}
{{e_\alpha \wedge f_\alpha}\over{ 1- e^{-(\chi,\alpha)}}}
$$
for any $\chi\in \h^*$ such that the denominators are nonzero.  
 
\begin{lemma} \label{lemma:4:11}
There is a unique subtorus $Z^*_{\on{alg}} \subset H^*_{\on{alg}}$ 
with Lie algebra $\z^*$. If $\tt\in (\CC^\times)^{\bar I} \setminus
\cup_{\theta\in \Delta_+(\g;\l)} \on{Ker}(\chi_\theta)$, there exists a
regular $\chi\in \z^*$ such that $e^{-(\chi,\bar\alpha_s)} = t_s$
for any $s\in \bar I$. Set $\chi_{L^*} = e^\chi$. Then the Poisson 
structure of $L^*_{\on{alg}}$ vanishes at $\chi_{L^*}$, and $\rho_{\tt} = 
\sigma_\l^\g(\chi_{L^*}) - r_{0,\u}$. 
\end{lemma}

\begin{corollary}
The Poisson brackets $\{-,-\}_{\chi_{L^*}}$ and $\{-,-\}_{\rho_{\tt}}$
(see Section \ref{sect:basic}) coincide. 
\end{corollary} 

{\em Proof.} If $g_1,g_2\in \CC[[G]]$, then $\{g_1,g_2\}_G 
= m((\RR^{\otimes 2} - \LL^{\otimes 2})(r_0)(g_1 \otimes g_2))$, 
so 
\begin{align*}
& \{g_1,g_2\}_{\chi_{L^*}} = 
- m(\LL^{\otimes 2}(r_0)(g_1 \otimes g_2)) 
- m(\RR^{\otimes 2}(\sigma_\l^\g(e^\chi) - r_0)(g_1 \otimes g_2))
\\ & 
= \{g_1,g_2\}_{\rho_{\tt}} 
+ m(\RR^{\otimes 2}(r_{0,\l})(g_1 \otimes g_2))
= \{g_1,g_2\}_{\rho_{\tt}} 
\end{align*} 
since $g_1$ and $g_2$ are $\l$-invariant. \hfill \qed \medskip

\section{Quantization of Poisson homogeneous structures on $G/L$}
\label{sect:5}

In this section, we construct quantizations of all the Poisson homogeneous 
structures on $G/L$, corresponding to elements of $\cP_0$. For this, we prove
algebraicity results for $\Psi_\h^\g$ and $\Psi_\l^\g$, relying on the
computation of \cite{DCK} of the determinant of the Shapovalov pairing for
$U_q(\g)$. We can then construct the element $(\Psi_\l^\g)^{-1}(\chi)$, which
serves to quantize $G/L$ in the usual way. 

\subsection{Algebraicity of $\Psi_\h^\g$}

Let $U_\hbar(\h) \subset U_\hbar(\g)$ be the subalgebra 
generated by the $h_i$. Then $U_\hbar(\h) = \CC[h_i, i = 1,\ldots,r][[\hbar]]$. 
The intersection of this algebra with $\CC[[G^*]]_\hbar$ is 
$\CC[[H^*]]_\hbar = \CC[[\hbar h_i,\hbar]]$. 

Set $q = e^\hbar$, $k_i = e^{\hbar h_i}$. If $\beta = 
\sum_i n_i \alpha_i$ is a positive root of $\g$, we set 
$k_\beta = \prod_{i=1}^r k_i^{n_i}$. Set $\CC[H^*]_\g = 
\CC[k_i^{\pm 1}, i = 1,\ldots,r]_{\g}[[\hbar]]$, where the index 
means that we invert all the $k_\beta^2 - 1$, where $\beta\in\Delta_+(\g)$.  

\begin{proposition}
1) $\Upsilon_\h^\g$ belongs to $U_\hbar(\n_+) \wh\otimes U_\hbar(\n_-)
\wh\otimes \CC[H^*]_{\g}$.  

2) $\Psi_\h^\g$ belongs to $U_\hbar(\n_+) \wh\otimes U_\hbar(\b_-)
\wh\otimes \CC[H^*]_{\g}$.  
\end{proposition}

{\em Proof.} Let us prove 1). $\Upsilon_\h^\g$ is the sum 
$\sum_{\alpha\in\ZZ_{+}^r} (\Upsilon_\h^\g)_\alpha$, where 
$(\Upsilon_\h^\g)_\alpha$ is inverse to the quantum Shapovalov pairing 
$U_\hbar(\n_+)_\alpha \otimes U_\hbar(\n_-)_{-\alpha} \to U_\hbar(\h)$, 
$x\otimes y\mapsto h(yx)$. Here $U_\hbar(\n_\pm) \subset U_\hbar(\g)$
are the subalgebras generated by the $e_i$ (resp., the $f_i$), and 
the index $\pm\alpha$ means the degree $\pm\alpha$ part. 
 
Let $U_q(\g)\subset U_\hbar(\g)$ be the 
$\CC(q)$-subalgebra generated by the $e_i,f_i$ and $k_i$, and 
$U_q(\n_\pm)$, $U_q(\h)$ the $\CC(q)$-subalgebras of $U_q(\g)$
generated by the $e_i$ (resp., the $f_i$, the $k_i^{\pm 1}$).
Then the Shapovalov pairing restricts to a pairing 
$U_q(\n_+)_\alpha \otimes U_q(\n_-)_{-\alpha} \to U_q(\h)$, 
the determinant of 
which was computed in \cite{DCK}, Proposition 1.9, (a). 
Let $\on{det}_\alpha$ be this determinant; then 
$$
(\Upsilon_\h^\g)_\alpha
\in U_q(\n_+) \otimes U_q(\n_-) \otimes U_q(\h)[1/\on{det}_\alpha], 
$$
where the tensor products are over $\CC(q)$. 
It follows from \cite{DCK} that if we identify $\on{det}_\alpha$
with its image in $\CC[k_i^{\pm 1}, i = 1,\ldots,r][[\hbar]]$, then 
we have 
$$
\on{det}_\alpha = u_\alpha \hbar^{- \sum_{n\geq 0} P(\alpha - n\beta)}
\prod_{\beta\in\Delta_+(\g)} (k_\beta^2 - 1)^{\sum_{n\geq 0} P(\alpha - n\beta)}
(1+ O(\hbar)), 
$$ 
where $P(\alpha)$ is the number of decompositions $\alpha =
\sum_{\beta\in\Delta_+(\g)} k_\beta \beta$, with $k_\beta\in \ZZ_+$
and $u_\alpha$ is invertible. Therefore 
$\on{det}_\alpha$ is invertible in $\CC[H^*]_{\g}$.
So $(\Upsilon_\h^\g)_\alpha \in U_\hbar(\n_+)_\alpha \wh\otimes 
U_\hbar(\n_-)_{-\alpha} \wh\otimes \CC[H^*]_{\on{loc}}$. We also know that
the $\hbar$-adic valuation of $(\Upsilon_\h^\g)_\alpha$ tends to infinity with
$\alpha$, which implies 1). 2) then follows from the fact that the "left 
coproduct" $\Delta_l$ takes 
$\CC[H^*]_{\g}$ to $\CC[H^*]_{\g} \bar\otimes 
\CC[[\hbar h_i,\hbar]]$.    
\hfill \qed\medskip 

\subsection{$(\Psi_\l^\g)^{-1}(\chi)$ and its properties}

Set $\l_+ = \oplus_{\alpha\in\Delta_+(\l)} \CC e_\alpha$, 
$\l_- = \oplus_{\alpha\in\Delta_+(\l)} \CC f_\alpha$. 
Denote by $\wh H_\l^\g : \CC[[L^*]]_{P_\l^\g,\hbar} \to 
\CC[[H^*]]_{(P_\l^\g)_{|H^*},\hbar}$
the Harish-Chandra map as in Section \ref{sect:3:3}. 

Set $\CC[H^*]_{\g;\l} := \CC[k_i^{\pm 1}, i = 1,\ldots,r]_{\g;\l}[[\hbar]]$,
where the index means that we invert all the $k_\beta^2 - 1$, $\beta \in
\Delta_+(\g) \setminus \Delta_+(\l)$. 

\begin{lemma}
Set $\Xi := (\id\otimes \id \otimes \wh H_\h^\l)((\Psi_\l^\g)^{-1})$. 
Then $\Xi \in U_\hbar(\u_+) \wh\otimes U_\hbar(\p_-) \wh\otimes 
\CC[H^*]_{\g;\l}$. 
\end{lemma}

{\em Proof.} We already know that $(\Psi_\l^\g)^{-1} \in U_\hbar(\u_+)
\wh\otimes U_\hbar(\p_-) \bar\otimes \CC[[L^*]]_{P_\l^\g,\hbar}$. 
Therefore 
$$
(\id\otimes \id\otimes \wh H_\h^\l)((\Psi_\l^\g)^{-1} )
\in U_\hbar(\u_+) \wh\otimes U_\hbar(\p_-) \bar\otimes 
\CC[[H^*]]_{(P_\l^\g)_{|H^*},\hbar},  
$$
where we recall that $(P_\l^\g)_{|H^*} = \prod_{\beta \in \Delta_+(\g) \setminus
\Delta_+(\l)} (k_\beta^2 - 1)$. 

On the other hand, Proposition \ref{quantum:comp} implies that 
$\eta((\Psi_\l^\g)^{-1})$ belongs to $U_\hbar(\g)^{\wh\otimes 2}
\wh\otimes \CC[H^*]_\g$. Therefore $\eta'\circ\eta((\Psi_\l^\g)^{-1})
= (\id\otimes \id\otimes \wh H_\l^\g)((\Psi_\l^\g)^{-1})$ 
belongs to the same space.

Now $\CC[[H^*]]_{(P_\l^\g)_{|H^*},\hbar} \cap \CC[H^*]_\g = 
\CC[H^*]_{\g;\l}$, which implies the lemma. \hfill \qed \medskip 

Let us denote by $\bar H : (\CC[[L_-^*]]_\hbar)_0 \wh\otimes 
(\CC[[L_+^*]]_\hbar)_0 
\to \CC[[H^*]]_\hbar$ the map $x\otimes y \mapsto H_\h^\l(xy)$ 
(here $A_0$ denotes the augmentation ideal of an augmented algebra $A$). 

\begin{lemma} \label{lemma:id:Xi}
The identity 
\begin{equation} \label{id:Xi}
\Xi^{2,3,4} \Xi^{1,23,4} = (\id^{\otimes 3} \otimes H_\h^\l)
\Big( \big( (\Psi_\l^\g)^{-1} \big)^{1,2,34 }\Big) \Xi^{12,3,4} 
\end{equation}
holds in $U_\hbar(\g)^{\wh\otimes 3} \wh\otimes 
\big( \CC[H^*]_{\g;\l}/J\big)$. Here $I \subset 
\CC[[H^*]]_{(P_\l^\g)_{|H^*},\hbar}$ is the ideal generated 
by $\on{Im}(\bar H)$ and $J := I \cap \CC[H^*]_{\g;\l}$. 
\end{lemma}

{\em Proof.} Apply $\id^{\otimes 3} \otimes \wh H_\h^\l$ 
to the equation derived from
the dynamical twist equation (\ref{dyn:twist}) by taking the 
inverses of both sides. We have 
$(\Psi_\l^\g)^{-1} = \Xi + \Theta$, where $\Theta = \sum_i 
a_i \otimes b_i \otimes c_i \ell_i d_i$, and $a_i,b_i\in U_\hbar(\g)$, 
$\ell_i \in \CC[[H^*]]_{(P_\l^\g)_{|H^*},\hbar}$, 
$c_i\in \CC[[L_+^*]]_\hbar$, $d_i
\in \CC[[L_-^*]]_\hbar$, and $\eps(c_i) \eps(d_i) = 0$. Then if 
$\xi \in \CC[H^*]_{\g;\l}$, then $H_\h^\l(\xi c_i \ell_i d_i) 
= H_\h^\l(c_i \ell_i d_i \xi)
= \eps(c_i) \eps(d_i) \xi\ell_i = 0$, and 
$H_\h^\l(c_i \ell_i d_i c_j \ell_j d_j) = \eps(c_i) \eps(d_j)
\ell_i H_\h^\l(d_i c_j) \ell_j$. So the "$\Theta$-$\Xi$" contributions 
vanish, and the "$\Theta$-$\Theta$" contribution involves only pairs 
$(i,j)$ such that $\eps(c_i)\eps(d_j)\neq 0$; if $(i,j)$ is such a pair, 
then $\eps(d_i) = \eps(c_j) = 0$; so this contribution belongs to 
$I = (\on{Im}(\bar H))$. So (\ref{id:Xi}) holds in 
$U_\hbar(\g)^{\wh\otimes 3} \bar\otimes \big( \CC[[H^*]]_{
(P_\l^\g)_{|H^*},\hbar}/I\big)$. Now both sides of (\ref{id:Xi})
belong to $U_\hbar(\g)^{\wh\otimes 3} \wh\otimes \CC[H^*]_{\g;\l}$, which
implies the lemma. \hfill \qed \medskip 

Let $\chi$ be a character of $\l$. As $\chi$ is uniquely determined by
its restriction of $\h$, we identify it with an element of $\h^*$. 
Recall that $\chi$
is regular iff for any $\beta \in \Delta_+(\g) \setminus \Delta_+(\l)$, we have
$\chi(h_\beta) \notin 2\pi i \ZZ$. Here for $\beta = \sum_{s=1}^r n_s \alpha_s
\in \Delta_+(\g)$, we set $h_\beta = \sum_{s=1}^r n_s h_s$. 
It follows that $\chi$ is regular iff $\sigma_\l^\g$ 
is defined at $e^\chi$. 

Then if $\chi$
is regular, there is a unique character $\bar\chi : \CC[H^*]_{\g;\l} \to
\CC[[\hbar]]$, taking each $k_s$ to $\on{exp}(\chi(h_s)/2)$. We generalize 
these definitions to the case of a character $\chi : \l\to\CC[[\hbar]]$ 
as follows: $\chi$ is regular iff $(\chi(h_\beta)\on{\ mod\ }\hbar)
\notin 2\pi i \ZZ$ for any $\beta \in \Delta_+(\g) \setminus \Delta_+(\l)$;
then $\bar\chi$ is defined in the same way.

Recall from Lemma \ref{lemma:id:Xi} that $J$ is an ideal of 
$\CC[H^*]_{\g;\l}$.

\begin{lemma} \label{lemma:J:chi}
If $\chi$ is a regular character of $\l$, then $J \subset
\on{Ker}(\bar\chi)$.  
\end{lemma} 

{\em Proof.} Let $\z$ be the center of $\l$, and set $\h_\l := \l' \cap \h$,
where $\l' = [\l,\l]$. Then $\h = \h_\l \oplus \z$. $\h_\l$
is spanned by the $h_s$, where $s\in [1,r]$ is an index of $\l$.  

$\CC[[H^*]]_{(P_\l^\g)_{|H^*},\hbar}$ is the $\hbar$-adic completion of 
$\CC[[h_s,\hbar, s = 1,\ldots,r]][(k_\beta^2-1)^{-1},\beta\in \Delta_+(\g)
\setminus \Delta_+(\l)]$, and $I$ is the complete ideal of 
$\CC[[H^*]]_{(P_\l^\g)_{|H^*},\hbar}$ generated by $\h_\l$, i.e., by the $h_s$, 
$s$ an index of $\l$. 

Recall that $\CC[H^*]_{\g;\l} = \CC[k_s^{\pm 1}, s = 1,\ldots,r]
[(k_\beta^2-1)^{-1},\beta\in\Delta_+(\g) \setminus \Delta_+(\l)][[\hbar]]$; then
$J$ is a complete ideal of $\CC[H^*]_{\g;\l}$, generated by the $k_s - 1$, $s$
an index of $\l$. Then if $\chi$ is a regular character of $\l$, and $s$ is an
index of $\l$, then $h_s\in [\l,\l]$, therefore  
$\chi(h_s) = 0$ and so $\bar\chi(k_s) = 
\on{exp}(\chi(h_s)/2) = 1$, so $\bar\chi(J) = 0$.  
\hfill \qed \medskip

\begin{proposition}  \label{prop:pre:assoc}
If $\chi$ is a regular character of $\l$, set $\Xi(\chi) := 
(\id\otimes \id\otimes \bar\chi)(\Xi)$. Then 
$$
\Xi(\chi)^{2,3} \Xi(\chi)^{1,23} = 
\big( \id^{\otimes 3} \otimes (\bar\chi\circ H_\h^\l)\big) 
\Big( \big( (\Psi_\l^\g)^{-1}\big)^{1,2,34}\Big) \Xi(\chi)^{12,3}. 
$$
\end{proposition}

{\em Proof.} This follows immediately from Lemma \ref{lemma:id:Xi}
and Lemma \ref{lemma:J:chi}. 
\hfill \qed \medskip 

\subsection{Quantized Poisson homogeneous structures} \label{sect:qP}

In this section, we construct quantizations of all the Poisson structures of
Section \ref{sect:4:2}. 

As before, $\g$ is a semisimple Lie algebra, $\l\subset \g$ is
a Levi subalgebra,
and $\chi\in\l^*$ is a regular character of $\l$. A quantization of the
formal Poisson manifold $L\setminus G$, equivariant under $U_\hbar(\g)$, is
constructed as follows. 

Set $\CC[[G]]_\hbar := \on{Hom}(U_\hbar(\g),\CC[[\hbar]])$ and define the left and
right actions of $U_\hbar(\g)$ on $\CC[[G]]_\hbar$ by $({\mathbf L}_a f)(x) = 
f(S(a) x)$, $({\mathbf R}_a f)(x) = f(xa)$ for any $a,x\in U_\hbar(\g)$. Then we
set 
$$
\CC[[G/L]]_\hbar = \{ f\in \CC[[G]]_\hbar | \forall a\in U_\hbar(\l), 
{\mathbf R}_a f = \eps(a) f\}. 
$$
$\CC[[G/L]]_\hbar$ is a subalgebra of $\CC[[G]]_\hbar$; we denote its
product by $m_\hbar$. 

The right action ${\mathbf R}$ restricts to an action of $U_\hbar(\g)$ on 
$\CC[[G/L]]_\hbar$. If $f_1,f_2\in \CC[[G/L]]_\hbar$, we set 
$$
f_1 * f_2 = m_\hbar \big( (\RR\circ S)^{\otimes 2}(\Xi(\chi))
(f_1 \otimes f_2)\big).   
$$

\begin{theorem} $\big( \CC[[G/L]]_\hbar,*\big)$ is a quantization of
$G/L$ with its Poisson structure described in Section \ref{sect:4:2}, 
equivariant under $U_\hbar(\g)$. 
\end{theorem}

{\em Proof.} Let us prove that $*$ is associative. If 
$x\in \CC[[L^*]]_{P_\l^\g,\hbar}$ is such
that $(\id\otimes \wh H_\l^\g) \circ \Delta(x) \in U_\hbar(\l) \wh\otimes
\CC[H^*]_{\g;\l}$, then applying $\eps$ to the first factor we get 
$\wh H_\l^\g(x)
\in \CC[H^*]_{\g;\l}$. Then $\big( \on{id} \otimes 
(\bar\chi \circ \wh H_\l^\g)\big) \circ \Delta(x) \in U_\hbar(\l)$. 
Therefore, if $h\in \CC[[G/L]]_\hbar$,
we get 
$$
{\mathbf R} \Big( \big( \id  \otimes (\bar\chi \circ \wh H_\l^\g)\big) 
\circ \Delta(x) \Big) (h) = \big(\eps\otimes (\bar\chi \circ \wh H_\l^\g)\big) 
\circ \Delta(x) h = (\bar\chi \circ \wh H_\l^\g)(x) h. 
$$
It follows that if $f,g,h \in \CC[[G/L]]_\hbar$, then 
\begin{align*}
& ( (\RR \circ S)^{\otimes 2} \otimes \RR\otimes \id)
\Big( \big(\id^{\otimes 3} \otimes (\bar\chi\circ \wh H_\l^\g)\big)
( (\Psi_\l^\g)^{-1} )^{1,2,34} \Big) 
(f\otimes g \otimes h)
\\ & = ( (\RR \circ S)^{\otimes 2} \otimes \id)
\Big( \big(\id\otimes \id \otimes (\bar\chi\circ \wh H_\l^\g)\big)
((\Psi_\l^\g)^{-1})\Big) 
(f\otimes g \otimes h)
\\ & = ( (\RR \circ S)^{\otimes 2} \otimes \id)
(\Xi(\chi)^{1,2})(f\otimes g \otimes h).
\end{align*}

Using this identity together with Proposition \ref{prop:pre:assoc}, we obtain
the associativity of $*$. The classical limit of $*$ is then the Poisson
structure of Section \ref{sect:4:2}. One shows easily that ${\mathbf R}$
restricts to an action of $U_\hbar(\g)$ on $\CC[[G/L]]_\hbar$, compatible
with $*$. 
\hfill \qed \medskip

\section{Relation with \cite{EE2}} \label{sect:6}

In this section, we compare the quantizations of the Poisson homogeneous 
structures on $G/L$ obtained in Section \ref{quant:homog} and in \cite{EE2}. 
We first recall the construction of
\cite{EE2} (Section \ref{recall:EE2}). We then prove a rigidity result
satisfied by the $\sigma_\l^\g(\chi)$ (Section \ref{sect:rigidity}). 
Finally, we compare both constructions (Section \ref{sect:compare}). 
As before, $\g$ is a semisimple Lie algebra and $\l\subset\g$ is a 
Levi subalgebra. 

\subsection{The construction of \cite{EE2}} \label{recall:EE2}

Let $\Phi$ be a Drinfeld associator, $t_0\in S^2(\g)^\g$ be the Casimir 
element, and set $\Phi_\g := \Phi(\hbar t_0^{1,2},\hbar t_0^{2,3})$. 
According to \cite{Dr:QH}, there exists $J_0 \in U(\g)^{\otimes 2}[[\hbar]]$, 
such that $J_0 = 1 + \hbar r_0 + O(\hbar^2)$, and 
\begin{equation} \label{J0:Phi}
J_0^{12,3} J_0^{1,2} = \Phi_\g^{-1} J_0^{1,23} J_0^{2,3}.  
\end{equation}
Then $(U(\g)[[\hbar]],\on{Ad}(J_0^{-1}) \circ \Delta_0)$ is 
isomorphic to $(U_\hbar(\g),\Delta_\hbar)$. 

In \cite{EE2}, we construct a solution $J \in U(\g)^{\otimes 2}
\wh\otimes \wh{U(\l)}[[\hbar]]$ of the dynamical pseudotwist equation 
$$
J^{12,3,4}J^{1,2,34} = (\Phi_\g^{-1})^{1,2,3} 
J^{1,23,4}J^{2,3,4}. 
$$
Here $\wh{U(\l)}$ is the microlocalization of $U(\l)$ w.r.t. $D_\l^\g$, which
may be defined as the unique element of $S^{\on{dim}(\u_+)}(\l)$ whose
restriction to $\h$ coincides with $\prod_{\beta \in \Delta_+(\g) \setminus 
\Delta_+(\l)}\beta$. 

Recall that $\{$characters of $\l\} = \z^*$. Then if $\chi$ belongs to 
an analytic open subset $U_{\z^*} \subset \z^*$, then 
$J(\chi) = (\id^{\otimes 2}
\otimes \hbar^{-1}\chi)(J)$ is well defined. For $f_1,f_2 \in 
\CC[[G/L]][[\hbar]]$, we set 
$$
f_1 * f_2 = m_0 \big( \LL^{\otimes 2}(J_0)
\RR^{\prime \otimes 2}(J(\chi)^{-1}) (f_1 \otimes f_2)\big).  
$$
(we define $\RR'$ by $\RR \circ S$). 
This defines a star-product on $G/L$, equivariant under the left
action of $(U(\g)[[\hbar]],\on{Ad}(J_0^{-1}) \circ \Delta_0) \simeq 
U_\hbar(\g)$, quantizing the Poisson structure of Section \ref{sect:basic}. 

Let us study the effect of a gauge transformation of $J_0$ on this
construction. Let $u\in U(\g)[[\hbar]]$ be an invertible element, 
and set $u * J_0 = \Delta_0(u) J_0 (u\otimes u)^{-1}$. Then 
$u * J_0$ is a solution of (\ref{J0:Phi}); we denote by $*_u$ the corresponding
star-product. Then $f\mapsto \RR'_u f$ is an isomorphism between the
star-products $*$ and $*_u$. 

We now prove: 

\begin{proposition} \label{J0:eq}
The solutions $J_0$ of (\ref{J0:Phi}) are all gauge-equivalent. 
\end{proposition}

{\em Proof.} According to \cite{EK}, there exists a universal series
$\rho\mapsto J^\Phi_{\on{EK}}(\rho)$, such that if $\rho$ is a solution of the 
CYBE, then $J_\rho := J_{\on{EK}}^\Phi(\rho)$ satisfies 
$J_\rho^{12,3}J_\rho^{1,2} = \Phi_\rho^{-1} J_\rho^{1,23}J_\rho^{2,3}$. 
Here $\Phi_\rho = \Phi(\hbar t_\rho^{1,2},\hbar t_\rho^{2,3})$, where
$t_\rho = \rho + \rho^{2,1}$. 
Then according to 
\cite{EE3}, we get $J_0 = v * J^\Phi_{\on{EK}}(\rho_0)$, where $\rho_0
\in \g^{\otimes 2}[[\hbar]]$ is such that $(\on{mod}\ \hbar)(\rho_0) = r_0$, 
$\rho_0 + \rho_0^{2,1} = t_0$, $\rho_0$ satisfies the CYBE, and $v\in
U(\g)[[\hbar]]$ is invertible. According to the Belavin-Drinfeld classification
of solutions of the CYBE, $\rho_0$ is conjugate to $r_0 + h_0$, where 
$h_0 \in \hbar \wedge^2(\h)[[\hbar]]$. Therefore 
$J_0 = v_1 * J_{\on{EK}}^\Phi(r_0 + h_0)$, 
where $v_1 \in U(\g)[[\hbar]]^\times$. Now the QUE algebra
$\big( U(\g)[[\hbar]], 
\on{Ad}(J_{\on{EK}}^\Phi(r_0 + h_0)^{-1})\circ \Delta_0\big)$
is isomorphic to $(U_\hbar(\g),\Delta_\hbar)$, and therefore contains $r 
= \on{rank}(\g)$ QUE subalgebras isomorphic to $U_\hbar({{\mathfrak sl}}_2)$, 
corresponding to each simple root. 

Applying the dequantization functor corresponding to the associator $\Phi$, 
we obtain the Lie bialgebra $(\g[[\hbar]],\mu,\on{ad}(r_0 + h_0))$; this Lie
bialgebra should contain $r$ Lie subbialgebras isomorphic to $({{\mathfrak
sl}}_2[[\hbar]],$ standard structure), whose reductions modulo $\hbar$
are the subalgebras corresponding to the simple roots of $\g$. One checks that
such Lie subbialgebras can exist only if $h_0 = 0$. Therefore $J_0$
is gauge-equivalent to $J^\Phi_{\on{EK}}(r_0)$. \hfill \qed \medskip  

It follows that the star-products $*$ and $*_u$ are gauge-equivalent. 
We will now assume that $J_0 = J^\Phi_{\on{EK}}(r_0)$. 

\begin{remark} Assume that $\g$ is simple. 
According to \cite{Dr:QH}, we have a bijection 
$\{\g$-invariant solutions of the pentagon equation in 
$U(\g)^{\otimes 3}[[\hbar]]\}/$(invariant 
twists)$\simeq \wedge^3(\g)^\g[[\hbar]] \simeq
\CC[[\hbar]]$, since $\g$ is simple. It follows that all the $\Phi_\g$
coincide, when $\Phi$ runs over all associators, up to a rescaling of $\hbar$. 
Moreover, the $\Phi_\g$ all satisfy the same hexagon relation, so they all
coincide up to twist. Then Proposition \ref{J0:eq} implies that if $\Phi,\Phi'$ 
are two associators, there exists an invertible 
$T\in (U(\g)^{\otimes 2})^\g[[\hbar]]$, such 
that $J_{\on{EK}}^\Phi(r_0)$ and $TJ_{\on{EK}}^{\Phi'}(r_0)$
are gauge-equivalent. 
\end{remark}

\subsection{Comparison of Levi subalgebras} \label{comp:Levi}

\begin{theorem} \label{thm:levi}
There exists an invertible element $u\in U(\g)^\h[[\hbar]]$, such that 
the subalgebra $\on{Ad}(u)(U(\l)[[\hbar]]) \subset U(\g)[[\hbar]]$ identifies
with $U_\hbar(\l) \subset U_\hbar(\g)$ under $\big( U(\g)[[\hbar]], 
\on{Ad}(J_{\on{EK}}^\Phi(r_0)^{-1}) \circ \Delta_0 \big) \simeq
(U_\hbar(\g),\Delta_\hbar)$. 
\end{theorem}

{\em Proof.} Let $Q$ be the quantization functor associated to $\Phi$ by
\cite{EK}. Apply $Q$ to the sequence of inclusions of Lie bialgebras $\h\subset
\l \subset \g$. We get the sequence of inclusions $U_\hbar(\h) \subset
U_\hbar(\l) \subset U_\hbar(\g)$. We have an isomorphism $U_\hbar(\l) \simeq
U(\l)[[\hbar]]$, such that the composed map $\h\subset U_\hbar(\l) \simeq
U(\l)[[\hbar]]$ is the standard inclusion (see \cite{EK6}). We have also 
an isomorphism $U_\hbar(\g) \simeq
U(\g)[[\hbar]]$, such that $\h\subset U_\hbar(\g) \simeq
U(\g)[[\hbar]]$ is the standard inclusion. 

Using these isomorphisms, $U_\hbar(\l) \subset U_\hbar(\g)$ identifies with 
an injection $U(\l)[[\hbar]] \hookrightarrow U(\g)[[\hbar]]$, such that the 
composed map $\h\subset U(\l)[[\hbar]] \hookrightarrow U(\g)[[\hbar]]$ is the 
standard inclusion. The
theorem now follows from the fact that such an injection is necessarily of the
form $\on{Ad}(u) \circ i$, where $u\in U(\g)^\h[[\hbar]]$ is invertible and $i$
is the inclusion $U(\l)[[\hbar]] \subset U(\g)[[\hbar]]$. This last statement is
a consequence of: 

\begin{lemma}
Any derivation $\delta : \l\to U(\g)$, such that $\delta_{|\h} = 0$, is inner,
i.e., of the form $x\mapsto [v,x]$, where $v\in U(\g)^\h$. 
\end{lemma}

{\em Proof of Lemma.} Since $\l = \z\oplus \l'$, we have an injection
$H^1(\l,U(\g)) \subset H^1(\z,U(\g)) \oplus H^1(\l',U(\g))$. Since $\l'$ is
semisimple and $U(\g)$ is a semisimple $\l'$-module, the last cohomology group
vanishes. It follows that there exists $w\in U(\g)$ and $\delta_0 \in
H^1(\z,U(\g)) = \z^* \otimes U(\g)^\z$, such that $\delta(z+l') = \delta_0(z) +
[w,l']$ for any $z\in \z$ and $l'\in \l'$. Now $\delta_{|\h} = 0$, so
$\delta_{|\z} = 0$, so $\delta_0 = 0$. On the other hand, $\delta$ is a
derivation of $\l$, hence $0 = \delta([z,l']) = [z,\delta(l')] = [z,[w,l']] =
[[z,w],l']$ for any $z\in \z$, $l'\in \l'$, i.e., $[z,w]\in U(\l)^{\l'}$ for any
$z\in \z$. 

The decomposition of $U(\g)$ as a $\l'$-module has the form 
$U(\g) = U(\g)^{\l'} \oplus \bigoplus_{\rho\in \on{Irr}(\l'), 
\rho \neq {\mathbf 1}}
M_\rho \otimes V_\rho$, where $\on{Irr}(\l')$ is the set of irreducible finite
dimensional $\l'$-modules, and ${\mathbf 1}$ is the trivial $\l'$-module. 
Let $w = w_0 + \sum_{\rho\in \on{Irr}(\l'), \rho\neq {\mathbf 1}} w_\rho$
be the corresponding decomposition of $w$. Then if $z\in \z$, we have
$[z,w_0]\in U(\g)^{\l'}$, and $[z,w_\rho] \in M_\rho \otimes V_\rho$, so
$[z,w_\rho] = 0$ for any $\rho\neq {\mathbf 1}$. We then set $v := \sum_{\rho
\in \on{Irr}(\l'), \rho \neq {\mathbf 1}} w_\rho$. We have for any $l'\in \l'$,
$\delta(l') = [w,l'] = [w-w_0,l'] = [v,l']$ since $w_0\in U(\g)^{\l'}$, 
and $\delta(z) = 0 = \sum_{\rho\in \on{Irr}(\l'), \rho\neq {\mathbf 1}}
[w_\rho,z] = [v,z]$ for any $z\in \z$. Hence $\delta(l) = [v,l]$ for any $l\in
\l$. Since $\delta_{|\h} = 0$, we also have $v\in U(\g)^\h$. 
\hfill \qed \qed \medskip 

We will now assume that $J_0 = u^{-1} * J_{\on{EK}}^\Phi(r_0)$,  
therefore 
$U(\l)[[\hbar]]$ is a QUE subalgebra of $\big( U(\g)[[\hbar]], 
\on{Ad}(J_0^{-1}) \circ \Delta_0\big)$, which is the image of 
$U_\hbar(\l) \subset U_\hbar(\g)$ under $U_\hbar(\g)
\simeq (U(\g)[[\hbar]],\Ad(J_0^{-1}) \circ \Delta_0)$.

\subsection{Deformations of $\sigma_\l^\g(\chi)$} \label{sect:rigidity}

Proposition \ref{isom:vars} can be generalized to
a scheme-theoretic setup as follows. 
If $R$ is a ring containing $\CC$, define $\cP(R)$ as the set of 
all $(c(\theta))_{\theta\in\Delta_+(\g;\l)}
\in R^{\Delta_+(\g;\l)}$, such that for any $\theta,\theta'\in 
\Delta_+(\g;\l)$ such that $\theta + \theta'\in \Delta_+(\g;\l)$, 
we have $c(\theta + \theta') \big( c(\theta) + c(\theta') \big) 
= c(\theta)c(\theta') + {1\over 4}$. Define $\cP_0(R)$ as the subset of 
$\cP(R)$ of all $(c(\theta))_{\theta\in\Delta_+(\g,\l)}$, 
such that for any $\theta\in
\Delta_+(\g;\l)$, $c(\theta) \pm {1\over 2}$ is invertible in $R$. 

Define $\chi_\theta : (R^\times)^{\bar I} \to R^\times$ by 
$\chi_\theta(\tt) = \prod_{s\in \bar I} t_s^{n_s(\theta)}$ if 
$\tt = (t_s)_{s\in\bar I}$, and define
$X(R) \subset (R^\times)^{\bar I}$ as the set of all 
$\tt$, such that $\chi_\theta(\tt)-1$ is invertible in $R$ for any $\theta 
\in \Delta_+(\g;\l)$. 

Then the isomorphism $\iota$ generalizes to an isomorphism 
between $\cP_0(R)$ and $X(R)$. 

When $R$ is the ring $\CC[\eps]/(\eps^2)$ of dual numbers, 
$\cP_0(R)$ is the preimage of $\cP_0(\CC)$ under the natural map 
$\cP(R)\to \cP(\CC)$. In the same way, $X(R)$ is the preimage of
$(\CC^\times)^{\bar I} \setminus \cup_{\theta\in\Delta_+(\g;\l)}
\on{Ker}(\chi_\theta)$ under $(R^\times)^{\bar I} \to (\CC^\times)^{\bar I}$. 
 
We denote by $(a,b)\mapsto \on{CYB}(a,b)$ the bilinear map
derived from the quadratic map $a\mapsto \on{CYB}(a)$.  

\begin{proposition} \label{rho:rigid}
If $\tt \in (\CC^\times)^{\bar I} \setminus 
\cup_{\theta\in\Delta_+(\g;\l)} \on{Ker}(\chi_\theta)$, 
and $\rho'\in \wedge^2(\u)^\l$ is such that 
$\on{CYB}(\rho_{\tt},\rho')= 0$ in $\wedge^3(\g/\l)$, 
then there exists a unique family $(\tau_s)_{s\in \bar I}$
in $\CC^{\bar I}$, such that 
$$
\rho' = \sum_{s\in \bar I}
{{\partial \rho_{\tt}}\over {\partial t_s}} \tau_s.
$$  
\end{proposition}

{\em Proof.} $\rho_{\tt} + \eps \rho'$ is an element of 
$\cP_0(R)$, where $R = \CC[\eps]/(\eps^2)$. Its image by the 
generalization of $\iota$ is a family of $(R^\times)^{\bar I}$, 
which has the form $(t_s + \eps \tau_s)_{s\in \bar I}$, with 
$\tau_s\in\CC$. 
\hfill \qed \medskip 

\begin{corollary} \label{rigidity:r:l:g}
If $\chi\in \z^*$ is a regular character of $\l$, and $\sigma'\in
\wedge^2(\u)^\l$ is such that $\on{CYB}(\sigma_\l^\g(e^\chi),\sigma') = 0$
in $\wedge^3(\g/\l)$, then there exists a unique element 
$\chi'\in \z^*$, such that 
$$
\sigma' = {{\on{d}}\over{\on{d}\eps}}_{|\eps = 0} 
\sigma_\l^\g(e^{\chi + \eps \chi'}). 
$$
\end{corollary}

{\em Proof.} We have a sequence of coverings $\z^* \stackrel{\on{exp}}{\to} 
Z^*_{\on{alg}} \to (\CC^\times)^{\bar I}$, with composition $\chi \mapsto 
\tt(\chi)$, such that $t_s = e^{-(\chi,\bar\alpha_s)}$. Then 
$\sigma_\l^\g(e^{\chi}) = \rho_{\tt(\chi)}$. The statement is then a 
consequence of Proposition \ref{rho:rigid}. 
\hfill \qed \medskip 

\subsection{Comparison with \cite{EE2}} \label{sect:compare}

Let us reexpress the algebra $(\CC[[G/L]]_\hbar,*)$ defined in 
Section \ref{sect:qP}. The isomorphism of $U_\hbar(\g)$
with $(U(\g)[[\hbar]],\Ad(J_0^{-1}) \circ \Delta_0)$ 
allows to express the product of $\CC[[G]]_\hbar$ as follows: 
$\CC[[G]]_\hbar$ is isomorphic to $\CC[[G]][[\hbar]]$, and 
$f * g = m(\LL^{\otimes 2}(J_0) \RR^{\prime \otimes 2}(J_0^{-1})
(f\otimes g))$. 

Moreover, by virtue of the remark following Theorem \ref{thm:levi}, 
the image of $\CC[[G/L]]_\hbar$ under $\CC[[G]]_\hbar \simeq \CC[[G]][[\hbar]]$
is $\CC[[G/L]][[\hbar]]$. Under this isomorphism, the product of 
$(\CC[[G/L]]_\hbar,*)$ is transported to the product 
$$
f_1 *_{\chi} f_2 = m_0 \big( \LL^{\otimes 2}(J_0) 
\RR^{\prime \otimes 2}( S^{\otimes 2}(\Xi(\chi)) J_0^{-1})
(f_1 \otimes f_2)\big),  
$$
for any $f_1,f_2 \in \CC[[G/L]][[\hbar]]$ (here we view 
$S^{\otimes 2}(\Xi(\chi))$ as an element of $U(\g)^{\otimes 2}[[\hbar]]$). 

\begin{lemma} \label{lemma:6:7}
If $A(\chi)$ is either of the elements $J(\chi)^{-1}$ and 
$S^{\otimes 2}(\Xi(\chi)) J_0^{-1}$ of $U(\g)^{\otimes 2}[[\hbar]]$, 
then $A(\chi)$ satisfies the identities:  

(i) the image of $A(\chi)$ in $(U(\g)/\l U(\g))^{\otimes 2}[[\hbar]]$
is $\l$-invariant; 

(ii) the image of $A(\chi)^{1,2} A(\chi)^{12,3} - A(\chi)^{2,3} A(\chi)^{1,23}
\Phi_\g$ is zero under the projection $U(\g)^{\otimes 3}[[\hbar]]
\to (U(\g)/\l U(\g))^{\otimes 3}[[\hbar]]$

(iii) the image of $A(\chi) - (1 + \hbar r_\l^\g(\chi))$ is $O(\hbar^2)$
under the projection $U(\g)^{\otimes 2} \to 
(U(\g)/\l U(\g))^{\otimes 2}[[\hbar]]$. 
\end{lemma}

{\em Proof.} The arguments are the same as those used to prove that 
$*$ and $*_\chi$ are associative products on $\CC[[G/L]]$. 
\hfill \qed \medskip

\begin{proposition}
Let $A(\chi)$, $B(\chi)$ be analytic functions $U_{\z^*} \to U(\g)^{\otimes
2}[[\hbar]]$, satisfying the conditions (i), (ii), (iii) of Lemma
\ref{lemma:6:7}. Then there exists an analytic function $u(\chi) : 
U_{\z^*} \to 1 + \hbar U(\g)[[\hbar]]$, and analytic function $\chi_n : 
U_{\z^*} \to \z^*$ ($n\geq 1$), such that the image of 
$$
A(\chi) - \big( u(\chi) \otimes u(\chi) \big) 
B(\chi'(\chi)) \Delta_0(u(\chi))^{-1}
$$  
in $(U(\g)/\l U(\g))^{\otimes 2}[[\hbar]]$ is zero. Here we set 
$\chi'(\chi) = \chi + \sum_{n\geq 1} \hbar^n \chi_n(\chi)$. 
\end{proposition}

{\em Proof.} There is a unique map 
$$
\dd : 
\big( (U(\g)/\l U(\g))^{\otimes 2} \big)^\l[[\hbar]]
\to \big( (U(\g)/\l U(\g))^{\otimes 3} \big)^\l[[\hbar]], 
$$
taking the class of $A$ to the class of $A^{1,2} A^{12,3} - A^{2,3} A^{1,23}
\Phi_\g$, where $A\in U(\g)^{\otimes 2}[[]\hbar]$ is a representative of $A$. 
The group $1 + \hbar (U(\g)/\l U(\g))^\l[[\hbar]]$ acts on 
$\big( (U(\g)/\l U(\g))^{\otimes 2} \big)^\l[[\hbar]]$ as follows: $\bar u *
\bar A$ is the class of $(u\otimes u)A \Delta_0(u)^{-1}$, where $u\in
U(\g)[[\hbar]]$, $A$ are representatives of $\bar u$, $\bar A$. Then 
$\dd(\bar u * \bar A)$ is the class of 
$u^{\otimes 3} \wt{\dd(\bar A)} \Delta_0^{(2)}(u)^{-1}$, where 
$\wt{\dd(\bar A)}$ is a representative of $\dd(\bar A)$; in particular, 
$\dd(\bar A) = 0$ iff $\dd(\bar u * \bar A) = 0$. 

The classes $\bar A,\bar B$ are analytic functions  
$U_{\z^*} \to \big( (U(\g) / \l U(\g) )^{\otimes 2}\big)^\l[[\hbar]]$, 
such that $\dd(\bar A(\chi)) = \dd(\bar B(\chi)) = 0$. Assume that we have
constructed $\bar u_p(\chi) : U_{\z^*} \to (U(\g)/\l U(\g))^\l$
and $\chi_p(\chi) : U_{\z^*} \to \z^*$, ($p = 1 ,\ldots,n-1$), such that 
if $\chi'_{(n-1)}(\chi) = \chi + \sum_{p=1}^{n-1} \hbar^p \chi'_p(\chi)$
and $u_{(n-1)}(\chi) = 1 + \sum_{p=1}^{n-1} \hbar^p u_p(\chi)$, then 
\begin{equation} \label{*}
\bar A(\chi) - u_{(n-1)}(\chi) * \bar B(\chi'_{(n-1)}(\chi)) 
= O(\hbar^{n-1}). 
\end{equation}
The co-Hochschild differential $\on{d}_{\on{co-Hoch}}(A) = A^{12,3} - A^{1,23}
- A^{2,3} + A^{1,2}$ induces a differential 
$((U(\g)/\l U(\g))^{\otimes 2})^{\l} \to ((U(\g)/\l U(\g))^{\otimes 3})^{\l}$, 
and identifies with the standard co-Hochschild differential 
$(S^\cdot(\u)^{\otimes 2})^\l \to (S^\cdot(\u)^{\otimes 3})^\l$. 

Let $\bar C(\chi)$ be the class modulo $\hbar$ of 
$\hbar^{1-n} \big( \bar A(\chi) - 
u_{(n-1)}(\chi) * \bar B(\chi'_{(n-1)}(\chi)) \big)$. Then 
$\on{d}_{\on{co-Hoch}}(\bar C(\chi)) =0$. So $\bar C(\chi) = 
\on{d}_{\on{co-Hoch}}(u_{(n)}(\chi)) + \rho(\chi)$, where 
$u_{(n)}(\chi)$, $\rho(\chi)$ are analytic functions 
$U_{\z^*} \to (U(\g) / \l U(\g))^\l$ and $U_{\z^*} \to \wedge^2(\g/\l)^\l$. 

Since $r_\l^\g(\chi)$ is $\l$-invariant, the map $\rho \mapsto
\on{CYB}(r_\l^\g(\chi),\rho)$ induces a linear map $\wedge^2(\g/\l)^\l
\to \wedge^3(\g/\l)^\l$. Proposition \ref{rigidity:r:l:g} says that the 
kernel of this map consists of the ${{\on{d}}\over{\on{d}\eps}}_{|\eps=0}
r_\l^\g(\chi + \eps \chi')$, $\chi'\in\z^*$. 
  
Now the vanishing of the coefficient of $\hbar^n$
in $\dd(\bar A(\chi)) = \dd(u_{(n-1)}(\chi) * \bar B(\chi'_{(n-1)}(\chi)))$
yields $\on{CYB}(r_\l^\g(\chi),\rho(\chi)) 
+ \on{d}_{\on{co-Hoch}}(\bar D(\chi)) = 0$, where $\bar D(\chi)$ is the
coefficient of $\hbar^n$ in $\bar A(\chi) - \bar u_{(n-1)} * 
\bar B(\chi'_{(n-1)}(\chi))$. It follows that 
$\on{CYB}(r_\l^\g(\chi),\rho(\chi)) =0$, so 
$\rho(\chi) = {{\on{d}}\over{\on{d}\eps}}_{|\eps=0}
r_\l^\g(\chi + \eps \chi_n(\chi))$, $\chi_n(\chi)\in\z^*$. 
So we have determined $u_{n}(\chi)$ and $\chi_n(\chi)$, such that 
(\ref{*}) holds at order $n$. 
\hfill \qed \medskip 

\begin{corollary}
$f\mapsto \RR'_{u(\chi)}(f)$ is an isomorphism between the products 
$*$ and $*_\chi$ on $\CC[[G/L]][[\hbar]]$.  
\end{corollary}

\section{Quantization of $G/L \hookrightarrow \g^*$, 
$G/L\hookrightarrow G^*$ and Verma modules} \label{sect:7}

Let $\g$ be a semisimple Lie algebra, $\l\subset \g$
be a Levi subalgebra, $\g = \l\oplus \u_+\oplus \u_-$ be a
decomposition of $\g$. Let $\chi\in \l^*$ be a character
of $\l$; we view $\chi$ as an element of $\g^*$ by defining 
it to be zero on $\u_\pm$. We assume that $\chi$ is regular, 
which means that the centralizer of $\chi$ is $\l$. 

Let $G_{\alg}$ be the adjoint group of $\g$, $L_\alg\subset G_\alg$ the 
subgroup corresponding to $\l$. Then the orbit 
$G_\alg\chi\subset \g^*$ identifies with $G_\alg/L_\alg$. We have a
sequence of maps 
\begin{equation} \label{seq:maps}
G_\alg \to G_\alg/L_\alg = G_\alg\chi 
\subset p^{-1}(p(\chi)) \subset \g^*; 
\end{equation}
here $p : \g^* \to \g^*// G_\alg = \h^* /W$ is the natural projection. 
$p^{-1}(p(\chi))$ is a union of coadjoint orbits on $\g^*$, and 
$G_\alg\chi$ is its only closed orbit (\cite{Ko}). 
In (\ref{seq:maps}), all varieties and maps except the first one
are Poisson.  

We now construct the quantum version of (\ref{seq:maps}). 

\subsection{Poisson algebras}

The function algebra of $G_\alg$ is $\CC[G] = \oplus_{V\in \Irr}
V^* \otimes V$. Here $\Irr$ is the set of simple objects in the 
subcategory of the tensor category of finite dimensional 
$\g$-modules, generated by the adjoint representation
\footnote{Irr may also be described as the set of simple $\g$-modules
$V$ with highest weight in the root lattice, or such that $V[0] \neq 0$}. 
Its product is defined by $(\xi\otimes v) * (\eta\otimes w) = 
\varpi((\xi\otimes \eta) \otimes (v\otimes w))$. Here $\xi\otimes v
\in V^* \otimes V$, $\eta\otimes w \in W^* \otimes W$, and 
$\varpi : (V\otimes W)^* \otimes (V\otimes W) \to \oplus_{Z\in\Irr}
Z^* \otimes Z$ is the composition of the map induced from 
$V\otimes W = \oplus_{Z\in \Irr}  Z\otimes M_{V,W}^Z$ with the maps 
$\on{tr} : M_{V,W}^Z \otimes (M_{V,W}^Z)^* \to \CC$
and the zero map on $M_{V,W}^Z \otimes (M_{V,W}^{Z'})^*$ if
$Z'\neq Z$. 

The map $\CC[G] \subset \CC[[G]] = U(\g)^*$ takes $\xi\otimes v$ to 
$T\mapsto \langle \xi, \rho_V(T)(v) \rangle$. 

The function algebra of $G_\alg/L_\alg$ is $\CC[G/L]
= \CC[G]^\l = \oplus_{V\in\Irr} V^* \otimes V^\l$. Its Poisson 
structure is defined by 
$$
\{\xi\otimes v,\eta\otimes w\} = 
\varpi \big( (\xi\otimes \eta) \otimes 
\rho_{V\otimes W}(r_\l^\g(\chi))(v\otimes w)\big) . 
$$
This is a Poisson subalgebra of $\CC[[G/L]] = (U(\g)/U(\g)\l)^*$. 

The dual of the map $G_\alg/L_\alg\to \g^*$
is the map $\CC[\g^*] = S^\cdot(\g) \to \CC[G/L]$, defined 
as the unique algebra morphism taking 
$x\in \g$ to $x\otimes \chi \in \g \otimes (\g^*)^\l$. 
This is a Poisson algebra morphism, which factors through 
a morphism 
$$
S^\cdot(\g) \to S^\cdot(\g)\otimes_{Z(S^\cdot(\g))} 
\CC_{p(\chi)} \to \CC[G/L]. 
$$ 
Here $Z(S^\cdot(\g)) = S^\cdot(\g)^\g$ is the Poisson center 
of $S^\cdot(\g)$, and $p(\chi)$ is the character of 
$Z(S^\cdot(\g))$ obtained by restricting the character 
$S^\cdot(\g) \to \CC$ induced by $\chi$. 

Indeed, an element $P\in S^l(\g)^\g$ is mapped by 
$P \otimes \chi^l$; $\varpi$ then maps it to 
$P(\chi) 1\otimes 1$. 

The composed map $S^\cdot(\g) \to 
(U(\g)/U(\g)\l)^*$ is given by the 
algebra-coalgebra pairing 
$S^\cdot(\g) \otimes (U(\g)/U(\g)\l) \to \CC$,  
taking $x\otimes T$ to $\chi \big( \ad(S(T))(x)\big)$
(here $S$ is the antipode, $\ad$ is the adjoint action of 
$U(\g)$ on $S^\cdot(\g)$, and $\chi : S^\cdot(\g) \to \CC$
is the character corresponding to $\chi$). 

$G_\alg/L_\alg$ is a subvariety of $\g^*$, and its function ring is 
$\CC[G/L]$, hence the map 
$S^\cdot(\g) \to \CC[G/L]$ is surjective.

\subsection{Quantization of $G/L \to G_\alg/L_\alg \to \g^*$}

In \cite{EE2}, we constructed an $\l$-invariant 
element $J := J_\l^\g(\hbar^{-1}\chi)
\in U(\u_+) \otimes U(\p_-)[[\hbar]]$, such that 
$J = 1 + O(\hbar)$, $\hbar^{-1}(J - J^{2,1}) = r_\l^\g(\chi)$ 
mod $\hbar$, and $\big( U(\g)/U(\g)\l \big) [[\hbar]]$, equipped with 
$\Delta(f) := \Delta_0(f)J^{2,1}$ is a coassociative 
coalgebra (we set $\p_\pm = \l\oplus \u_\pm$). Its dual algebra 
$\big( U(\g)/U(\g)\l \big)^*[[\hbar]]$ is a quantization of 
$\CC[[G/L]]$, which we denote by $\CC[[G/L]]_\hbar$. 

Define $\CC[G/L]_\hbar$ as $\oplus_{V\in \Irr}
V^* \otimes V^\l[[\hbar]]$, equipped with the product 
$$
(\xi\otimes v) * (\eta\otimes w) = \varpi((\xi\otimes \eta)\otimes 
\rho_{V\otimes W}(J)(v\otimes w)). 
$$ 
Then $\xi\otimes v \mapsto (T\mapsto \langle \xi, \rho_V(T)(v) \rangle)$
is an algebra embedding $\CC[G/L]_\hbar \subset \CC[[G/L]]_\hbar$. 

Define $\CC[\g^*]_\hbar$ as the enveloping algebra of $\g[[\hbar]]$, 
equipped with the bracket $[x,y]_\hbar := \hbar [x,y]$. Then 
$\CC[\g^*]_\hbar$ is a flat deformation of $\CC[\g^*] = S^\cdot(\g)$. 
$\CC[\g^*]_\hbar$ is a subalgebra of $U(\g)[[\hbar]]$ under 
$x\mapsto \hbar x$. 

If $\chi_0\in \l^*$ is a character, define $\overline\chi_0 : U(\g)
\to \CC$ as the unique linear map such that $\overline\chi_0(x_+ x_0 x_-)
= \eps(x_+)\eps(x_-) \chi_0(x_0)$, where $x_\pm \in U(\u_\pm)$
and $x_0 \in U(\l)$. This definition extends to the case when 
$\chi_0\in \l^*((\hbar))$. In particular, 
$\overline{\hbar^{-1}\chi}$ is a character $U(\g)[[\hbar]] \to \CC((\hbar))$; 
one checks that it restricts to a character $\wt\chi : 
\CC[\g^*]_\hbar \to \CC[[\hbar]]$. 

\begin{proposition}
There is a unique pairing 
$$
\big( U(\g)/U(\g)\l \big) \otimes \CC[\g^*]_\hbar \to \CC[[\hbar]], 
$$ 
taking $[T] \otimes x$ to $\wt\chi(\on{ad}(S^{-1}(T))(x))$. It induces an 
algebra morphism $\CC[\g^*]_\hbar^{\on{op}} \to \CC[[G/L]]_\hbar$. 
This morphism factors through a morphism 
$\CC[\g^*]_\hbar^{\on{op}} \to \CC[G/L]_\hbar 
\subset \CC[[G/L]]_\hbar$, equivariant under the adjoint action of 
$U(\g)$ on $\CC[\g^*]_\hbar$ and its left action on the two 
other algebras. 
\end{proposition} 

{\em Proof.} The fact that the pairing is well-defined follows from the 
$\l$-invariance of $\wt\chi$. Let us denote it by $[T]\otimes x 
\mapsto (T,x)$ and check that it is a coalgebra-algebra pairing. 
We first prove that if $x,y\in \CC[\g^*]_\hbar$, then 
\begin{equation} \label{id:chi}
\wt\chi(xy) = \sum_i \wt\chi(\on{ad}(S^{-1}(\alpha_i))(x)) 
\wt\chi(\on{ad}(S^{-1}(\beta_i))(y)),  
\end{equation}
where $J = \sum_i \alpha_i \otimes \beta_i$.  
Recall that $J = \sum_i a_i \otimes S(b_i)S(\ell_i^{(2)})
\wt\chi(\ell_i^{(1)})$, where 
$K = \sum_i a_i \otimes b_i \otimes \ell_i$ is the element of 
$U(\u_+) \otimes U(\u_-) \wh\otimes \wh{U(\l)}$ such that 
the identity $\sum_i H(xa_i)\ell_i H(b_i y) = H(xy)$
for any $x,y\in U(\g)$ holds; here $H : U(\g)\to U(\l)$
is the linear map taking $x_+ x_0 x_-$ to $\eps(x_+)\eps(x_-)
x_0$, where $x_\pm\in U(\u_\pm),x_0\in U(\l)$; we again denote by $\wt\chi$
the unique extension of $\wt\chi$ to a character of $\wh{U(\l)}$
(it exists because $\chi$ is nondegenerate). 

Using the identities $\wt\chi(x_+x)
= \eps(x_+)\wt\chi(x)$, $\wt\chi(xx_-) = \wt\chi(x)\eps(x_-)$
for $x\in U(\g)$, $x_\pm\in U(\u_\pm)$, and the $\l$-invariance of 
$\wt\chi$, we identify
the r.h.s. of (\ref{id:chi}) with $\sum_i \wt\chi(x\alpha_i)
\wt\chi(S^{-1}(\beta_i)y)$. This is 
$\sum_i \wt\chi(x a_i) \wt\chi(\ell_i^{(1)}) 
\wt\chi(\ell_i^{(2)} b_i y)$, i.e., $\sum_i \wt\chi(xa_i)
\wt\chi(\ell_i)\wt\chi(b_iy)$ (as $\chi$ is a regular character of $\l$, 
we have $(\wt\chi \otimes \wt\chi) \circ \Delta_r = \wt\chi$, where 
$\Delta_r : \wh{U(\l)} \to \wh{U(\l)} \wh\otimes U(\l)$ is the 
"right coproduct" of $\wh{U(\l)}$), 
i.e., $\wt\chi(\sum_i H(xa_i)\ell_i H(b_i y))$, i.e., $\wt\chi(H(xy))
= \wt\chi(xy)$. This proves (\ref{id:chi}). 

If $T\in U(\g)$ and $x,y\in \CC[\g^*]_\hbar$, we have 
\begin{align*}
(T,xy) & = \wt\chi \big( \on{ad}(S^{-1}(T^{(2)}))(x) 
\on{ad}(S^{-1}(T^{(1)}))(y)\big) 
\\ & 
= \sum_i \wt\chi 
\big( \on{ad}(S^{-1}(\alpha_i)) \on{ad}(S^{-1}(T^{(2)}))(x) \big) 
\wt\chi
\big( \on{ad}(S^{-1}(\beta_i)) \on{ad}(S^{-1}(T^{(1)}))(y) \big) 
\on{\ (by\ (\ref{id:chi}))}
\\ & 
= \sum_i (T^{(2)}\alpha_i,x)(T^{(1)}\beta_i,y). 
\end{align*}

If $x\in U(\g)$ is an element of degree $d$, then 
the map $T\mapsto \wt\chi(\on{ad}(S^{-1}(x)))$ is the matrix 
coefficient of $x \otimes \wt\chi_{|U(\g)_{\leq d}}
\in U(\g)_{\leq g} \otimes (U(\g)^*_{\leq d})^\l$; 
so the image of $\CC[\g^*]_\hbar^{\on{op}} \to \CC[[G/L]]_\hbar$
is contained in $\CC[G/L]_\hbar$. 
\hfill \qed \medskip 

\subsection{Relation with generalized Verma modules}

Let $(M_\chi,\pi_\chi)$ be the $U(\g)$-module 
$\on{Ind}_{\p_-}^\g(\CC((\hbar))_{\wt\chi})$, where $\wt\chi$
is the character $\p_- \to \l \stackrel{\hbar^{-1}\chi}{\to} \CC((\hbar))$.  

%Let us denote by $\End(M_\chi)_{\on{ft}}$ the space of 
%finite type endomorphisms of $M_\chi$; this is the direct 
%sum of all finite dimensional $U(\g)$-submodules of $\End(M_\chi)$
%(on which $U(\g)$ acts by the adjoint action 
%x * \phi = \sum \pi_\chi(x^{(1)}) \circ \phi 
%\circ \pi_\chi(S(x^{(2)}))$). 

\begin{proposition}
The morphism $\CC[\g^*]_\hbar \to \End(M_\chi)[[\hbar]]$ induced by 
$\pi_\chi$ factors through $\CC[\g^*]_\hbar 
\stackrel{\alpha}{\twoheadrightarrow}
\CC[G/L]_\hbar \stackrel{\beta}{\hookrightarrow} \End(M_\chi)[[\hbar]]$,
where $\alpha$ is surjective and $\beta$ is injective. It follows that 
$\CC[G/L]_\hbar \simeq \CC[\g^*]_\hbar / (\CC[\g^*]_\hbar \cap 
\on{Ker}(\pi_\chi))$. 
\end{proposition}

{\em Proof.} Let us construct an algebra morphism 
$\beta : \CC[G/L]_\hbar \to \End(M_\chi)[[\hbar]]$. 
If $V\in \Irr$, we define a linear map $V^* \otimes V^\l \to 
\End(M_\chi)[[\hbar]]$ by dualizing the linear map 
$M_\chi \to M_\chi \otimes V \otimes (V^\l)^*[[\hbar]]$
equal to $\sum_{b\in B} \Phi_b \otimes b^*$. Here $B$ is a 
basis of $V^\l$, $(b^*)$ is the dual 
basis of $(V^\l)^* \simeq (V^*)^\l$, and 
$\Phi_b : M_\chi \to M_\chi\otimes V[[\hbar]]$
is the intertwiner with expectation value $b$. 
Using the fact that $\Phi_b(1_\chi) 
= \sum_i \alpha_i 1_\la \otimes \beta_i b$, one can prove that 
$\CC[G/L]_\hbar \to \End(M_\chi)[[\hbar]]$ is an algebra morphism. 

Let us prove that $\beta$ is injective. Recall that $\beta$ takes
$v\otimes v^* \in V^\l \otimes V^*$ to the endomorphism 
$m\mapsto (\id\otimes v^*)\circ \Phi_v(m)$ of $\End(M_\chi)[[\hbar]]$. 
Therefore, $\beta$ extends to a linear map $\beta' : 
(U(\g)/U(\g)\l)^* \to \End(M_\chi)[[\hbar]]$, taking
$\xi$ to the endomorphism $x_+ 1_\chi \mapsto \sum_i \langle \xi, 
x_+^{(1)} \beta_i \rangle x_+^{(2)} \alpha_i 1_\chi$
(here $x_+ \in U(\u_+)$). 

We have a linear isomorphism $(U(\g)/ U(\g)\l)^* \simeq 
U(\u_+)^* \bar\otimes U(\u_-)^*$, whose inverse takes 
$\xi_+ \otimes \xi_-$ to the form $x_+ x_- \mapsto 
\langle \xi_+, x_+\rangle \langle \xi_-, x_-\rangle$, where
$x_\pm \in U(\u_\pm)$. 

Let $\pi : U(\p_+) \to U(\u_+)$ be the linear map taking 
$x_+ x_0$ to $x_+ \eps(x_0)$, where $x\in U(\u_+)$, $x_0 \in U(\l)$. 

The map $\beta'' : U(\u_+)^* \bar\otimes U(\u_-)^* \to 
\End(M_\chi)[[\hbar]]$ induced by $\beta'$ takes $\xi_+ \otimes \xi_-$
to the endomorphism 
$$
x_+ 1_\chi \mapsto \sum_i \langle \xi_+,
x_+^{(1)} \rangle x_+^{(2)} \langle \xi_-, \pi(\beta_i)\rangle
\alpha_i 1_\chi . 
$$ 
Assume that $\sum_\alpha \xi_{+,\alpha} \otimes \xi_{-,\alpha}$
is in $\on{Ker}(\beta'')$. Then 
\begin{equation} \label{ker:beta''}
\forall x\in U(\u_+), \; 
\sum_{i,\alpha} 
\langle \xi_{+,\alpha}, x_+^{(1)} \rangle x_+^{(2)} 
\langle \xi_{-,\alpha}, \pi(\beta_i)\rangle
\alpha_i =0. 
\end{equation}
Since $\sum_i \alpha_i \otimes 
\pi(\beta_i)$ coincides with $\on{exp}(r_\l^\g(\chi))$ up to 
lower degree terms,  the linear map $U(\u_-)^* \to U(\u_+)$, 
$\xi\mapsto \sum_i \langle \xi, \pi(\beta_i)\rangle \alpha_i$
is injective. 

Then using (\ref{ker:beta''}) with $x=1$, we find $\sum_\alpha 
\langle \xi_{+,\alpha}, x\rangle \xi_{-,\alpha} = 0$. 
Let us prove by induction
on $d$ that for any $x\in U(\u_+)_{\leq d}$, $\sum_\alpha
\langle \xi_{+,\alpha}, x\rangle \xi_{-,\alpha} =0$. 
Assume that this holds at order $d$ and let $x\in U(\u_+)_{\leq d+1}$. 
Then the induction hypothesis implies 
$\sum_{\alpha,i} \langle \xi_{+,\alpha}, x\rangle \langle \xi_{-,\alpha},
\pi(\beta_i) \rangle \alpha_i =0$. The injectivity of 
$\xi\mapsto \sum_i \langle \xi, \pi(\beta_i)\rangle \alpha_i$
then implies that 
$\sum_{\alpha} \langle \xi_{+,\alpha}, x\rangle \xi_{-,\alpha}
=0$, as wanted. Therefore $\sum_\alpha \xi_{+,\alpha} \otimes \xi_{-,\alpha}=0$.
Hence $\beta$ is  injective. 

Let us now show that $\beta \circ \alpha = \pi_{\chi | \CC[\g^*]_\hbar}$. 
Let $x\in U(\g)$, then $\alpha(x)$ is the form 
$U(\g)/U(\g)\l\to \CC$, $T\mapsto \wt\chi(\on{ad}(S^{-1}(T))(x))$. 
Then $\beta''(\alpha(x))$ is the endomorphism of $M_\chi[[\hbar]]$
taking $x_+ 1_\chi$ to 
\begin{align} \label{action:1}
\nonumber \sum_i \langle \alpha(x), x_+^{(1)} \beta_i\rangle 
x_+^{(2)} \alpha_i 1_\chi
& = 
\sum_i \wt\chi\big( S^{-1}(\beta_i^{(2)}) S^{-1}(x_+^{(2)}) x x_+^{(1)}
\beta_i^{(1)}\big) x_+^{(3)} \alpha_i 1_\chi
\\ & 
= 
\sum_i \wt\chi\big( S^{-1}(\beta_i) S^{-1}(x_+^{(2)}) x x_+^{(1)} \big) 
x_+^{(3)} \alpha_i 1_\chi. 
\end{align}
Now we have $\sum_i \alpha_i \wt\chi(S^{-1}(\beta_i) X_+ X_0 X_-) = 
X_+ \wt\chi(X_0) \eps(X_-)$ for $X_\pm \in U(\u_\pm)$ and $X_0\in U(\l)$. 
Hence if $S^{-1}(x_+^{(2)}) x x_+^{(1)} = \sum_\alpha X_{+,\alpha} X_{0,\alpha}
X_{-,\alpha}$, then 
$$
(\ref{action:1}) = x_+^{(3)} \big(
\sum_\alpha X_{+,\alpha} \wt\chi(X_{0,\alpha})
\eps(X_{-,\alpha}) \big) 1_\chi. 
$$
On the other hand, 
\begin{align*}
xx_+ 1_\chi & = x_+^{(3)} S^{-1}(x_+^{(2)}) x x_+^{(1)} 1_\chi
= x_+^{(3)} (\sum_\alpha X_{+,\alpha} X_{0,\alpha} X_{-,\alpha}) 1_\chi
\\ & = x_+^{(3)} \sum_\alpha X_{+,\alpha} \wt\chi(X_{0,\alpha}) 
\eps(X_{-,\alpha}) 1_\chi,  
\end{align*}
which is (\ref{action:1}). Hence $(\ref{action:1}) = xx_+ 1_\chi$, so 
$\beta \circ \alpha = \pi_{\chi|\CC[\g^*]_\hbar}$.  
\hfill \qed \medskip 

\subsection{PL versions}

Let $\tt, \chi,\chi_{L^*}$ be as in Lemma \ref{lemma:4:11}. The dressing
orbit of $\chi_{L^*} \in G^*_{\on{alg}}$ is $G_{\alg} \chi_{L^*} \simeq
G_\alg/L_\alg$. Moreover, the Poisson structure of $G_\alg\chi_{L^*}$ identifies
with $\{-,-\}_{\chi_{L^*}}$. We have therefore a Poisson $G$-space embedding 
$G_\alg/L_\alg \simeq G_\alg \chi_{L^*} \hookrightarrow G_\alg^*$. As before, we
can construct a morphism $\CC[G^*]_\hbar \to \CC[G/L]_\hbar$
quantizing this embedding; here $\CC[G^*]_\hbar$ is a subalgebra of 
$U_\hbar(\g)$.
Let $M_\chi^\hbar$ be the generalized Verma module over $U_\hbar(\g)$
corresponding to $\chi$. As before, the morphism $\pi_\chi^\hbar : 
U_\hbar(\g) \to
\End(M_\chi^\hbar)$ restricts to a morphism $\CC[G^*] \to
\End(M_\chi^\hbar)$, which factors though an injective morphism
$\CC[G/L]_\hbar \hookrightarrow \End(M_\chi^\hbar)$ which can be defined
using intertwiners. So $\CC[G/L]_\hbar \simeq \CC[G^*]_\hbar / 
(\CC[G^*]_\hbar \cap \on{Ker}(\pi^\hbar_\chi))$.

\end{document}